\documentclass[11pt]{article}

\usepackage{calc,enumerate}
\usepackage{amscd}
\usepackage{latexsym}
\usepackage{amsmath,amsfonts,amssymb}
\usepackage{diagrams}
  \diagramstyle[balance,righteqno,PostScript=dvips]
  \newarrow{Equalsto}=====
  \newarrow{Mapsto}|--->
   \newarrow{Congruent}33333 
\usepackage{t-angles}
  \providecommand{\obj}[1]{\object{\scriptstyle #1}}
  \providecommand{\Oh}[1]%
  {\step[-0.7]{\hstr{140}\obox1{\scriptstyle #1}}\step[-0.7]}

\numberwithin{equation}{section}
\let\sstyle\scriptstyle

\newcommand{\Label}[1]{\label{#1}\ \\ \makebox[0mm][r]
  {\fontsize{9}{11pt}\fontseries{bx}\fontshape{n}
  \selectfont #1\hspace*{8mm}}\!\!}
 \newcommand{\mathlabel}[1]{\label{#1}\makebox[0mm][r]
  {\fontsize{9}{11pt}\fontseries{bx}\fontshape{n}
  \selectfont #1\hspace*{8mm}}\!\!}
\renewcommand{\mathlabel}{\label}         
\renewcommand{\Label}{\label}             

\newcommand\nowrite[1]{}

\newtheorem{theorem}{{\bf{T}}{\fontsize{7pt}{11}\fontshape{n}
\fontseries{bx}\selectfont{$\!\!\!$HEOREM}}}[section]
\newtheorem{proposition}[theorem]{{\bf{P}}{\fontsize{7pt}{11}\fontshape{n}
\fontseries{bx}\selectfont{$\!\!\!$ROPOSITION}}}
\newtheorem{corollary}[theorem]{{\bf{C}}{\fontsize{7pt}{11}\fontshape{n}
\fontseries{bx}\selectfont{$\!\!\!$OROLLARY}}}
\newtheorem{lemma}[theorem]{{\bf{L}}{\fontsize{7pt}{11}\fontshape{n}
\fontseries{bx}\selectfont{$\!\!\!$EMMA}}}
\newtheorem{definition}[theorem]{\sc Definition}
\newtheorem{remark}{\sc Remark}

\newenvironment{proof}{\par\noindent{\sc Proof.}}{\endproof}
\renewcommand\endproof{\lfl$\sstyle\blacksquare$}

\newenvironment{entry}
 {\begin{list}{}%
  {%
   \setlength{\labelwidth}{12pt}%
   \setlength{\leftmargin}{\labelwidth+\labelsep}%
  }%
 }%
  {\end{list}}

  {\end{list}}

\newcommand\mDelta[8]{
\begin{array}{c}
{\ }\\
\hstretch 150\divide\unitlens by 3%
\begin{tangle}#7\end{tangle}\enspace,\enspace%
\begin{tangle}#8\end{tangle}\\[0pt]\\%
\end{array}}

\def\n-{\nobreakdash-\hspace{0pt}}    

\newcommand\abs{\par\vskip 0.3cm\goodbreak\noindent}

\newcommand\C{\mathcal{C}}

\newcount\siga
\newcount\sigb
\newcount\sigc
\newcount\sigd
\newcount\sige
\newcount\sigf

\newcommand\classbox[6]{\def\Siga{\ }\def\Sigb{\ }\def\Sigc{\ }
 \def\Sigd{\ }\def\Sige{\ }\def\Sigf{\ }
 \siga=#1 \sigb=#2 \sigc=#3 \sigd=#4 \sige=#5 \sigf=#6
 \ifnum\siga=1\def\Siga{\bullet}\fi
 \ifnum\sigb=1\def\Sigb{\bullet}\fi
 \ifnum\sigc=1\def\Sigc{\blacktriangledown}\fi
 \ifnum\sigd=1\def\Sigd{\blacktriangle}\fi
 \ifnum\sige=1\def\Sige{\bullet}\fi
 \ifnum\sigf=1\def\Sigf{\bullet}\fi
 \Put(0,0)[cb]{\begin{tangle}
 \hh\ffbox1{\normalfont\Siga}%
  \ffbox1{\normalfont\Sigb}\\
 \hh\ffbox1{\normalfont\Sigc}\ffbox1{\normalfont\Sigd}\\
 \hh\ffbox1{\normalfont\Sige}\ffbox1{\normalfont\Sigf}
 \end{tangle}}}

\newcommand\classification[2]{\medskip \\ {\em #1:\ } #2}

\newcommand\cpcyalg{\mathop{\bowtie^{\hat\sigma}_{\varphi_{2,1}}}}
\newcommand\cpcyco{\mathop{{}_{\varphi_{1,2}}^{\ \ \,\hat\rho}\!\bowtie}}
\newcommand\cpcybi[1][ ]{\mathop{{}^{\ \ \,\,\hat\rho}_{\varphi_{1,2}}
                {\overset{{}_{#1}}\bowtie}{}_{\varphi_{2,1}}^{\hat\sigma}}}

\newcommand\D{\mathcal{D}}

\newcommand\E{{1\mkern-5mu {\mathrm I}}}
\newcommand\End{\operatorname{End}}

\newcommand\Hom{\operatorname{Hom}}
\newcommand\id{\mathrm{id}}

\newcommand\inj{\mathrm{i}}

\newcommand\ixi[8]{\left\langle\!\!\begin{smallmatrix}%
#3#4\\ #5#6#7#8 \\ #1#2\end{smallmatrix}\!\!\right\rangle}
\newcommand\iy[5]{\left\langle\!\begin{smallmatrix}#5\end{smallmatrix}\!\!%
\begin{smallmatrix}#3#4\\ #1#2\end{smallmatrix}\!\right\rvert}
\newcommand\yi[5]{\left\lvert\!\begin{smallmatrix}#3#4\\ #1#2%
\end{smallmatrix}\!\!\begin{smallmatrix}#5\end{smallmatrix}\!\right\rangle}

\newcommand\lfl{\leaders\hbox to 1em{\hss \hss}\hfill}
\newcommand\m{\mathrm{m}}

\newcommand\nl{\par\noindent}

\newcommand\proj{\mathrm{p}}

\newcounter{Part}
\begin{document}

\parindent=10pt
\author{Yuri Bespalov \and Bernhard Drabant}
\title{\vskip -2truecm{\hfill\normalsize DAMTP-98-117}\\[2truecm]
  {\bfseries\LARGE Cross Product Bialgebras}\\
  {\bfseries\large Part \Roman{Part}}}
\date{July 1998/March 1999}
\maketitle

\bibliographystyle{amsalpha}

\begin{abstract}
\noindent
This is the central article of a series of three papers on cross product
bialgebras. We present a universal theory of bialgebra factorizations
(or cross product bialgebras) with cocycles and dual cocycles. 
We also provide an equivalent (co\n-)modular co-cyclic 
formulation. All known examples as for instance bi-
or smash, doublecross and bicross product bialgebras as well as double
biproduct bialgebras and bicrossed or cocycle bicross product
bialgebras are now united within a single theory. Furthermore our construction
yields various novel types of cross product bialgebras.
\nl
\classification{1991 Mathematics Subject Classification}
               {16W30, 16S40, 18D10, 16B50, 81R50}
\end{abstract}

\section{Introduction}

By definition cross product bialgebras are isomorphic to 
factorizations $B=B_1\otimes B_2$, in a braided category say, such that the
given isomorphisms can be characterized universally in terms of certain
projections and injections from the bialgebra into the particular
tensor factors. In this context tensor product bialgebras, bi- or bismash
product bialgebras \cite{Rad1:85,Tak1:81}, doublecross product bialgebras,
bicross product bialgebras and double biproduct bialgebras
\cite{Ma1:90,Ma1:95} are cross product bialgebras. But also the
bicrossed or cocycle bicross product bialgebras
\cite{MS:94,Ma6:94} are cross product
bialgebras in our terminology; their
universal characterization is given in terms of cleft extensions
and coextensions (see \cite{BCM:86,DT:86,Mon1:92,Sch1:94} for a
comprehensive study of cleft or normal basis Hopf-Galois extensions
on the algebra level). The notion of cross product bialgebras
therefore as well includes constructions involving cocycles and dual
cocycles (or simply ``cycles'').
Well known examples of cross product bialgebras are
Drinfel'd's quantum double, the quantum Poincar\'e group, Radford's
4-parameter Hopf algebra, Lusztig's construction of the quantum
enveloping algebra, the quantum Weyl group, the affine
quantum groups $U_q(\hat g)$, the Connes-Moscovici Hopf algebra in
transverse differential geometry \cite{CM:98}, etc.
The universal constructions of cross product bialgebras in
\cite{Rad1:85,Ma1:90,Ma6:94,MS:94,Tak1:81} additionally admit an equivalent
characterization in terms of mutual weak (co\n-)modular, co-cyclic
relations of the tensor factors $B_1$ and $B_2$. 
This means that there exist certain weak (co\n-)multiplications,
(co\n-)cycles and (co\n-)actions which are morphisms with tensor
rank $(2,1)$ or $(1,2)$ defined on $B_1$ and $B_2$ and their 
two-fold tensor products, such that the multiplication $\m_B$ and the 
comultiplication $\Delta_B$ of the bialgebra $B$ can be composed monoidally
by these structure morphisms. For example we encounter a weak left coaction
$\nu_l:B_1\to B_2\otimes B_1$ which has tensor rank $(1,2)$, 
a weak cocycle $\sigma:B_2\otimes B_2\to B_1$ with rank $(2,1)$,
etc. All these structure morphisms are subject to certain 
relations with rank $(2,2)$, $(3,1)$ and $(1,3)$ which 
on the one hand follow from the bialgebra structure of $B=B_1\otimes B_2$ 
and on the other hand determine the bialgebra structure of $B_1\otimes B_2$ 
uniquely.

Despite those common properties of the cross product bialgebra constructions
in \cite{Rad1:85,Ma1:90,Ma6:94,MS:94,Tak1:81} no theory
exists which characterizes all of them
as special versions of a single universal construction which in turn
is uniquely determined by its (co\n-)modular co-cyclic structure. A first step 
towards this objective has been achieved in
\cite{BD1:98} where we discovered a method to describe cross product bialgebras
without co-cycles, generalizing and uniting \cite{Ma1:90,Ma1:95,Rad1:85}. 

The present article considerably extends the results of \cite{BD1:98}.
The main outcome is a universal and (co\n-)modular, co-cyclic
theory of cross product bialgebras with cocycles and cycles which unites
all known constructions within a single setting
and provides several new families of cross product
bialgebras. This is the most comprising
framework for cross product bialgebras so far, which equivalently takes into
account both universal and (co\n-)modular co-cyclic aspects.
Furthermore our construction is designed to work in arbitrary braided
categories\footnote{Therefore we have been able to apply the results
of \cite{BD1:98} to categories of Hopf bimodules and
Yetter-Drinfel'd modules \cite{BD3:98}; we characterized two-sided bi-
or bismash product bialgebras as certain cross product Hopf bimodule
bialgebras. In particular we showed that the double biproduct is a twisted
Hopf bimodule tensor product bialgebra.}.

To derive these results we will proceed in two steps.
We begin with the general definition of cross product bialgebras
and find a universal characterization for them.
We show that so-called cocycle cross product bialgebras
possess a certain (co\n-)modular co\n-cyclic structure.
Then we restrict our consideration to strong
cross product bialgebras. They also admit a universal
characterization and will turn out to be the central objects in our
theory of cross product bialgebras.

In a second step we present a (co\n-)modular co-cyclic
construction method by so-called Hopf data\footnote{Observe that
Hopf data without co-cycles have been introduced in \cite{BD1:98} already.
There should be no confusion with the more general notation of the same 
name which will be defined in the present article. In a very similar
context the term ``Hopf datum" has been used in \cite{And1:96} for
special cases of bicrossed product bialgebras \cite{Ma6:94}.}.
A Hopf datum
consists of a pair of objects $(B_1,B_2)$ and weak (co\n-)\-ac\-tions,
(co\n-)\-mul\-ti\-pli\-cations and co-cycles defined on 
$B_1$ and $B_2$ which obey certain 
interrelated identities with rank $(2,2)$, $(3,1)$ and $(1,3)$. 
Then we introduce strong Hopf data for which additional
``strong'' relations hold. We show that strong Hopf data induce the
structure of a strong cross product bialgebra on the tensor
product $B_1\otimes B_2$. Strong Hopf data and strong cross
product bialgebras are different depictions of the same object.

Then we combine the universal characterization and the (co\n-)modular 
co-cyclic description in terms of strong Hopf data to
obtain the theory of strong cross product bialgebras.
This is the central result of the article. Eventually we apply our 
construction and investigate strong cross product bialgebras according 
to their (co\n-)modular co-cyclic structure. In particular we recover all known 
constructions \cite{Rad1:85,Ma1:90,MS:94,Ma6:94,BD1:98} and
find various new types of cross product bialgebras. All of them are
special versions of the most general strong cross
product bialgebra construction.
\abs
After introducing preliminary notations and definitions we study
cross product algebras and cross product coalgebras from a
universal point of view in Section \ref{sec-ccpa}. These results will be used
in Section \ref{sec-ccpb} to define cross product bialgebras. From
this general definition we extract cocycle cross product bialgebras
and strong cross product bialgebras. The structure of
cocycle cross product bialgebras
gives useful hints for the design of Hopf data. This will be done in
Section \ref{hopf-data} where we
investigate the connection of Hopf data and cocycle cross product bialgebras.
A Hopf datum is canonically assigned to every cocycle cross
product bialgebra. The converse is not true in general.
But Theorem \ref{hp-cp} states that strong Hopf data yield strong
cross product bialgebras. More precisely strong Hopf data
and strong cross product bialgebras are
equivalent constructions which are in one-to-one
correspondence. The universal (co\n-)modular co-cyclic theory of
strong cross product bialgebras culminates in the
central Theorems \ref{hp-cp} and \ref{central-theory}.
Finally we present a table of special strong 
cross product bialgebras which contains all known and several
new types of cross product bialgebras embedded in our most general framework.

The rather intricate proof of Theorem \ref{hp-cp} has been postponed to
Section \ref{proof-hp-cp}.
\abs
In the recent articles \cite{CIMZ:98,CZ:98,Sbg1:98}
cross product bialgebras have been studied as well. 
In \cite{CIMZ:98} the universal properties of co-cycle free cross 
product bialgebras have been considered. The article 
\cite{CZ:98} investigates certain cross product
bialgebras without co-cycles in universal and (co\n-)modular co-cyclic terms.
In \cite{Sbg1:98} a special type of our cocycle cross product
bialgebras with one cocycle has been studied in the special case where
all objects are vector spaces over a field.

\subsection{Preliminaries}

We are working throughout in (strict) braided monoidal categories
\cite{JS:93}. In our article we denote
categories by calligraphic letters $\C$, $\D$, etc. For a braided monoidal
category $\C$ the tensor product is denoted by $\otimes_\C$, the unit
object by $\E_\C$, and the braiding by ${}^\C\Psi$. If it is clear from
the context we omit the index `$\C$' at the various symbols.
We confine ourselves to braided categories which admit split idempotents
\cite{BD:95,Lyu1:95}; for each idempotent $\Pi=\Pi^2:M\to M$ of any
object $M$ in $\C$ there exists an object $M_\Pi$ and a pair of morphisms
$(\inj_\Pi,\proj_\Pi)$ such that $\proj_\Pi\circ\inj_\Pi=\id_{M_\Pi}$
and $\inj_\Pi\circ\proj_\Pi=\Pi$. This is not a severe restriction
of the categories since every braided category
can be canonically embedded into a braided category which admits split
idempotents \cite{BD:95,Lyu1:95}.

We use generalized algebraic
structures like algebra, bialgebra, module, comodule, etc.~in such
categories. We assume the reader is familiar with these
generalizations. A thorough introduction to braided algebraic
structures and the graphical calculus coming along with it
can be found in \cite{FY:89,JS:93,Kas1:95,Lyu1:95,Ma4:93,ReT:90,Tur1:94}.
We denote by $\m:A\otimes A\to A$ the multiplication and by
$\eta:\E\to A$ the unit of an algebra $A$ in $\C$.
$\Delta: C\to C\otimes C$ is meant for the comultiplication and
$\varepsilon:C\to \E$ for the counit of a coalgebra $C$ in $\C$,
$\mu_l:A\otimes M\to M$ is the left action of an algebra $A$ on a
module $M$, and $\nu_l:N\to C\otimes N$ denotes the left coaction of a
coalgebra $C$ on a comodule $N$. Right actions are denoted by $\mu_r$ and
right coactions by $\nu_r$. In the course of our work we often
encounter weak versions of (co\n-)multiplication and (co\n-)action
which do not necessarily obey the relations of `true' (co\n-)algebras
and (co\n-)modules. We will nevertheless use the same symbols as
for proper (co\n-)multiplications and (co\n-)actions -- it should
become clear from the context whether these structure morphisms are weak or
not.

In our article we make use of graphical calculus for (strict) braided monoidal
categories which simplifies intricate categorical equations of
morphisms and helps to uncover their intrinsic structure.
Morphisms will be composed from up to down, i.~e.~the domains of the morphisms
are at the top and the codomaines are at the bottom of the graphics.
Tensor products are represented by horizontal concatenation
in the corresponding order.
We present our own conventions \cite{BD:95,BD:97,BD1:98} in Figure
\ref{fig-conv}. If there is no fear of confusion
we omit the assignment of a specific object to the ends
of the particular strings in the graphics.
\begin{figure}
\begin{equation*}
\begin{array}{c}
\m=\divide\unitlens by 2
   \begin{tangle}\cu\end{tangle}
   \multiply\unitlens by 2
\quad ,\quad
\eta=\divide\unitlens by 2
   \ \begin{tangle}\unit\end{tangle}
   \multiply\unitlens by 2
\quad ,\quad
\Delta=\divide\unitlens by 2
   \begin{tangle}\cd\end{tangle}
   \multiply\unitlens by 2
\quad ,\quad
\varepsilon=\divide\unitlens by 2
   \ \begin{tangle}\counit\end{tangle}
   \multiply\unitlens by 2
\quad ,\quad
\\[8pt]
\mu_l=\divide\unitlens by 2
   \begin{tangle}\lu\end{tangle}
   \multiply\unitlens by 2
\quad ,\quad
\mu_r=\divide\unitlens by 2
   \begin{tangle}\ru\end{tangle}
   \multiply\unitlens by 2
\quad ,\quad
\nu_l=\divide\unitlens by 2
   \begin{tangle}\ld\end{tangle}
   \multiply\unitlens by 2
\quad ,\quad
\nu_r=\divide\unitlens by 2
   \begin{tangle}\rd\end{tangle}
   \multiply\unitlens by 2
\quad ,\quad
\\[8pt]
\Psi=\divide\unitlens by 2
   \begin{tangle}\x\end{tangle}
   \multiply\unitlens by 2
\quad ,\quad
\Psi^{-1}=\divide\unitlens by 2
   \begin{tangle}\xx\end{tangle}
   \multiply\unitlens by 2
\quad .\quad
\end{array}
\end{equation*}
\caption{Graphical presentation of (weak) multiplication $\m$, unit
$\eta$, (weak) comultiplication $\Delta$, counit $\varepsilon$,
(weak) left action $\mu_l$, (weak) right action $\mu_r$,
(weak) left coaction $\nu_l$, (weak) right coaction $\nu_r$,
braiding $\Psi$, and inverse braiding
$\Psi^{-1}$.}\Label{fig-conv}
\end{figure}
Below we elucidate the graphical calculus to readers unfamiliar
with this technique. The first example is the associativity of the
action of an algebra $A$ on a right $A$-module $M$. The second is the
naturality of the braiding in a braided category.
\begin{equation*}
\begin{array}{ccc}
\mu_r\circ(\mu_r\otimes\id_A)=\mu_r\circ(\id_M\otimes\m_A)
&\quad\text{corresponds to}\quad
&\divide\unitlens by 2
{\begin{tangle}
\hh\hru\hstep\id\\
\ru
\end{tangle}}
=
{\begin{tangle}
\hh\id\hstep\hcu\\
\ru
\end{tangle}}\,\,\,,
\\[10pt]
\Psi\circ (f\otimes g) = (g\otimes f)\circ\Psi
&\quad\text{corresponds to}\quad
&\divide\unitlens by 2
{\begin{tangle}
\O{\sstyle f}\step[2]\O{\sstyle g}\\
\x
\end{tangle}}
\,=\,
{\begin{tangle}
\x\\
\O{\sstyle g}\step[2]\O{\sstyle f}
\end{tangle}}\,\,\,\,\,.
\end{array}
\end{equation*}
Note that throughout the article there is (almost) no need to require
invertibility of the braiding. Therefore most of the results can be 
derived if we assume that the underlying category $\C$ is pre-braided.
We do not discuss these generalizations further and confine to braided
categories in what follows.

\section{Cross Product Algebras and Coalgebras}\Label{sec-ccpa}

In the first part of Section \ref{sec-ccpa} we study cross
product algebras. They have been considered also in \cite{Brz1:96}.
It turns out that cross product algebras are universal
constructions. They generalize crossed product algebras
\cite{DT:86,Mon1:92,Sch1:94}.
Cross product coalgebras will be studied
in the second part of Section \ref{sec-ccpa}. Since the results for
cross product coalgebras can be obtained easily by certain
categorical dualization, we omit all the proofs in this case and
refer to the analogous proofs for cross product algebras.
Both structures, cross product algebras and cross product
coalgebras, will be needed later in the definition of cross product
bialgebras.

\subsection{Cross Product Algebras}

\begin{definition}\Label{caat}
Let $(B_1,\m_1,\eta_1)$ be an algebra and $B_2$ be an object in $\C$.
Suppose there are morphisms $\eta_2:\E\to B_2$,
$\varphi_{2,1}:B_2\otimes B_1\to B_1\otimes B_2$, and
$\hat\sigma:B_2\otimes B_2\to B_1\otimes B_2$ such that
the relation
\begin{equation}\mathlabel{cond1-caat}
\varphi_{2,1}=(\m_1\otimes\id_{B_2})\circ(\id_{B_1}\otimes\hat\sigma)\circ
(\varphi_{2,1}\otimes\eta_2)
\end{equation}
holds. Suppose further that $B=B_1\otimes B_2$ is an algebra through
\begin{equation}\mathlabel{cond2-caat}
\m_B = (\m_1\otimes\id_{B_2})\circ(\m_1\otimes\hat\sigma)\circ
(\id_{B_1}\otimes\varphi_{2,1}\otimes\id_{B_2})\,,\quad
\eta_B =\eta_1\otimes\eta_2
\end{equation}
then $B$ is called cross product algebra and will be denoted by
$B_1\cpcyalg B_2$.
\end{definition}
\abs
In the subsequent proposition we find equivalent constructive
conditions for cross product algebras. Similar results
have been obtained in \cite{Brz1:96}.

\begin{proposition}\Label{cross-equiv}
The following statements are equivalent.
\begin{enumerate}
\item
$B_1\cpcyalg B_2$ is a cross product algebra.
\item
$(B_1,\m_1,\eta_1)$ is an algebra and $B_2$ is an object in $\C$. There
exist morphisms $\eta_2:\E\to B_2$, $\varphi_{2,1}:B_2\otimes B_1\to
B_1\otimes B_2$ and $\hat\sigma:B_2\otimes B_2\to B_1\otimes B_2$ for
which the subsequent identities hold.
\begin{equation}\mathlabel{cond3-caat}
\begin{gathered}
{\begin{split}
\varphi_{2,1}\circ(\eta_2\otimes\id_{B_1})&=\id_{B_1}\otimes\eta_2\,,\\
\varphi_{2,1}\circ(\id_{B_2}\otimes\eta_1)&=\eta_1\otimes\id_{B_2}\,,
\end{split}
\quad
\begin{split}
\hat\sigma\circ(\eta_2\otimes\id_{B_2})&=\eta_1\otimes\id_{B_2}\,,\\
\hat\sigma\circ(\id_{B_2}\otimes\eta_2)&=\eta_1\otimes\id_{B_2}\,,
\end{split}}\\[5pt]
{\begin{gathered}
\varphi_{2,1}\circ(\id_{B_2}\otimes\m_1)=
(\m_1\otimes\id_{B_2})\circ(\id_{B_1}\otimes\varphi_{2,1})\circ
(\varphi_{2,1}\otimes\id_{B_1})\,,\\[5pt]
(\m_1\otimes\id_{B_2})\circ(\id_{B_1}\otimes\hat\sigma)\circ
(\varphi_{2,1}\otimes\id_{B_2})\circ(\id_{B_2}\otimes\hat\sigma)\\
=(\m_1\otimes\id_{B_2})\circ(\id_{B_1}\otimes\hat\sigma)\circ
(\hat\sigma\otimes\id_{B_2})\,,\\[5pt]
(\m_1\otimes\id_{B_2})\circ(\id_{B_1}\otimes\hat\sigma)\circ
(\varphi_{2,1}\otimes\id_{B_2})\circ(\id_{B_2}\otimes\varphi_{2,1})\\
=(\m_1\otimes\id_{B_2})\circ(\id_{B_1}\otimes\varphi_{2,1})\circ
(\hat\sigma\otimes\id_{B_1})\,.
\end{gathered}}
\end{gathered}
\end{equation}
\end{enumerate}
\end{proposition}

\begin{proof}
Suppose that $B_1\cpcyalg B_2$ is a cross product algebra.
Since by assumption $(B_1,\m_1,\eta_1)$ is itself an algebra one
obtains the first and second identity of \eqref{cond3-caat}
by unitality of $\m_B$ and $\m_1$. Using these relations
yields the fourth identity of \eqref{cond3-caat}
by application of $\id_{B_2}\otimes\eta_1$ to
\eqref{cond1-caat}. Again using the left unitality of $\m_B$ one concludes
that $(\m_1\otimes\id_{B_2})\circ(\id_{B_1}\otimes
\hat\sigma\circ(\eta_2\otimes\id_{B_2}))=\id_{B_1}\otimes\id_{B_2}$
from which one immediately derives the third identity of
\eqref{cond3-caat}. Now the associativity
$\m_B\circ(\m_B\otimes\id_B)= \m_B\circ(\id_B\otimes\m_B)$ yields
the fifth identity of \eqref{cond3-caat} by application of
$\eta_1\otimes\id\otimes\id\otimes\eta_1\otimes\id\otimes\eta_2$.
The sixth relation is obtained by applying
$\eta_1\otimes\id\otimes\eta_1\otimes\id\otimes\eta_1\otimes\id$.
To derive the seventh identity one has to apply
$\eta_1\otimes\id\otimes\eta_1\otimes\id\otimes\id\otimes\eta_2$.
Conversely suppose the conditions of the second item of the proposition
are fulfilled. Define $\m_B$ according to \eqref{cond2-caat}. Then
${\m_B\circ(\m_B\otimes\id_B)} = (\m_1^{(5)}\otimes\id)\circ
 (\id^{(4)}\otimes\hat\sigma)\circ(\id^{(3)}\otimes\varphi_{2,1}\otimes\id)
 \circ(\id^{(2)}\otimes\hat\sigma\otimes\id^{(2)})\circ
 (\id\otimes\varphi_{2,1}\otimes\id^{(3)})$
where $\m_1^{(n)}:B_1^{\otimes_n}
\to B_1$ is the canonical $n$-fold multiplication, and
$\id^{(n)}$ is the abbreviation of the identity of an $n$-fold tensor
product of (combined) $B_1$'s and $B_2$'s.
On the other hand
\begin{equation*}
\begin{split}
\m_B\circ(\id_B\otimes\m_B) &= (\m_1^{(4)}\otimes\id)\circ
 (\id^{(4)}\otimes\hat\sigma)\circ(\id^{(3)}\otimes\varphi_{2,1}\otimes\id)
 \circ\\
&\quad\circ(\id^{(2)}\otimes\varphi_{2,1}\otimes\hat\sigma)\circ
 (\id\otimes\varphi_{2,1}\otimes\varphi_{2,1}\otimes\id)\\
&=(\m_1^{(5)}\otimes\id)\circ(\id^{(4)}\otimes\hat\sigma)\circ
 (\id^{(3)}\otimes\hat\sigma\otimes\id)\circ\\
&\quad\circ(\id^{(2)}\otimes\varphi_{2,1}\otimes\id^{(2)})\circ
 (\id\otimes\varphi_{2,1}\otimes\varphi_{2,1}\otimes\id)\\
&= (\m_1^{(5)}\otimes\id)\circ
 (\id^{(4)}\otimes\hat\sigma)\circ(\id^{(3)}\otimes\varphi_{2,1}\otimes\id)
 \circ\\
&\quad\circ(\id^{(2)}\otimes\hat\sigma\otimes\id^{(2)})\circ
 (\id\otimes\varphi_{2,1}\otimes\id^{(3)})
\end{split}
\end{equation*}
where the fourth condition of \eqref{cond3-caat} has been used twice
to obtain the first equation, the sixth relation of \eqref{cond3-caat}
has been applied in the second equation, and with the help of the seventh
relation of \eqref{cond3-caat} we obtained the third equation.
Hence associativity of $\m_B$ has been proven. Using the first
four identities of \eqref{cond3-caat} one easily proves
\eqref{cond1-caat} and unitality
$\m_B\circ(\eta_B\otimes\id_B)=\id_B=\m_B\circ(\id_B\otimes\eta_B)$.
\end{proof}
\abs
\begin{remark} {\normalfont Condition \eqref{cond1-caat}
in Definition \ref{caat} can be replaced equivalently by the identity
$(\id_{B_1}\otimes\hat\sigma)
=(\m_1\otimes\id_{B_2})\circ(\id_{B_1}\otimes\hat\sigma)\circ
(\varphi_{2,1}\circ(\id_{B_2}\otimes\eta_1)\otimes\id_{B_2})$.}
\end{remark}

\begin{remark} {\normalfont There is another equivalent definition of 
cross product algebras which has been kindly reported to us by
Gigel Militaru. Observe that the morphism $\Gamma$ below is our morphism
$\m_{0,20}$ in \eqref{aux-m-delta}. Then the following statements are 
equivalent.
\begin{enumerate}
\item
$B_1\cpcyalg B_2$ is a cross product algebra.
\item
$(B_1,\m_1,\eta_1)$ is an algebra and $B_2$ is an object in $\C$.
There are morphisms $\eta_2:\E\to B_2$ and
$\Gamma:B_2\otimes B_1\otimes B_2\to B_1\otimes B_2$ such that
$B=B_1\otimes B_2$ is an algebra through
\begin{equation}\mathlabel{cond4-caat}
\m_B = (\m_1\otimes\id_{B_2})\circ(\id_{B_1}\otimes\Gamma)\quad\text{and}
\quad\eta_B =\eta_1\otimes\eta_2\,.
\end{equation}
\item
$(B_1,\m_1,\eta_1)$ is an algebra and $B_2$ is an object in $\C$. There
exist morphisms $\eta_2:\E\to B_2$ and
$\Gamma:B_2\otimes B_1\otimes B_2\to B_1\otimes B_2$ such that
\begin{equation}\mathlabel{cond5-caat}
\begin{gathered}
\Gamma\circ(\eta_2\otimes\id_{B_1\otimes B_2})=\id_{B_1\otimes B_2}\,,
\quad
\Gamma\circ(\id_{B_2}\otimes\eta_1\otimes\eta_2)=
\eta_1\otimes\id_{B_2}\,,\\[5pt]
(\m_1\otimes\id_{B_2})\circ(\id_{B_1}\otimes\Gamma)\circ
(\Gamma\otimes\id_{B_1\otimes B_2})\\
=\Gamma\circ(\id_{B_2}\otimes(\m_1\otimes\id_{B_2})\circ
(\id_{B_1}\otimes\Gamma))\,.
\end{gathered}
\end{equation}
\end{enumerate}
The one-to-one correspondence is given by 
$\Gamma= (\m_1\otimes\id_{B_2})\circ(\id_{B_1}\otimes\hat\sigma)\circ
(\varphi_{2,1}\otimes\id_{B_2})$ and its inverse
$\varphi_{2,1}=\Gamma\circ(\id_{B_2\otimes B_1}\otimes\eta_2)$ and
$\hat\sigma= \Gamma\circ(\id_{B_2}\otimes\eta_1\otimes\id_{B_2})$.}
\end{remark}
\abs
We will show that cross product algebras are universal
constructions. The first proposition describes equivalent projection and
injection conditions for algebras isomorphic to cross product
algebras (see also \cite{BD1:98} for the cocycle-free case).
The second proposition is closely related to the first one and
characterizes the universal construction explicitely.

\begin{proposition}\Label{univ-constr}
Let $A$ be an algebra in $\C$. Then it holds equivalently
\begin{enumerate}
\item
$A$ is algebra isomorphic to a cross product algebra
$B_1\cpcyalg B_2$.
\item
There is an algebra $(B_1,\m_1,\eta_1)$, an object $B_2$,
morphisms $B_1\overset{\inj_1}\to A\overset{\inj_2}\leftarrow B_2$ and
$\eta_2:\E\to B_2$ where $\inj_1$ is an algebra morphism and
$\inj_2\circ\eta_2=\eta_A$ such that
$\m_A\circ(\inj_1\otimes\inj_2): B_1\otimes B_2\to A$
is an isomorphism in $\C$.
\end{enumerate}
\end{proposition}

\begin{proof}
If $\Lambda:B_1\cpcyalg B_2\to A$ is an algebra isomorphism then
define $\inj_1:=\Lambda\circ(\id_{B_1}\otimes\eta_2)$
and $\inj_2:=\Lambda\circ(\eta_1\otimes\id_{B_2})$.
Using the particular definition of $\m_{B_1\cpcyalg B_2}$ in
\eqref{cond2-caat} and the identities \eqref{cond1-caat},
it is verified immediately that $\inj_1$ is an algebra morphism
since $\Lambda$ is an algebra morphism. Similarly one proves that
$\inj_2\circ\eta_2=\eta_A$ and $\m_A\circ(\inj_1\otimes\inj_2)=\Lambda$
which is therefore an isomorphism.
If on the other hand the conditions of the second statement of
Proposition \ref{univ-constr} is fulfilled, then
$\Lambda:=\m_A\otimes(\inj_1\otimes\inj_2$) is an isomorphism by
assumption. Therefore $B=B_1\otimes B_2$ is canonically an algebra
through $\m_B:=\Lambda^{-1}\circ\m_A\circ(\Lambda\otimes\Lambda)$ and
$\eta_B:=\Lambda^{-1}\circ\eta_A$. We will show that this defines a
cross product algebra structure on $B$ by
$\varphi_{2,1}:=\Lambda^{-1}\circ\m_A\circ(\inj_2\otimes\inj_1)$
and $\hat\sigma:=\Lambda^{-1}\circ\m_A\circ(\inj_2\otimes\inj_2)$.
Using the explicit expression for $\Lambda$,
inserting several times $\Lambda\circ\Lambda^{-1}$ and using the above
definitions, as well as the fact that $\inj_1$ is
algebra morphism, one obtains the following identities.
\begin{equation*}
\begin{split}
\m_B &=\Lambda^{-1}\circ\m_A^{(3)}\circ
       (\id\otimes\Lambda\circ\Lambda^{-1}\otimes\id)\circ
       (\inj_1\otimes\m_A\circ(\inj_2\otimes\inj_1)\otimes\inj_2)\\
   &=\Lambda^{-1}\circ\m_A^{(4)}\circ
     (\inj_1\otimes(\inj_1\otimes\inj_2)\circ\varphi_{2,1}\otimes\inj_2)\\
   &=\Lambda^{-1}\circ\m_A^{(3)}\circ(\inj_1\otimes\inj_1\otimes\inj_2)\circ
     (\m_1\otimes\hat\sigma)\circ(\id\otimes\varphi_{2,1}\otimes\id)\\
   &=(\m_1\otimes\id)\circ(\m_1\otimes\hat\sigma)\circ
     (\id\otimes\varphi_{2,1}\otimes\id)\,.
\end{split}
\end{equation*}
Similarly one derives $\Lambda\circ(\eta_1\otimes\eta_2)=\eta_A$, hence
$\eta_B=\eta_1\otimes\eta_2$. Therefore
$(B,\m_B,\eta_B)=B_1\cpcyalg B_2$.\end{proof}
\abs
\begin{proposition}\Label{univ-constr1}
Let $B_1\cpcyalg B_2$ be a cross product algebra, and $A$ be
an algebra. Suppose there are morphisms $\alpha:B_1\to A$ and
$\beta:B_2\to A$ such that
\begin{enumerate}
\item
$\alpha$ is an algebra morphism.
\item
$\beta\circ\eta_2=\eta_A\,.$
\item
$\m_A\circ(\alpha\otimes\beta)\circ\varphi_{2,1}
 = \m_A\circ(\beta\otimes\alpha)\,.$
\item
$\m_A\circ(\alpha\otimes\beta)\circ\hat\sigma
 = \m_A\circ(\beta\otimes\beta)\,.$
\end{enumerate}
Then there exists a unique algebra morphism $\gamma: B_1\cpcyalg B_2\to A$
obeying the identities $\gamma\circ(\id_{B_1}\otimes\eta_2)=\alpha$ and
$\gamma\circ(\eta_1\otimes\id_{B_2})=\beta$.
\end{proposition}

\begin{proof}
Define $\gamma:=\m_A\circ(\alpha\otimes\beta)$. Then by assumption 1
and 2 of Proposition \ref{univ-constr1} it follows
$\gamma\circ(\id_{B_1}\otimes\eta_2)=\alpha$,
$\gamma\circ(\eta_1\otimes\id_{B_2})=\beta$, and
$\gamma\circ(\eta_1\otimes\eta_2)=\eta_A$.
Using consecutively that $\alpha$ is an algebra morphism,
assumption 4, and assumption 3 of Proposition \ref{univ-constr1} then
yields $\gamma\circ\m_{B_1\cpcyalg B_2} = \m_A^{(4)}\circ
(\alpha\otimes\alpha\otimes(\alpha\otimes\beta)\circ\hat\sigma)\circ
(\id\otimes\varphi_{2,1}\otimes\id)=\m_A^{(4)}\circ
(\alpha\otimes(\alpha\otimes\beta)\circ\varphi_{2,1}\otimes\beta)
=\m_A\circ(\m_A\circ(\alpha\otimes\beta)\otimes\m_A\circ(\alpha\otimes\beta))
=\m_A\circ(\gamma\otimes\gamma)$.
Therefore $\gamma$ is an algebra morphism obeying the conditions
of Proposition \ref{univ-constr1}. Suppose that there is another
$\gamma'$ meeting the same stipulations as $\gamma$. Then
$\gamma=\m_A\circ(\alpha\otimes\beta)=\m_A\circ
(\gamma'\circ(\id\otimes\eta_2)\otimes\gamma'\circ(\eta_1\otimes\id))
=\gamma'\circ\m_{B_1\cpcyalg B_2}\circ(\id\otimes\eta_2\otimes
\eta_1\otimes\id)=\gamma'$.
In the last equation unital properties of cross product algebra
have been used. This proves uniqueness of $\gamma$.\end{proof}
\abs
\begin{remark}{\normalfont The conditions of Proposition \ref{univ-constr1}
can be realized on $A=B_1\cpcyalg B_2$ with $\alpha=\id_{B_1}\otimes\eta_2$
and $\beta=\eta_1\otimes\id_{B_2}$. Then of course $\gamma=\id$.}
\end{remark}
\abs
In the following we consider generalized smash product algebras \cite{BD1:98}
and demonstrate that they are special cases of cross product
algebras\footnote{We adopt the naming "smash product algebras" from
\cite{CIMZ:98}.}. A (generalized) {\it smash product algebra}
$B_1\times_{\phi_{2,1}}B_2$ consists of
algebras $B_1$ and $B_2$, and a morphism
$\phi_{2,1}:B_2\otimes B_1\to B_1\otimes B_2$,
such that $B=B_1\otimes B_2$ is an algebra through
$\m_B=(\m_1\otimes\m_2)\circ(\id_{B_1}\otimes\phi_{2,1}\otimes\id_{B_2})$
and $\eta_B=\eta_1\otimes\eta_2$.
The next proposition shows that smash product algebras are indeed special
cases of cross product algebras under certain natural conditions.
\abs
\begin{proposition}\Label{smash-ccpa}
Suppose there is a morphism $\varepsilon_1:B_1\to \E$ with
$\varepsilon_1\circ\eta_1=\id_{\E}$. Then
\begin{entry}
\item[\ ]
$B_1\cpcyalg B_2$ is a cross product algebra
with $\hat\sigma = \eta_1\otimes\m_2$ for some morphism
$\m_2:B_2\otimes B_2\to B_2$.
\item[$\Leftrightarrow$]
$B_1\times_{\varphi_{2,1}}B_2$ is a smash product algebra.
\end{entry}
\end{proposition}

\begin{proof}
The proposition can be proven easily by
\eqref{cond3-caat} and triviality of $\hat\sigma$.\end{proof}
\abs
\begin{remark}{\normalfont Crossed product algebras
$A\#_\sigma H$ \cite{DT:86,Mon1:92}
are special examples of cross product algebras through
$\varphi_{2,1}:=(\mu_l\otimes\id_H)\circ(\Delta_H\otimes\id_A)$ and
$\hat\sigma :=(\sigma\otimes\m_H)\circ(\id_H\otimes\Psi_{H,H}\otimes\id_H)
\circ(\Delta_H\otimes\Delta_H)$
where $A$ is an algebra, $\mu_l:H\otimes A\to A$ is a left $H$-measure
on $A$, and $\sigma:H\otimes H\to A$ is a (convolution invertible)
cocycle.}
\end{remark}
\abs
The next proposition yields criteria under which a morphism
$f:B_1\cpcyalg B_1\to A$ is multiplicative.

\begin{proposition}\Label{alg-morph-equiv}
Let $A$ be an algebra and $B=B_1\cpcyalg B_2$ be a cross product
algebra. Then the subsequent statements are equivalent.
\begin{enumerate}
\item\Label{alg-morph-equiv1}
The morphism $f:B\to A$ is multiplicative,
i.~e.~$f\circ\m_B=\m_A\circ(f\otimes f)$.
\item\Label{alg-morph-equiv2}
The identities
\begin{equation*}
\begin{split}
\m_A\circ(f\otimes f)\circ
(\id_{B_1}\otimes\eta_2\otimes\id_{B_1}\otimes\id_{B_2}) &=
f\circ\m_B\circ(\id_{B_1}\otimes\eta_2\otimes\id_{B_1}\otimes\id_{B_2})\,,\\
\m_A\circ(f\otimes f)\circ
(\eta_1\otimes\id_{B_2}\otimes\id_{B_1}\otimes\id_{B_2}) &=
f\circ\m_B\circ(\eta_1\otimes\id_{B_2}\otimes\id_{B_1}\otimes\id_{B_2})\\[-12pt]
\end{split}
\end{equation*}
are satisfied.
\item\Label{alg-morph-equiv3}
The identities
\begin{equation*}
\begin{split}
\m_A\circ(f\otimes f)\circ
(\id_{B_1}\otimes\id_{B_2}\otimes\id_{B_1}\otimes\eta_2) &=
f\circ\m_B\circ(\id_{B_1}\otimes\id_{B_2}\otimes\id_{B_1}\otimes\eta_2)\,,\\
\m_A\circ(f\otimes f)\circ
(\id_{B_1}\otimes\id_{B_2}\otimes\eta_1\otimes\id_{B_2}) &=
f\circ\m_B\circ(\id_{B_1}\otimes\id_{B_2}\otimes\eta_1\otimes\id_{B_2})\\[-12pt]
\end{split}
\end{equation*}
are satisfied.
\item\Label{alg-morph-equiv4}
The identities
\begin{equation*}
\begin{split}
\m_A\circ(f\otimes f)\circ
(\id_{B_1}\otimes\eta_2\otimes\id_{B_1}\otimes\id_{B_2}) &=
f\circ\m_B\circ(\id_{B_1}\otimes\eta_2\otimes\id_{B_1}\otimes\id_{B_2})\,,\\
\m_A\circ(f\otimes f)\circ
(\id_{B_1}\otimes\id_{B_2}\otimes\eta_1\otimes\id_{B_2}) &=
f\circ\m_B\circ(\id_{B_1}\otimes\id_{B_2}\otimes\eta_1\otimes\id_{B_2})\,,\\
\m_A\circ(f\otimes f)\circ
(\eta_1\otimes\id_{B_2}\otimes\id_{B_1}\otimes\eta_2) &=
f\circ\m_B\circ(\eta_1\otimes\id_{B_2}\otimes\id_{B_1}\otimes\eta_2)\\[-12pt]
\end{split}
\end{equation*}
are satisfied.
\end{enumerate}
\end{proposition}

\begin{proof}
The non-trivial part of the proposition can be verified easily with the help
of the identity
$\m_B\circ(\id\otimes\eta_2\otimes\eta_1\otimes\id)=\id\otimes\id$,
the associativity of $A$ and $B$, and the assumptions
\ref{alg-morph-equiv}.\ref{alg-morph-equiv2},
\ref{alg-morph-equiv}.\ref{alg-morph-equiv3}
or \ref{alg-morph-equiv}.\ref{alg-morph-equiv4} respectively.
\end{proof}

\subsection{Cross Product Coalgebras}

Cross product coalgebras are somehow dual constructions to
cross product algebras. We will omit proofs since they can be
obtained from the corresponding proofs for cross product
algebras by an obvious kind of dualization.
\abs
\begin{definition}\Label{ccat}
Let $(C_2,\Delta_2,\varepsilon_2)$ be a coalgebra and $C_1$ be an
object in $\C$. Suppose there exist morphisms $\varepsilon_1:B_1\to\E$,
$\varphi_{1,2}:C_1\otimes C_2\to C_2\otimes C_1$,
$\hat\rho:C_1\otimes C_2\to C_1\otimes C_1$ such that the relation
$\varphi_{1,2}=(\varepsilon_1\otimes\varphi_{1,2})\circ
(\hat\rho\otimes\id_{C_2})\circ(\id_{C_1}\otimes\Delta_2)$
holds. Then $C=C_1\otimes C_2$ is called cross product coalgebra if
it is a coalgebra through $\Delta_C= (\id_{B_1}\otimes
\varphi_{1,2}\otimes\id_{B_2})\circ(\hat\rho\otimes\Delta_2)\circ
(\id_{B_1}\otimes\Delta_2)$ and $\varepsilon_C=
\varepsilon_1\otimes\varepsilon_2$. We denote $C$ by $C_1\cpcyco C_2$.
\end{definition}
\abs
\begin{proposition}
The following statements are equivalent.
\begin{enumerate}
\item
$C_1\cpcyco C_2$ is a cross product coalgebra.
\item
$(C_2,\Delta_2,\varepsilon_2)$ is a coalgebra and $C_1$ is an object in
$\C$. There are morphisms $\varepsilon_1:B_1\to\E$ and
$\varphi_{1,2} : C_1\otimes C_2\to C_2\otimes C_1$ such that
\begin{equation*}
\begin{gathered}
{\begin{split}
(\id_{C_2}\otimes\varepsilon_1)\circ\varphi_{1,2}&=
\varepsilon_1\otimes\id_{C_2}\,,\\
(\varepsilon_1\otimes\id_{C_2})\circ\varphi_{1,2}&=\id_{C_2}\otimes
\varepsilon_1\,,
\end{split}
\quad
\begin{split}
(\id_{C_1}\otimes\varepsilon_1)\circ\hat\rho &=\id_{C_1}\otimes
\varepsilon_2\,,\\
(\varepsilon_1\otimes\id_{C_1})\circ\hat\rho &=
\id_{C_1}\otimes\varepsilon_2\,,
\end{split}}\\[5pt]
{\begin{gathered}
(\Delta_2\otimes\id_{C_1})\circ\varphi_{1,2}=
(\id_{C_2}\otimes\varphi_{1,2})\circ(\varphi_{1,2}\otimes\id_{C_2})\circ
(\id_{C_1}\otimes\Delta_2)\,,\\[5pt]
(\hat\rho\otimes\id_{C_1})\circ(\id_{C_1}\otimes\varphi_{1,2})\circ
(\hat\rho\otimes\id_{C_2})\circ(\id_{C_1}\otimes\Delta_2)\\
= (\id_{C_1}\otimes\hat\rho)\circ(\hat\rho\otimes\id_{C_2})\circ
(\id_{C_1}\otimes\Delta_2)\,,\\[5pt]
(\varphi_{1,2}\otimes\id_{C_1})\circ(\id_{C_1}\otimes\varphi_{1,2})\circ
(\hat\rho\otimes\id_{C_2})\circ(\id_{C_1}\otimes\Delta_2)\\
=(\id_{C_2}\otimes\hat\rho)\circ(\varphi_{1,2}\otimes\id_{C_2})\circ
(\id_{C_1}\otimes\Delta_2)\,.
\end{gathered}}
\end{gathered}
\end{equation*}
\end{enumerate}
~\endproof
\end{proposition}
\abs
Like cross product algebras the cross product coalgebras are
universal constructions, too.

\begin{proposition}\Label{couniv-constr}
Let $C$ be a coalgebra. Then
\begin{entry}
\item[\ ]
$C$ is coalgebra isomorphic to a cross product coalgebra
$C_1\cpcyco C_2$.
\item[$\Leftrightarrow$]
There is a coalgebra $(C_2,\Delta_2,\varepsilon_2)$ and
an object $C_1$ in $\C$, morphisms $C_1\overset{\proj_1}\leftarrow C
\overset{\proj_2}\to C_2$ and $\varepsilon_1:C_1\to \E$ where
$\proj_2$ is a coalgebra morphism and
$\varepsilon_1\circ\proj_1=\varepsilon_C$ such that
$(\proj_1\otimes\proj_2)\circ\Delta_C:C\to C_1\otimes C_2$
is an isomorphism in $\C$.\endproof
\end{entry}
\end{proposition}
\abs
\begin{proposition}
Let $C_1\cpcyco C_2$ be a cross product coalgebra and
$C$ be a coalgebra. Suppose that there exist morphisms $a$ and $b$
such that
\begin{enumerate}
\item
$a:C\to C_1$ and $\varepsilon_1\circ a=\varepsilon_C\,.$
\item
$b:C\to C_2$ is a coalgebra morphism.
\item
$\varphi_{1,2}\circ(a\otimes b)\circ\Delta_C=(b\otimes a)\circ\Delta_C\,.$
\item
$\hat\rho\circ(a\otimes b)\circ\Delta_C=(a\otimes a)\circ\Delta_C\,.$
\end{enumerate}
Then there exists a unique coalgebra morphism $c:C\to C_1\cpcyco C_2$
obeying the identities $a=(\id_{C_1}\otimes\varepsilon_2)\circ c$ and
$b=(\varepsilon_1\otimes\id_{C_2})\circ c$.\endproof
\end{proposition}
\abs
Co-smash product coalgebras are special cases of cross
product coalgebras. More precisely, a {\it co-smash product coalgebra}
$C=C_1\times^{\phi_{1,2}}C_2$ is given by coalgebras $C_1$ and $C_2$,
and a morphism $\phi_{1,2}:C_1\otimes C_2\to C_2\otimes C_1$,
such that $C=C_1\otimes C_2$ is a coalgebra through
$\Delta_C=(\id_{C_1}\otimes\phi_{1,2}\otimes\id_{C_2})\circ
(\Delta_1\otimes\Delta_2)$ and
$\varepsilon_C=\varepsilon_{C_1}\otimes\varepsilon_{C_2}$.
Then the corresponding result of Proposition \ref{smash-ccpa} holds
for co-smash product coalgebras.

\begin{proposition}
Suppose $C_2$ is a coalgebra and there is a morphism $\eta_2:\E\to C_2$
with $\varepsilon_2\circ\eta_2=\id_{\E}$. Then
\begin{entry}
\item[\ ]
$(C_1,C_2,\varphi_{1,2},\hat\rho)$ is a cross product coalgebra
with $\hat\rho=\Delta_1\otimes\varepsilon_2$ for some morphism
$\Delta_2:C_2\to C_2\otimes C_2$.
\item[$\Leftrightarrow$]
$C_1\times^{\varphi_{1,2}}C_2$ is a co-smash product coalgebra.\endproof
\end{entry}
\end{proposition}
\abs
\begin{remark}[$\pi$-Symmetry]\Label{pi-sym}{\normalfont We would like
to point out the following important observation to the reader. Definition
\ref{caat} is not dual to Definition \ref{ccat}.
Rather, both definitions can be obtained from each other by
a combination of duality (followed by the usual exchanges of
multiplication and comultiplication, unit and counit, cocycle and cycle,
etc.) and use of the opposite tensor product (followed by exchange of indices
"1 $\leftrightarrow$ 2", and later also by exchange of "left/right
(coaction) $\leftrightarrow$ right/left (action)").
This is not an accidential fact but has to do
with the subsequent definition of cross product bialgebras where we
use precisely these algebra and coalgebra structures. Only then are we
able to get all the known cross product bialgebras
\cite{Rad1:85,Ma1:90,MS:94,Ma6:94,BD1:98}.
This kind of symmetry between cross product algebra and
cross product coalgebra will be called 
\textbf{$\mathbf{\pi}$-symmetry} henceforth. In terms of graphical calculus
$\pi$-symmetry can be interpreted easily as rotation of the graphic
by the angle $\pi$ along an axis normal to the planar graphic
followed by the above mentioned exchanges of
indices and morphism types. We will often apply this principle to obtain
$\pi$-symmetric results simply by $\pi$-rotation.
Hence cross product bialgebras - as defined below - are therefore
$\pi$-symmetric invariant rather than dual symmetric invariant.
Dualization of our definition would lead to another kind of cross 
product bialgebras. The corresponding results for these dual
versions follow from our results straightforwardly by dualization.}
\end{remark}

\section{Cross Product Bialgebras}\Label{sec-ccpb}

In Section \ref{sec-ccpb} we investigate cross product bialgebras
from a universal point of view.
They are simultaneously cross product algebras and
cross product coalgebras with compatible bialgebra structure.
We find a universal description for the isomorphism classes of
cross product bialgebras. Then we will consider
special cases called cocycle and strong cross product bialgebras
which are universal constructions as well. But in addition the given
isomorphism class of a cocycle/strong cross product bialgebra
can be characterized by certain interrelated
(co\n-)modular co-cyclic structures on the tensor factors.
Strong cross product bialgebras are the pivotal objects
for the studies in Section \ref{hopf-data}.

\subsection{Cross Product Bialgebras -- General Definition}

\begin{definition}\Label{cbat}
A bialgebra $B$ is called cross product bialgebra
if its underlying algebra is a cross product algebra
$B_1\cpcyalg B_2$, and its underlying coalgebra is a
cross product coalgebra $B_1\cpcyco B_2$ on the same objects.
The cross product bialgebra $B$ will be denoted by $B_1\cpcybi B_2$.
A cross product bialgebra is called normalized if
$\varepsilon_1\circ\eta_1=\id_{\E}$ (and then equivalently
$\varepsilon_2\circ\eta_2=\id_{\E}$).
\end{definition}
\abs
Cross product bialgebras are universal in the following sense.

\begin{theorem}\Label{proj-cycle}
Let $B$ be a bialgebra in $\C$.
Then the subsequent equivalent conditions hold.
\begin{enumerate}
\item\Label{proj-cycle3}
$B$ is bialgebra isomorphic to a normalized cross product
bialgebra $B_1\cpcybi B_2$.
\item\Label{proj-cycle1}
There are idempotents $\Pi_1,\Pi_2\in\End(B)$ such that
\begin{equation}\mathlabel{pi-rel1}
\begin{gathered}
\m_B\circ(\Pi_1\otimes\Pi_1)=\Pi_1\circ\m_B\circ(\Pi_1\otimes\Pi_1)\,,\qquad
\Pi_1\circ\eta_B=\eta_B\,,\\
(\Pi_2\otimes\Pi_2)\circ\Delta_B
 =(\Pi_2\otimes\Pi_2)\circ\Delta_B\circ\Pi_2\,,\qquad
\varepsilon_B\circ\Pi_2=\varepsilon_B\,,
\end{gathered}
\end{equation}
and the sequence
$B\otimes B\xrightarrow{\m_B\circ(\Pi_1\otimes\Pi_2)}B
         \xrightarrow{(\Pi_1\otimes\Pi_2)\circ\Delta_B} B\otimes B$
is a splitting of the idempotent $\Pi_1\otimes\Pi_2$ of $B\otimes B$.
\item\Label{proj-cycle2}
There are objects $B_1$ and $B_2$ and
morphisms $B_1\overset{\inj_1}\to B\overset{\proj_1}\to B_1$
and $B_2\overset{\inj_2}\to B\overset{\proj_2}\to B_2$ where
$\inj_1$ is algebra morphism, $\proj_2$ is coalgebra morphism and
$\proj_j\circ\inj_j=\id_{B_j}$ for $j\in\{1,2\}$ such that
$\m_B\circ(\inj_1\otimes\inj_2):B_1\otimes B_2\to B$ and
$(\proj_1\otimes\proj_2)\circ\Delta_B:B\to B_1\otimes B_2$
are mutually inverse isomorphisms.
\end{enumerate}
\end{theorem}

\begin{proof}
``\ref{proj-cycle1}.$\Rightarrow$\ref{proj-cycle2}.'':
For $j\in\{1,2\}$ we define $\inj_j$, $\proj_j$ to be the morphisms
splitting the idempotent $\Pi_j$ as $\Pi_j=\inj_j\circ\proj_j$ and
$\proj_j\circ\inj_j=\id_{B_j}$ for some objects $B_j$. Then with the
help of Theorem \ref{proj-cycle}.\ref{proj-cycle1} it follows
$\m_B\circ(\inj_1\otimes\inj_2)\circ(\proj_1\otimes\proj_2)\circ\Delta_B
=\m_B\circ(\Pi_1\otimes\Pi_2)\circ\Delta_B=\id_B$ and
$(\proj_1\otimes\proj_2)\circ\Delta_B\circ\m_B\circ(\inj_1\otimes\inj_2)
= (\proj_1\otimes\proj_2)\circ(\Pi_1\otimes\Pi_2)\circ
 \Delta_B\circ\m_B\circ(\Pi_1\otimes\Pi_2)\circ(\inj_1\otimes\inj_2)
=(\proj_1\circ\inj_1\otimes\proj_2\circ\inj_2)=\id_{B_1\otimes B_2}$. Hence
$(\m_B\circ(\inj_1\otimes\inj_2))^{-1}=(\proj_1\otimes\proj_2)\circ\Delta_B$.
Now we define $\m_1:=\proj_1\circ\m_B\circ(\inj_1\otimes\inj_1)$,
$\eta_1:=\proj_1\circ\eta_B$,
$\Delta_2:=(\proj_2\otimes\proj_2)\circ\Delta_B\circ\inj_2$, and
$\varepsilon_2:=\varepsilon_B\circ\inj_2$. Then using
\eqref{pi-rel1} and the (co\n-)algebra property of $B$ one verifies
easily that $(B_1,\m_1,\eta_1)$ is an algebra and
$(B_2,\Delta_2,\varepsilon_2)$ is a coalgebra.
Furthermore $\inj_1$ is an algebra morphism because
$\inj_1\circ\m_1 = \Pi_1\circ\m_B\circ(\inj_1\otimes\inj_1)
= \m_B\circ(\inj_1\otimes\inj_1)$, and
$\inj_1\circ\eta_1=\Pi_1\circ\eta_B = \eta_B$. In a $\pi$-symmetric manner
one proves that $\proj_2$ is a coalgebra morphism.
\nl
``\ref{proj-cycle2}.$\Rightarrow$\ref{proj-cycle1}.'':
Conversely we define $\Pi_j:=\inj_j\circ\proj_j$ for $j\in\{1,2\}$.
Since by assumption $\inj_1$ is an algebra morphism it follows
$\Pi_1\circ\eta_B=\Pi_1\circ\inj_1\circ\eta_1=\inj_1\circ\eta_1=\eta_B$
and using $\proj_1\circ\inj_1=\id_{B_1}$ it holds
$\Pi_1\circ\m_B\circ(\Pi_1\otimes\Pi_1) =
 \inj_1\circ\proj_1\circ\inj_1\circ\m_1\circ(\proj_1\otimes\proj_1)
= \m_B\circ(\Pi_1\otimes\Pi_1)$.
In $\pi$-symmetrical analogy it will be proven that
$(\Pi_2\otimes\Pi_2)\circ\Delta_B=(\Pi_2\otimes\Pi_2)\circ\Delta_B\circ\Pi_2$
and $\varepsilon_B\circ\Pi_2=\varepsilon_B$.
From Theorem \ref{proj-cycle}.\ref{proj-cycle2} we conclude
$(\Pi_1\otimes\Pi_2)\circ\Delta_B\circ\m_B\circ(\Pi_1\otimes\Pi_2)
=(\inj_1\otimes\inj_2)\circ\id_{B_1\otimes B_2}\circ
(\proj_1\otimes\proj_2)=\Pi_1\otimes\Pi_2$ and
$\id_B=\m_B\circ(\Pi_1\otimes\Pi_2)\circ\Delta_B=
\m_B\circ(\Pi_1\otimes\Pi_2)\circ(\Pi_1\otimes\Pi_2)\circ\Delta_B$.
Thus $\m_B\circ(\Pi_1\otimes\Pi_2)$ and $(\Pi_1\otimes\Pi_2)\circ\Delta_B$
split the idempotent $(\Pi_1\otimes\Pi_2)$.
\nl
``\ref{proj-cycle2}.$\Rightarrow$\ref{proj-cycle3}.'':
By assumption $(B_1,\m_1,\eta_1)$ is an algebra and $\inj_1$ is
an algebra morphism, $\inj_2$ is a morphism in $\C$, and
$\m_B\circ(\inj_1\otimes\inj_2)$ is an isomorphism.
Define $\eta_2:=\proj_2\circ\eta_B:\E\to B_2$. Then
$\inj_2\circ\eta_2 =\Pi_2\circ\eta_B =\m_B\circ(\Pi_1\circ\eta_B\otimes
\Pi_2\circ\eta_B) =\m_B\circ(\inj_1\otimes\inj_2)\circ(\proj_1\otimes\proj_2)
\circ\Delta_B\circ\eta_B =\eta_B$. Therefore all conditions of
Proposition \ref{univ-constr}.2 are fulfilled implying that
$B$ is algebra isomorphic to a cross product algebra
$B_1\cpcyalg B_2$. From the proof of Proposition \ref{univ-constr} one
reads off that the isomorphism is given by
$\Lambda=\m_B\circ(\inj_1\otimes\inj_2)$. In a $\pi$-symmetric way
one shows that $B$ is coalgebra isomorphic to a cross product
coalgebra $B_1\cpcyco B_2$ on the same tensor product with
isomorphism $\tilde\Lambda=(\proj_1\otimes\proj_2)\circ\Delta_B$. Thus
by assumption $\tilde\Lambda=\Lambda^{-1}$. Since $B$ is a bialgebra
it follows then from Definition \ref{cbat} that $B_1\otimes B_2$ is
cross product bialgebra $B_1\cpcybi B_2$ and
$\Lambda:B_1\cpcybi B_2\to B$ is bialgebra isomorphism. Furthermore
the identities $\id_\E=\varepsilon_B\circ\eta_B=\varepsilon_B\circ
\inj_2\circ\eta_2=\varepsilon_2\circ\proj_2\circ\inj_2\circ\eta_2=
\varepsilon_2\circ\eta_2$ show that $B_1\cpcybi B_2$ is normalized.
\nl
``\ref{proj-cycle3}.$\Rightarrow$\ref{proj-cycle2}.'':
Given a bialgebra isomorphism $\Lambda:B_1\cpcybi B_2\to B$.
Then in particular $\Lambda:B_1\cpcyalg B_2\to B$ is an algebra
isomorphism and $\Lambda:B_1\cpcyco B_2\to B$ is a coalgebra isomorphism.
We define $\inj_1 :=\Lambda\circ(\id_{B_1}\otimes\eta_2)$,
$\inj_2 :=\Lambda\circ(\eta_1\otimes\id_{B_2})$,
$\proj_1 :=(\id_{B_1}\otimes\varepsilon_2)\circ\Lambda^{-1}$, and
$\proj_2 :=(\varepsilon_1\otimes\id_{B_2})\circ\Lambda^{-1}$.
From Propositions \ref{univ-constr} and \ref{couniv-constr} one
derives that $\inj_1$ is an algebra morphism, $\proj_2$ is a coalgebra
morphism, and $\Lambda=\m_B\circ(\inj_1\otimes\inj_2)=
\big((\proj_1\otimes\proj_2)\circ\Delta_B\big)^{-1}$.
By assumption $B_1\cpcybi B_2$ is normalized from which follows
$\proj_1\circ\inj_1=(\id_{B_1}\otimes\varepsilon_2)\circ
\Lambda^{-1}\circ\Lambda\circ(\id_{B_1}\otimes\eta_2)=\id_{B_1}$ and
similarly $\proj_2\circ\inj_2=\id_{B_2}$.\end{proof}
\abs
The next corollary of Theorem \ref{proj-cycle} will be used
frequently in the following.

\begin{corollary}\Label{proj-cycle-eqs}
Under the conditions of Theorem \ref{proj-cycle} it holds
\begin{equation}\mathlabel{cons-proj-cycle1}
\begin{gathered}
\begin{split}
\Pi_2\circ\eta_B&=\eta_B\,,\\
\varepsilon_B\circ\Pi_1&=\varepsilon_B\,,
\end{split}
\quad
\begin{split}
\Pi_1\circ\Pi_2&=\eta_B\circ \varepsilon_B\,,\\
\Pi_2\circ\Pi_1&=\eta_B\circ \varepsilon_B\,,
\end{split}
\quad
\begin{split}
\Pi_2\circ\m_B &=\Pi_2\circ\m_B\circ(\Pi_2\otimes\id_B)\,,\\
\Delta_B\circ \Pi_1&=(\id_B\otimes\Pi_1)\circ\Delta_B\circ\Pi_1\,,
\end{split}\\[5pt]
{\begin{split}
\Pi_1\circ\m_B\circ(\Pi_1\otimes\id_B)&=
\Pi_1\circ\m_B\circ(\Pi_1\otimes\Pi_1)\,,\\
(\id_B\otimes\Pi_2)\circ\Delta_B\circ \Pi_2&=
(\Pi_2\otimes\Pi_2)\circ\Delta_B\circ \Pi_2\,.
\end{split}}
\end{gathered}
\end{equation}
The structure morphisms $\varphi_{1,2}$, $\varphi_{2,1}$, $\hat\rho$,
and $\hat\sigma$ of the cross product bialgebra
$B_1\cpcybi B_2$ can be expressed through the projections
and injections on $B$ as
\begin{equation}\mathlabel{struc-morph}
\begin{gathered}
\inj_1\circ\eta_1=\eta_B\,,\quad
\varepsilon_1\circ\proj_1=\varepsilon_B\,,\quad
\inj_2\circ\eta_2=\eta_B\,,\quad
\varepsilon_2\circ\proj_2= \varepsilon_B\,,\\[5pt]
{\begin{split}
\m_1&=\proj_1\circ\m_B\circ(\inj_1\otimes\inj_1)\,,\\
\varphi_{1,2} &=(\proj_2\otimes\proj_1)\circ\Delta_B\circ\m_B\circ
      (\inj_1\otimes\inj_2)\,,\\
\varphi_{2,1} &=(\proj_1\otimes\proj_2)\circ\Delta_B\circ\m_B\circ
      (\inj_2\otimes\inj_1)\,,
\end{split}}
\quad
{\begin{split}
\Delta_2&=(\proj_2\otimes\proj_2)\circ\Delta_B\circ\inj_2\,,\\
\hat\rho &=(\proj_1\otimes\proj_1)\circ\Delta_B\circ\m_B\circ
      (\inj_1\otimes\inj_2)\,,\\
\hat\sigma &=(\proj_1\otimes\proj_2)\circ\Delta_B\circ\m_B\circ
      (\inj_2\otimes\inj_2)\,.
\end{split}}
\end{gathered}
\end{equation}
\end{corollary}

\begin{proof}
The identities $\inj_1\circ\eta_1=\eta_B$ and
$\varepsilon_2\circ\proj_2=\varepsilon_B$ hold since $\inj_1$ is an
algebra morphism and $\proj_2$ is a coalgebra morphism, whereas
the identities $\inj_2\circ\eta_2=\eta_B$ and
$\varepsilon_1\circ\proj_1=\varepsilon_B$, as well as the equations
$\Pi_2\circ\eta_B=\eta_B$ and $\varepsilon_B\circ\Pi_1=\varepsilon_B$
have been shown in the proof of Theorem \ref{proj-cycle}.
From the proofs of Propositions \ref{univ-constr} and \ref{couniv-constr}
one can directly read off the structure of the morphisms
$\varphi_{1,2}$, $\varphi_{2,1}$, $\hat\rho$, and $\hat\sigma$
given in \eqref{struc-morph} using
$\Lambda=\m_B\circ(\inj_1\otimes\inj_2)$ and
$\Lambda^{-1}=(\proj_1\otimes\proj_2)\circ\Delta_B$.
Since $\inj_1$ is algebra morphism, $\proj_2$ is coalgebra
morphism and $\proj_j\circ\inj_j=\id_{B_j}$ for $j\in\{1,2\}$ it
follows $\m_1=\proj_1\circ\m_B\circ(\inj_1\otimes\inj_1)$ and
$\Delta_2=(\proj_2\otimes\proj_2)\circ\Delta_B\circ\inj_2$.
Because of \eqref{pi-rel1} and Theorem \ref{proj-cycle}.\ref{proj-cycle1}
it follows
$\Pi_1\circ\Pi_2=(\id_B\otimes\varepsilon_B)\circ(\Pi_1\otimes\Pi_2)
\circ\Delta_B\circ\m_B\circ{(\Pi_1\otimes\Pi_2)}\circ(\eta_B\otimes\id_B)=
(\Pi_1\circ\eta_B\otimes\varepsilon_B\circ\Pi_2)=\eta_B\circ\varepsilon_B$.
Using $\Pi_2\circ\eta_B=\eta_B$ and $\varepsilon_B\circ\Pi_1=\varepsilon_B$
from above we obtain in a similar way
$\Pi_2\circ\Pi_1=\eta_B\circ\varepsilon_B$.
Using consecutively Theorem \ref{proj-cycle}.\ref{proj-cycle1} two times,
\eqref{pi-rel1}, the identity
$\Pi_2=(\varepsilon_B\circ\Pi_1\otimes\Pi_2)\circ\Delta_B$,
again Theorem \ref{proj-cycle}.\ref{proj-cycle1}, and the identity
$\varepsilon_B\circ\Pi_1=\varepsilon_B$
one obtains the following series of equations.
\begin{equation*}
\begin{split}
\Pi_2\circ\m_B&=\Pi_2\circ\m_B\circ(\m_B\circ(\Pi_1\otimes\Pi_2)\circ\Delta_B
 \otimes\id_B)\\
&=\Pi_2\circ\m_B^{(3)}\circ(\Pi_1\otimes\Pi_1\otimes\Pi_2)\circ
 (\id_B\otimes\Delta_B\circ\m_B)\circ\\
&\quad\circ(\Pi_1\otimes\Pi_2\otimes\id_B)\circ(\Delta_B\otimes\id_B)\\
&=\Pi_2\circ\m_B\circ(\Pi_1\otimes\Pi_2)\circ
 (\m_B\circ(\Pi_1\otimes\Pi_1)\otimes\id_B)\circ\\
&\quad\circ(\id_B\otimes\Delta_B\circ\m_B)\circ(\id_B\otimes\Pi_2\otimes\id_B)
 \circ(\Delta_B\otimes\id_B)\\
&=\Pi_2\circ\m_B\circ(\Pi_2\otimes\id_B)
\end{split}
\end{equation*}
and $\pi$-symmetrically
$\Delta_B\circ\Pi_1= (\id_B\otimes\Pi_1)\circ\Delta_B\circ\Pi_1$.
Finally one obtains in a similar manner
\begin{equation*}
\begin{split}
\Pi_1\circ\m_B\circ(\Pi_1\otimes\id_B) &= \Pi_1\circ\m_B\circ
 (\Pi_1\otimes\m_B\circ(\Pi_1\otimes\Pi_2)\circ\Delta_B)\\
&=\Pi_1\circ\m_B^{(3)}\circ(\Pi_1\otimes\Pi_1\otimes\Pi_2)\circ
  (\id_B\otimes\Delta_B)\\
&=\Pi_1\circ\m_B\circ(\Pi_1\otimes\Pi_2)\circ
 (\m_B\circ(\Pi_1\otimes\Pi_1)\otimes\Pi_2)\circ(\id_B\otimes\Delta_B)\\
&=(\Pi_1\otimes\varepsilon_B\circ\Pi_2)\circ\Delta_B\circ\m_B
  \circ(\Pi_1\otimes\Pi_2)\circ\\
&\quad\circ(\m_B\circ(\Pi_1\otimes\Pi_1)\otimes\Pi_2)
 \circ(\id_B\otimes\Delta_B)\\
&=\Pi_1\circ\m_B\circ(\Pi_1\otimes\Pi_1)
\end{split}
\end{equation*}
and $(\id_B\otimes\Pi_2)\circ\Delta_B\circ\Pi_2=
(\Pi_2\otimes\Pi_2)\circ\Delta_B\circ\Pi_2$ by $\pi$-symmetry.\end{proof}

\begin{remark}\Label{alpha-ip}
{\normalfont The statement of Theorem \ref{proj-cycle} can be
refined in the following sense. For a given bialgebra $B$
the tuples $(B,\Lambda,B_1\cpcybi B_2)$
and $(B,\inj_1,\inj_2,\proj_1,\proj_2)$ obeying the conditions of
Theorem \ref{proj-cycle}.1 and \ref{proj-cycle}.3
respectively, are in one-to-one correspondence. This correspondence 
has been constructed in the proof of Theorem \ref{proj-cycle}.}
\end{remark}
\abs
Theorem \ref{proj-cycle} and the previous remark imply a useful corollary.
Given a normalized cross product bialgebra $B_1\cpcybi B_2$.
We set $\inj_{\alpha\,1}:= \id_{B_1}\otimes\eta_2\,,$ $\inj_{\alpha\,2}:=
\eta_1\otimes\,\id_{B_2}\,,$ $\proj_{\alpha\,1}:=\id_{B_1}\otimes\,\varepsilon_2$ and
$\proj_{\alpha\,2}:=\varepsilon_1\otimes\id_{B_2}$. Then we define
\begin{equation}\mathlabel{m-delta-ijk}
\begin{split}
\m_{\bowtie\,i,jk}&:=\proj_{\alpha\,i}\circ\m_{B_1\cpcybi B_2}\circ
 (\inj_{\alpha\,j}\otimes\inj_{\alpha\,k})\,,\\
\Delta_{\bowtie\, ij,k}&:=(\proj_{\alpha\,i}\otimes\proj_{\alpha\,j})\circ
\Delta_{B_1\cpcybi B_2}\circ\inj_{\alpha\,k}
\end{split}
\end{equation}
for $i,j,k\in\{1,2\}$. On the other hand suppose that $B$ is a bialgebra which
obeys the conditions of Theorem \ref{proj-cycle}.3. Then we define
\begin{equation}\mathlabel{m-delta-ijk1}
\m_{B\,i,jk}:=\proj_i\circ\m_B\circ(\inj_j\otimes\inj_k)\,,\quad
\Delta_{B\,ij,k}:=(\proj_i\otimes\proj_j)\circ\Delta_B\circ\inj_k
\end{equation}
for $i,j,k\in\{1,2\}$.

\begin{corollary}\Label{co-act-inv}
Let $B$ be bialgebra and suppose that the tuples 
$(B,\Lambda,B_1\cpcybi B_2)$ and
$(B,\inj_1,\inj_2,\proj_1,\proj_2)$ are related according to
Theorem \ref{proj-cycle}. Then $\m_{\bowtie\,i,jk}=\m_{B\,i,jk}$ and
$\Delta_{\bowtie\,ij,k}= \Delta_{B\,ij,k}$.
\end{corollary}

\begin{proof} If the conditions of Theorem \ref{bialg-cocycle2} hold,
then by construction (see the proof of the theorem) the injections and
projections of $B$ and the morphisms $\inj_{\alpha\,1}$, $\inj_{\alpha\,2}$,
$\proj_{\alpha\,1}$ and $\proj_{\alpha\,2}$ are related by
$\inj_1 :=\Lambda\circ\inj_{\alpha\,1}$,
$\inj_2 :=\Lambda\circ\inj_{\alpha\,2}$,
$\proj_1 :=\proj_{\alpha\,1}\circ\Lambda^{-1}$, and
$\proj_2 :=\proj_{\alpha\,2}\circ\Lambda^{-1}$, where
$\Lambda:B_1\otimes B_2\to B$ is the given bialgebra isomorphism.
Thus $\m_{\bowtie\,i,jk}=\m_{B\,i,jk}$ and $\Delta_{\bowtie\,ij,k}=
\Delta_{B\,ij,k}$ follow directly.\end{proof}
\abs
Hence we will use the notation $\m_{i,jk}$ and $\Delta_{ij,k}$
for the corresponding morphisms henceforth.

\subsection{Cocycle Cross Product Bialgebras}

Up to now a (co\n-)modular co-cyclic structure of
cross product bialgebras and their tensor factors $B_1$ and $B_2$ did not
emerge. In Definition \ref{type-alpha} below we will restrict our 
considerations to cocycle cross product bialgebras for which (co\n-)modular
co-cyclic structures appear in a very natural way. The universal property
of cocycle cross product bialgebras will be discussed
subsequently. We define strong cross product
bialgebras in Definition \ref{strong-type-alpha}. They are the basic
objects which eventually constitute the universal,
(co\n-)modular co-cyclic theory of cross product bialgebras.

In Definition \ref{cbat} the morphisms $\m_1$,
$\eta_1$, $\varepsilon_1$, $\Delta_2$, $\eta_2$, $\varepsilon_2$, etc.~occur. 
Definition \ref{type-alpha} below requires additional morphisms
$\Delta_1$, $\m_2$, $\mu_l$, $\mu_r$, $\nu_l$, $\nu_r$, $\sigma$, and $\rho$.
This is the stage where we find it convenient to start working with
graphical calculus. All $\m$, $\Delta$, $\eta$, $\varepsilon$,
$\mu$, and $\nu$ will be presented by the graphics displayed
in Figure \ref{fig-conv} respectively\footnote{Note that some of these
morphisms might be weak. That is, weak (co\n-)\-mul\-ti\-plications, weak
(co\n-)\-act\-ions, etc.}. The cocycle and cycle
morphisms $\sigma$ and $\rho$ will be presented henceforth graphically
by $\sigma=\begin{tangle}\hh \cu* \end{tangle}:B_2\otimes B_2\to B_1$ and
$\rho=\begin{tangle}\hh \cd* \end{tangle}:B_2\to B_1\otimes B_1$.

\begin{definition}\Label{type-alpha}
A normalized cross product bialgebra $B_1\cpcybi B_2$ is called
cocycle cross product bialgebra if there exist additional morphisms
\begin{equation*}
\begin{gathered}
\Delta_1:B_1\to B_1\otimes B_1\,,\quad\m_2:B_2\otimes B_2\to B_2\,,\\
\mu_l:B_2\otimes B_1\to B_1\,,\quad\mu_r:B_2\otimes B_1\to B_2\,,\\
\nu_l:B_1\to B_2\otimes B_1\,,\quad\nu_r:B_2\to B_2\otimes B_1\,,\\
\sigma:B_2\otimes B_2\to  B_1\,,\quad\rho:B_2\to B_1\otimes B_1
\end{gathered}
\end{equation*}
such that $\varphi_{1,2}$, $\varphi_{2,1}$, $\hat\sigma$ and $\hat\rho$
are of the form
\begin{equation*}
\begin{split}
\varphi_{1,2}&=(\m_2\otimes\m_1)\circ(\id_{B_1}\otimes\Psi_{B_1,B_2}\otimes
\id_{B_2})\circ(\nu_l\otimes\nu_r)\,,\\
\varphi_{2,1}&=(\mu_l\otimes\mu_r)\circ(\id_{B_2}\otimes\Psi_{B_2,B_1}
\otimes\id_{B_1})\circ(\Delta_2\otimes\Delta_1)\,,\\
\hat\sigma&=(\sigma\otimes\m_2\circ(\mu_r\otimes\id_{B_2}))\circ
(\id_{B_2}\otimes\Psi_{B_2,B_2}\otimes\id_{B_1\otimes B_2})\circ
(\Delta_2\otimes(\nu_r\otimes\id_{B_2})\circ\Delta_2)\,,\\
\hat\rho&=(\m_1\circ(\id_{B_1}\otimes\mu_l)\otimes\m_1)\circ
(\id_{B_1\otimes B_2}\otimes\Psi_{B_1,B_1}\otimes\id_{B_1})\circ
((\id_{B_1}\otimes\nu_l)\circ\Delta_1\otimes\rho)\,.
\end{split}
\end{equation*}
Furthermore we require the following (co\n-)unital identities
\begin{equation}\mathlabel{type-alpha1}
\begin{gathered}
\Delta_1\circ\eta_1=\eta_1\otimes\eta_2\,,\quad
(\id_{B_1}\otimes\varepsilon_1)\circ\Delta_1=
(\varepsilon_1\otimes\id_{B_1})\circ\Delta_1=\id_{B_1}\,,\\
\varepsilon_2\circ\m_2=\varepsilon_2\otimes\varepsilon_2\,,\quad
\m_2\circ(\id_{B_2}\otimes\eta_2)=\m_2\circ(\eta_2\otimes\id_{B_2})=
\id_{B_2}\,,\\[5pt]
\mu_l\circ(\eta_2\otimes\id_{B_1})=\id_{B_1}\,,\quad
\mu_l\circ(\id_{B_2}\otimes\eta_1)=\eta_1\circ\varepsilon_2\,,\quad
\varepsilon_1\circ\mu_l=\varepsilon_2\otimes\varepsilon_1\,,\\
\mu_r\circ(\id_{B_2}\otimes\eta_1)=\id_{B_2}\,,\quad
\mu_r\circ(\eta_2\otimes\id_{B_1})=\eta_2\circ\varepsilon_1\,,\quad
\varepsilon_2\circ\mu_r=\varepsilon_2\otimes\varepsilon_1\,,\\
(\varepsilon_2\otimes\id_{B_1})\circ\nu_l=\id_{B_1}\,,\quad
(\id_{B_2}\otimes\varepsilon_1)\circ\nu_l=\eta_2\circ\varepsilon_1\,,
\quad\nu_l\circ\eta_1=\eta_2\otimes\eta_1\,,\\
(\id_{B_2}\otimes\varepsilon_1)\circ\nu_r=\id_{B_2}\,,\quad
(\varepsilon_2\otimes\id_{B_1})\circ\nu_r=\eta_1\circ\varepsilon_2\,,
\quad\nu_r\circ\eta_2=\eta_2\otimes\eta_1
\end{gathered}
\end{equation}
and the ``projection" relations
\begin{equation}\mathlabel{type-alpha2}
\begin{split}
(\m_2\otimes\id_{B_1})\circ(\id_{B_2}\otimes\nu_r)&=
(\m_2\otimes\id_{B_1})\circ(\mu_r\otimes\varphi_{1,2})\circ
(\id_{B_2}\otimes(\rho\otimes\id_{B_2})\circ\Delta_2)\,,\\
(\mu_l\otimes\id_{B_1})\circ(\id_{B_2}\otimes\Delta_1)&=
(\m_1\circ(\id_{B_1}\otimes\sigma)\otimes\id_{B_1})\circ
(\varphi_{2,1}\otimes\nu_l)\circ(\id_{B_2}\otimes\Delta_1)\,.
\end{split}
\end{equation}
Cocycle cross product bialgebras will be denoted subsequently
by $B_1\cpcybi[\mathfrak{c}] B_2$.
\end{definition}

\begin{remark}\Label{unit-prop}
{\normalfont For a cocycle cross product bialgebra
the following (co\n-)unital identities can be derived easily.
\begin{equation*}
\begin{gathered}
\Delta_2\circ\eta_2=\eta_2\otimes\eta_2\,,\quad\varepsilon_1\circ\m_1=
\varepsilon_1\otimes\varepsilon_1\,,\\
\varepsilon_1\circ\sigma=\varepsilon_2\otimes\varepsilon_2\,,\quad
\rho\circ\eta_2=\eta_1\otimes\eta_1\,,\\
\sigma\circ(\eta_2\otimes\id_{B_2})=\sigma\circ(\id_{B_2}\otimes\eta_2)
= \eta_1\circ\varepsilon_2\,,\\
(\id_{B_1}\otimes\varepsilon_1)\circ\rho=
(\varepsilon_1\otimes\id_{B_1})\circ\rho=
\eta_1\circ\varepsilon_2\,.
\end{gathered}
\end{equation*}}
\end{remark}

\begin{remark}\Label{alpha-struc}{\normalfont
By definition, cocycle cross product bialgebras always come with
certain structure morphisms $\m_1$, $\m_2$, $\Delta_1$, $\Delta_2$, etc.
There might be two different cocycle cross product bialgebras
with the same underlying bialgebra structure. It is understood
henceforth, that the notation $B_1\cpcybi[\mathfrak{c}] B_2$ always refers to
the complete structure of cocycle cross product bialgebras.}
\end{remark}
\abs
The graphics of the morphisms $\varphi_{1,2}$, $\varphi_{2,1}$,
$\hat\rho$ and $\hat\sigma$ of a cocycle cross product
bialgebra $B_1\cpcybi[\mathfrak{c}] B_1$ are given by
\begin{equation}\mathlabel{phi-sigma-rho}
\begin{array}{c}
\hstretch 125
\vstretch 85
\varphi_{1,2}
:=
\divide\unitlens by 3
\begin{tangle}
\hh\step\id\step\id\\
\hh\ld\step\rd\step\\
\id\step\hx\step\id\\
\hh\cu\step\cu\\
\hh\hstep\id\Step\id
\end{tangle}
\multiply\unitlens by 3
\quad ,\quad
\varphi_{2,1}
:=
\divide\unitlens by 3
\begin{tangle}
\hh\hstep\id\Step\id\\
\hh\cd\step\cd\\
\id\step\hx\step\id\\
\hh\lu\step\ru\\
\hh\step\id\step\id
\end{tangle}
\multiply\unitlens by 3
\quad ,\quad
\hat\rho
:=
\hstretch 150
\vstretch 85
\divide\unitlens by 3
\begin{tangle}
\hh\hstep\id\Step\id\\
\hh\cd\step\cd*\\
\hh \id\hstep\hld\step\id\step\id\\
\id\hstep\id\hstep\hx\step\id\\
\hh \id\hstep\hlu\step\id\step\id\\
\hh \cu\step\cu\\
\hh\hstep\id\Step\id
\end{tangle}
\multiply\unitlens by 3
\quad ,\quad
\hat\sigma
:=
\divide\unitlens by 3
\begin{tangle}
\hh\hstep\id\Step\id\\
\hh\cd\step\cd \\
\hh\id\step\id\step\hrd\hstep\id \\
\id\step\hx\hstep\id\hstep\id \\
\hh\id\step\id\step\hru\hstep\id \\
\hh\cu*\step\cu\\
\hh\hstep\id\Step\id
\end{tangle}
\multiply\unitlens by 3
\end{array}
\end{equation}
\abs
Then the graphical shapes of the multiplication and the
comultiplication look like
\begin{equation}\mathlabel{type-alpha3}
\begin{array}{c}
\m_{B_1\cpcybi[\mathfrak{c}] B_1}
=
\hstretch 150
\divide\unitlens by 3
\begin{tangle}
\hh \id\step\cd\step\cd\step\id \\
\id\step\id\step\hx\step\id\step\id \\
\hh \id\step\lu\step\ru\step\id \\
\hh \id\Step\id\hstep\cd\step\cd \\
\hh \id\Step\id\hstep\id\step\id\step\hrd\hstep\id \\
\cu\hstep\id\step\hx\hstep\id\hstep\id \\
\hh \step\id\step\hstep\id\step\id\step\hru\hstep\id \\
\hh \step\id\step[1.5]\cu*\step\cu \\
\step\cu\Step\id
\end{tangle}
\multiply\unitlens by 3
\qquad\text{and}\qquad
\Delta_{B_1\cpcybi[\mathfrak{c}] B_1}
=
\divide\unitlens by 3
\begin{tangle}
\hstep\id\Step\cd \\
\hh \cd\step\cd*\step[1.5]\id \\
\hh \id\hstep\hld\step\id\step\id\hstep\step\id \\
\id\hstep\id\hstep\hx\step\id\hstep\cd \\
\hh \id\hstep\hlu\step\id\step\id\hstep\id\Step\id \\
\hh \cu\step\cu\hstep\id\Step\id \\
\hh \hstep\id\step\ld\step\rd\step\id \\
\hstep\id\step\id\step\hx\step\id\step\id \\
\hh \hstep\id\step\cu\step\cu\step\id \end{tangle}
\multiply\unitlens by 3
\end{array}
\end{equation}
\abs
It is an easy exercise to prove that Definition \ref{type-alpha} is compatible
with trivial (co\n-)actions $\mu_l$, $\mu_r$, $\nu_l$, $\nu_r$ and trivial
(co\n-)cycles $\sigma$ and $\rho$ respectively, in the sense that no
additional identities involving (co\n-)units have to be required for the
remaining morphisms which define this special cocycle
cross product bialgebra.

The following theorem is the central basic statement in Section \ref{sec-ccpb}.
It describes universality of cocycle cross product bialgebras. 

\begin{theorem}\Label{bialg-cocycle2}
Let $B$ be a bialgebra in $\C$. Then the following statements are
equivalent.

\begin{enumerate}[1.]
\item
$B$ is bialgebra isomorphic to a cocycle cross product
bialgebra $B_1\cpcybi[\mathfrak{c}] B_2$.
\item
There are idempotents $\Pi_1,\Pi_2\in \End(B)$ such that the
conditions of Theorem \ref{proj-cycle}.\ref{proj-cycle1} and the
following ``projection" relations hold.
\begin{align}
\notag
(\mathbf{\gamma_1})&:&(\id\otimes\Pi_1)\circ\Delta_B
&=\delta_{(1,2,1)}\circ\gamma_1\circ(\Pi_1\otimes\Pi_2)\circ\Delta_B\,,\\
\mathlabel{techn2-cond}
(\mathbf{\gamma_2})&:&\m_B\circ(\Pi_2\otimes\id)
&= \m_B\circ(\Pi_1\otimes\Pi_2)\circ\gamma_2\circ\delta_{(2,1,2)}\,,\\
\notag
(\mathbf{\delta_j})&:
&(\Pi_j\otimes\id)\circ\delta_{(0,0,0)}\circ(\id\otimes\Pi_j)
&=(\Pi_j\otimes\id)\circ\delta_{(0,j,0)}\circ(\id\otimes\Pi_j)
\end{align}
for $j\in\{1,2\}$. In \eqref{techn2-cond} we used the abbreviations
\begin{equation*}
\begin{split}
\gamma_j &:= (\m_B\otimes\m_B)\circ(\id\otimes\Pi_j\otimes\id\otimes\id)
\circ(\id\otimes\Psi_{B,B}\otimes\id)\circ(\Delta_B\otimes\Delta_B)\,,\\
\delta_{(i,j,k)} &:=(\m_B\otimes\id)\circ(\Pi_i\otimes\Pi_j\otimes\Pi_k)
\circ(\id\otimes\Delta_B)\,.
\end{split}
\end{equation*}
for $i,j,k\in\{0,1,2\}$, and we set $\Pi_0:=\id_A$.
\item
There are objects $B_1$ and $B_2$ and
morphisms $B_1\overset{\inj_1}\to B\overset{\proj_1}\to B_1$
and $B_2\overset{\inj_2}\to B\overset{\proj_2}\to B_2$ for which
the conditions of Theorem \ref{proj-cycle}.\ref{proj-cycle2} and the
``projection" relations \eqref{techn2-cond} hold, with
$\Pi_j:=\inj_j\circ\proj_j$, $j\in\{1,2\}$.
\end{enumerate}
\end{theorem}

\begin{proof}
Obviously statement $2.$ and $3.$ are equivalent due to Theorem
\ref{proj-cycle}. Hence suppose statements
$2.$ and $3.$ hold. Then by Theorem \ref{proj-cycle} there
exists a normalized cross product bialgebra $B_1\cpcybi B_2$ which is
isomorphic to $B$ via the isomorphism $\Lambda=\m_B\circ
(\inj_1\otimes\inj_2)$. We define
\begin{equation}\mathlabel{proj-struc}
{\begin{split}
\m_2 &:=\proj_2\circ\m_B\circ(\inj_2\otimes\inj_2)\,,\\
\sigma &:= \proj_1\circ\m_B\circ(\inj_2\otimes\inj_2)\,,\\
\mu_l &:= \proj_1\circ\m_B\circ(\inj_2\otimes\inj_1)\,,\\
\mu_r &:= \proj_2\circ\m_B\circ(\inj_2\otimes\inj_1)\,,
\end{split}}
\quad
{\begin{split}
\Delta_1 &:= (\proj_1\otimes\proj_1)\circ\Delta_B\circ\inj_1\,,\\
\rho &:= (\proj_1\otimes\proj_1)\circ\Delta_B\circ\inj_2\,,\\
\nu_l &:= (\proj_2\otimes\proj_1)\circ\Delta_B\circ\inj_1\,,\\
\nu_r &:= (\proj_2\otimes\proj_1)\circ\Delta_B\circ\inj_2\,.\\
\end{split}}
\end{equation}
The morphisms $\m_1$, $\Delta_2$, $\eta_1$, $\eta_2$, $\varepsilon_1$, and
$\varepsilon_2$ are given in terms of the projections and injections
by the corresponding relations in \eqref{struc-morph} in Corollary
\ref{proj-cycle-eqs}. With the help of these data we define
structure morphisms $\tilde{\varphi}_{1,2}$, $\tilde{\varphi}_{2,1}$,
$\hat{\tilde\sigma}$, and $\hat{\tilde\rho}$ analogous to
\eqref{phi-sigma-rho}. We will show that the morphisms
$\m_{B_1\otimes B_2}$ and $\Delta_{B_1\otimes B_2}$ defined
with these structure morphisms according to \eqref{type-alpha3}
precisely coincide with $\m_{B_1\cpcybi B_2}$ and
$\Delta_{B_1\cpcybi B_2}$ respectively, if the assumptions of
statement $2.$ (or $3.$) are satisfied. Before we prove this
we will provide several auxiliary identities.
From \eqref{pi-rel1} and Corollary \ref{proj-cycle-eqs} we obtain
\begin{gather}\mathlabel{aux1}
{\begin{split}
(\Pi_2\otimes\Pi_2)\circ\Delta_B\circ\m_B\circ(\Pi_2\otimes\id_B) &=
(\Pi_2\otimes\Pi_2)\circ\Delta_B\circ\Pi_2\circ\m_B\circ(\Pi_2\otimes\id_B)\\
&= (\Pi_2\otimes\Pi_2)\circ\Delta_B\circ\m_B\,,
\end{split}}
\\
\mathlabel{aux2}
{\begin{split}
\Pi_1\circ\m_B^{(3)}\circ(\Pi_1\otimes\Pi_1\otimes\Pi_1)
&=\Pi_1\circ\m_B\circ(\id_B\otimes\Pi_1\circ\m_B)\circ
   (\Pi_1\otimes\Pi_1\otimes\Pi_1)\\
&=\Pi_1\circ\m_B\circ(\id_B\otimes\Pi_1\circ\m_B)\circ
   (\Pi_1\otimes\Pi_1\otimes\id_B)\\
&=\Pi_1\circ\m_B^{(3)}\circ(\Pi_1\otimes\Pi_1\otimes\id_B)\,.
\end{split}}
\end{gather}
Furthermore it holds
\begin{equation}\mathlabel{aux3}
(\Pi_1\otimes\id_B)\circ\delta_{(1,0,2)}
= (\Pi_1\otimes\Pi_2)\circ\Delta_B\circ\m_B\circ(\Pi_1\otimes\id_B)
\end{equation}
since the following equations are satisfied.
\begin{equation*}
\begin{split}
\m_B\circ(\Pi_1\otimes\Pi_2)\circ(\m_B\otimes\id_B)\circ(\Pi_1\otimes
 \Delta_B)
&= \m_B\circ(\Pi_1\otimes\id_B)\circ(\m_B\circ(\Pi_1\otimes\Pi_1)\otimes
 \Pi_2)\circ\\
&\quad\circ(\id_B\otimes\Delta_B)\\
&= \m_B\circ(\Pi_1\otimes\m_B\circ(\Pi_1\otimes\Pi_2)\circ\Delta_B)\\
&= \m_B\circ(\Pi_1\otimes\id_B)\,.
\end{split}
\end{equation*}
Then \eqref{aux3} follows from Theorem \ref{proj-cycle}.\ref{proj-cycle2}.
Observe that Corollary \ref{proj-cycle-eqs} implies
the additional conditions
\eqref{techn2-cond} have not been used to derive \eqref{aux3}.
Subsequently we will prove $\m_{B_1\otimes B_2}=\m_{B_1\cpcybi B_2}$
graphically. Henceforth we use the notation
$\Pi_j=$\hbox to 0.5truecm{\hfil$\divide\unitlens by 2
\begin{tangle}\hh\O{\sstyle j}\end{tangle}$\hfil}
for $j\in\{1,2\}$.
\unitlens 17pt                                              %
\begin{equation}\mathlabel{aux7}
\begin{split}
\m_{B_1\otimes B_2}&=
\quad
\divide\unitlens by 2
\vstretch 60\hstretch 75
\begin{tangle}
{\hstr{90}\vstr{80}\O{\sstyle\inj_1}}\step[2]%
{\hstr{90}\vstr{80}\O{\sstyle\inj_2}}\step[4]%
{\hstr{90}\vstr{80}\O{\sstyle\inj_1}}\step[5]%
{\hstr{90}\vstr{80}\O{\sstyle\inj_2}}\\
\id\step\cd\step[2]\cd\step[3]\cd\\
\id\step\O{\sstyle 2}\step[2]\O{\sstyle 2}\step[2]\O{\sstyle 1}%
\step[2]\O{\sstyle 1}\step[3]\O{\sstyle 2}\step[2]\O{\sstyle 2}\\
\id\step\id\step[2]\x\step[2]\id\step[3]\id\step[2]\id\\
\id\step\cu\step[2]\cu\step[3]\id\step[2]\id\\
\id\step[2]\O{\sstyle 1}\step[4]\O{\sstyle 2}\step[4]\id\step[2]\id\\
\cu\step[3]\cd\step[2]\cd\step\id\\
\step\O{\sstyle 1}\step[4]\O{\sstyle 2}\step[2]\O{\sstyle 2}\step[2]%
\O{\sstyle 2}\step[2]\O{\sstyle 1}\step\id\\
\step\id\step[4]\id\step[2]\x\step[2]\id\step\id\\
\step\id\step[4]\cu\step[2]\cu\step\id\\
\step\nw1\step[4]\O{\sstyle 1}\step[4]\O{\sstyle 2}\step[2]\id\\
\step[2]\Cu\step[4]\cu\\
\step[4]{\hstr{90}\vstr{80}\O{\sstyle\proj_1}}\step[7]%
{\hstr{90}\vstr{80}\O{\sstyle\proj_2}}
\end{tangle}
\quad=\quad
\begin{tangle}
{\hstr{90}\vstr{80}\O{\sstyle\inj_1}}\step[2]%
{\hstr{90}\vstr{80}\O{\sstyle\inj_2}}\step[4]%
{\hstr{90}\vstr{80}\O{\sstyle\inj_1}}\step[4]%
{\hstr{90}\vstr{80}\O{\sstyle\inj_2}}\\
\id\step\cd\step[2]\cd\step[2]\cd\\
\id\step\id\step[2]\O{\sstyle 2}\step[2]\id%
\step[2]\id\step[2]\O{\sstyle 2}\step[2]\id\\
\id\step\id\step[2]\x\step[2]\id\step[2]\id\step[2]\id\\
\id\step\cu\step[2]\cu\step[2]\id\step\cd\\
\id\step[2]\O{\sstyle 1}\step[4]\id\step[3]\id\step\id\step[2]\id\\
\cu\step[3]\cd\step[2]\id\step\O{\sstyle 1}\step[2]\O{\sstyle 2}\\
\step\id\step[4]\O{\sstyle 2}\step[2]\O{\sstyle 2}\step[2]\id\step\cu\\
\step\id\step[4]\id\step[2]\x\step[2]\id\\
\step\id\step[4]\cu\step[2]\cu\\
\step\nw1\step[4]\O{\sstyle 1}\step[4]\id\\
\step[2]\Cu\step[4]\id\\
\step[4]{\hstr{90}\vstr{80}\O{\sstyle\proj_1}}\step[6]%
{\hstr{90}\vstr{80}\O{\sstyle\proj_2}}
\end{tangle}
\quad =\quad
\begin{tangle}
{\hstr{90}\vstr{80}\O{\sstyle\inj_1}}\step[2]%
{\hstr{90}\vstr{80}\O{\sstyle\inj_2}}\step[4]%
{\hstr{90}\vstr{80}\O{\sstyle\inj_1}}\step[4]%
{\hstr{90}\vstr{80}\O{\sstyle\inj_2}}\\
\id\step\cd\step[2]\cd\step[2]\cd\\
\id\step\id\step[2]\id\step[2]\id%
\step[2]\id\step[2]\O{\sstyle 2}\step[2]\id\\
\id\step\id\step[2]\x\step[2]\id\step[2]\id\step[2]\id\\
\id\step\cu\step[2]\cu\step[2]\id\step[2]\id\\
\id\step[2]\O{\sstyle 1}\step[4]\id\step[3]\id\step[2]\id\\
\cu\step[3]\cd\step[2]\id\step[2]\id\\
\step\id\step[4]\O{\sstyle 2}\step[2]\id\step[2]%
\id\step[2]\id\\
\step\id\step[4]\id\step[2]\x\step[2]\id\\
\step\id\step[4]\cu\step[2]\cu\\
\step\nw1\step[4]\O{\sstyle 1}\step[4]\id\\
\step[2]\Cu\step[4]\id\\
\step[4]{\hstr{90}\vstr{80}\O{\sstyle\proj_1}}\step[6]%
{\hstr{90}\vstr{80}\O{\sstyle\proj_2}}
\end{tangle}\\
&= \quad
\divide\unitlens by 2
\vstretch 60\hstretch 75
\begin{tangle}
{\hstr{90}\vstr{80}\O{\sstyle\inj_1}}\step[2]%
{\hstr{90}\vstr{80}\O{\sstyle\inj_2}}\step[2]%
{\hstr{90}\vstr{80}\O{\sstyle\inj_1}}\step[3]%
{\hstr{90}\vstr{80}\O{\sstyle\inj_2}}\\
\id\step[2]\cu\step[3]\id\\
\id\step[2]\cd\step[2]\cd\\
\id\step[2]\id\step[2]\id\step[2]\O{\sstyle2}\step[2]\id\\
\id\step[2]\id\step[2]\x\step[2]\id\\
\nw1\step\cu\step[2]\cu\\
\step\cu\step[4]\id\\
\step[2]{\hstr{90}\vstr{80}\O{\sstyle\proj_1}}\step[5]%
{\hstr{90}\vstr{80}\O{\sstyle\proj_2}}
\end{tangle}
\quad = \quad
\begin{tangle}
{\hstr{90}\vstr{80}\O{\sstyle\inj_1}}\step[2]%
{\hstr{90}\vstr{80}\O{\sstyle\inj_2}}\step[2]%
{\hstr{90}\vstr{80}\O{\sstyle\inj_1}}\step[2]%
{\hstr{90}\vstr{80}\O{\sstyle\inj_2}}\\
\id\step[2]\cu\step\dd\\
\id\step[3]\cu\\
\nw1\step[2]\cd\\
\step\cu\step[2]\d\\
\step[2]{\hstr{90}\vstr{80}\O{\sstyle\proj_1}}\step[4]%
{\hstr{90}\vstr{80}\O{\sstyle\proj_2}}
\end{tangle}
\quad = (\proj_1\otimes\proj_2)\circ\Delta_B\circ\m_B^{(4)}\circ
(\inj_1\otimes\inj_2\otimes\inj_1\otimes\inj_2)\\
&=\Lambda^{-1}\circ\m_B\circ(\Lambda\otimes\Lambda)=
\m_{B_1\cpcybi B_2}\,.
\end{split}
\end{equation}
\unitlens 15pt                                              %
To derive the first identity of \eqref{aux7} we used the specific form of
the structure morphisms $\tilde{\varphi}_{1,2}$, $\tilde{\varphi}_{2,1}$,
$\hat{\tilde\sigma}$, and $\hat{\tilde\rho}$.
Then \eqref{pi-rel1}, \eqref{cons-proj-cycle1}, and the third
identity of \eqref{techn2-cond} for $j=1$ yield the second equality
of \eqref{aux7}. With Theorem \ref{proj-cycle}.\ref{proj-cycle1},
\eqref{cons-proj-cycle1} and
\eqref{aux1} we derive the third equation, whereas \eqref{aux2} and
again use of Theorem \ref{proj-cycle}.\ref{proj-cycle1} yield the fourth
identity of \eqref{aux7}. The fifth equality comes from application of the
second ``projection" relation of \eqref{techn2-cond}, and for the derivation
of the sixth identity we used \eqref{aux3}. Finally the definition of
$\Lambda$, given in the proof of Theorem \ref{proj-cycle}, yields the
result. In a $\pi$-symmetric way the identity
$\Delta_{B_1\otimes B_2}=\Delta_{B_1\cpcybi B_2}$ will be shown.
The (co\n-)unital identities \eqref{type-alpha1} can be verified
straightforwardly from the definitions, the assumptions and
Corollary \ref{proj-cycle-eqs}. It remains to prove the ``projection''
relations \eqref{type-alpha2}. Observe that
\begin{equation}\mathlabel{aux4}
\Pi_1=\Lambda\circ(\id_{B_1}\otimes\eta_2\circ\varepsilon_2)
\circ\Lambda^{-1}\,,\quad
\Pi_2=\Lambda\circ(\eta_1\circ\varepsilon_1\otimes\id_{B_2})\circ
\Lambda^{-1}\,.
\end{equation}
Taking into account the relation
$\Delta_B=(\Lambda\otimes\Lambda)\circ\Delta_{B_1\cpcybi B_2}\circ
\Lambda^{-1}$ and the relation $\m_B=\Lambda\circ\m_{B_1\cpcybi B_2}\circ
(\Lambda^{-1}\otimes\Lambda^{-1})$ one obtains with the help of
\eqref{aux4} and \eqref{type-alpha1}
\begin{equation}\mathlabel{aux5}
\begin{split}
\Delta_B\circ\Pi_1 &= (\Lambda\otimes\Lambda\circ(\id_{B_1}\otimes\eta_2))
 \circ(\id_{B_1}\otimes\nu_l)\circ\Delta_1\circ(\id_{B_1}\otimes
 \varepsilon_2)\circ\Lambda\,,\\
\Pi_1\circ\m_B &=\Lambda\circ(\id\otimes\eta_2)\circ\m_1\circ
(\m_1\otimes\sigma)\circ(\id\otimes\varphi_{2,1}\otimes\id)\circ
(\Lambda^{-1}\otimes\Lambda^{-1})\,.
\end{split}
\end{equation}
Gluing the identities of \eqref{aux5} and $\Pi_1$ according to the
left and right hand side of the third equation of \eqref{techn2-cond}
(for $j=1$), using \eqref{aux4} again, and eventually multiplying both
resulting sides with
$\vert\circ(\Lambda\otimes\Lambda)\circ(\eta_1\otimes\id\otimes\id\otimes
\eta_2)$ and $(\id\otimes\varepsilon_2\otimes\id\otimes\varepsilon_2)\circ
(\Lambda^{-1}\otimes\Lambda^{-1})\circ\vert$ yields
$(\mu_l\otimes\id_{B_1})\circ(\id_{B_2}\otimes\Delta_1)=
(\m_1\circ(\id_{B_1}\otimes\sigma)\otimes\id_{B_1})\circ
(\varphi_{2,1}\otimes\nu_l)\circ(\id_{B_2}\otimes\Delta_1)$
which is the second identity of \eqref{type-alpha2}. Analogously
the first equation of \eqref{type-alpha2} can be derived.
This proves that $B_1\cpcybi B_2$ is cocycle cross product
bialgebra.

Conversely suppose that $\Lambda:B_1\cpcybi[\mathfrak{c}] B_2\to B$ is an isomorphism
of bialgebras. To prove statement $2.$ it suffices to verify relations
\eqref{techn2-cond}. Like in the proof
``\ref{proj-cycle3}.$\Rightarrow$\ref{proj-cycle2}.'' of Theorem
\ref{proj-cycle} we define
$\inj_1 :=\Lambda\circ(\id_{B_1}\otimes\eta_2)$,
$\inj_2 :=\Lambda\circ(\eta_1\otimes\id_{B_2})$,
$\proj_1 :=(\id_{B_1}\otimes\varepsilon_2)\circ\Lambda^{-1}$, and
$\proj_2 :=(\varepsilon_1\otimes\id_{B_2})\circ\Lambda^{-1}$.
Then the identities \eqref{aux5} will be proven for
$\Pi_j:=\inj_j\circ\proj_j$ in the way described above. And the
previously performed gluing of the identities \eqref{aux5} yields
\begin{equation*}
\begin{split}
(\Pi_1\otimes\id)\circ\delta_{(0,0,0)}\circ(\id\otimes\Pi_1) &=
(\Lambda\circ(\id\otimes\eta_2)\otimes\Lambda\circ(\id\otimes\eta_2))\circ
 (\m_1\circ(\m_1\otimes\sigma)\otimes\id)\circ\\
&\quad\circ(\id\otimes\varphi_{2,1}\otimes\nu_l)\circ
(\Lambda^{-1}\otimes\Delta_1\circ(\id\otimes\varepsilon_2)\circ\Lambda^{-1})\,,
\\
(\Pi_1\otimes\id)\circ\delta_{(0,1,0)}\circ(\id\otimes\Pi_1) &=
(\Lambda\circ(\id\otimes\eta_2)\otimes\Lambda\circ(\id\otimes\eta_2))\circ
 (\m_1\circ(\id\otimes\mu_l)\otimes\id)\circ\\
&\quad\circ(\Lambda^{-1}\otimes\Delta_1\circ
(\id\otimes\varepsilon_2)\circ\Lambda^{-1})\,.
\end{split}
\end{equation*}
Therefore by assumption \eqref{type-alpha2} the third identity
of \eqref{techn2-cond} for $j=1$ follows. Applying $\pi$-symmetry
the third equation of \eqref{techn2-cond} for $j=2$ will be derived.
Then all conditions are satisfied which have been needed to
derive the first four identities in \eqref{aux7}. Hence
\begin{equation*}
\begin{split}
\m_{B_1\cpcybi[\mathfrak{c}] B_2} &=
(\proj_1\circ\m_B^{(3)}\otimes\proj_1\circ\m_B)\circ
(\id^{(2)}\otimes\Psi_{B,B}\circ(\id\otimes\Pi_2)\otimes\id)\circ\\
&\quad\circ(\inj_1\otimes\Delta_B\circ\m_B\circ(\inj_2\otimes\inj_1)\otimes
\Delta_B\circ\inj_2)
\end{split}
\end{equation*}
On the other hand it holds by assumption
$\m_{B_1\cpcybi[\mathfrak{c}] B_2}
=\Lambda^{-1}\circ\m_B\circ(\Lambda\otimes\Lambda)
= (\proj_1\otimes\proj_2)\circ\Delta_B\circ\m_B^{(4)}\circ
 (\inj_1\otimes\inj_2\otimes\inj_1\otimes\inj_2)
= (\proj_1\otimes\proj_2)\circ\delta_{(0,0,0)}\circ
 (\inj_1\otimes\m_B^{(3)}\circ(\inj_2\otimes\inj_1\otimes\inj_2))$,
where we used \eqref{aux3} which has been derived under the more
general assumption of Theorem \ref{proj-cycle}. Multiplying both
expressions of $\m_{B_1\cpcybi[\mathfrak{c}] B_2}$ with
$\vert\circ(\eta_1\otimes\proj_2\otimes\proj_1\otimes\proj_2)$ and
$(\inj_1\otimes\inj_2)\circ\vert$, and using that $\inj_1$ is algebra
morphism yields
\begin{equation*}
\begin{split}
&(\Pi_1\circ\m_B\otimes\Pi_1\circ\m_B)\circ
(\id\otimes\Psi_{B,B}\circ(\id\otimes\Pi_2)\otimes\id)\circ
(\Delta_B\circ\m_B\circ(\Pi_2\otimes\Pi_1)\otimes\Delta_B\circ\Pi_2)\\
&=(\Pi_1\otimes\Pi_1)\circ\Delta_B\circ\m_B^{(3)}\circ
(\Pi_2\otimes\Pi_1\otimes\Pi_2)
\end{split}
\end{equation*}
from which the second identity of \eqref{techn2-cond} can be derived
easily with the help of Theorem \ref{proj-cycle}.\ref{proj-cycle1}.
Similarly by $\pi$-symmetry the first equation of \eqref{techn2-cond}
will be proven.\end{proof}
\abs
\begin{remark}\Label{yu-cond}
{\normalfont The ``projection" relations $(\gamma_1)$ and
$(\gamma_2)$ in \eqref{techn2-cond} can be derived from the
identities
\begin{equation}\mathlabel{alt-gamma}
(\mathbf{\beta_1}):\ \beta_1=\beta_1\circ(\id\otimes\Pi_1)
\quad\text{and}\quad
(\mathbf{\beta_2}):\ \beta_2=(\Pi_2\otimes\id)\circ\beta_2
\end{equation}
where $\beta_1:=(\Pi_1\circ\m_B\otimes\id)\circ
(\id\otimes\Psi_{B,B})\circ(\Delta_B\circ\Pi_1\otimes\id)$ and
$\beta_2:=(\id\otimes\Pi_2\circ\m_B)\circ
(\Psi_{B,B}\otimes\id)\circ(\id\otimes\Delta_B\circ\Pi_2)$.
However, in general the conditions $(\mathbf{\beta_1})$ and 
$(\mathbf{\beta_2})$ are not equivalent to the conditions
$(\mathbf{\gamma_1})$ and $(\mathbf{\gamma_2})$ in 
Theorem \ref{bialg-cocycle2}.}
\end{remark}

\begin{proposition}\Label{alg-coalg-triv}
Under the equivalent conditions of Theorem \ref{bialg-cocycle2}
the following statements hold.
\begin{enumerate}
\item
$\mu_l=\varepsilon_2\otimes\id_{B_1}$ is trivial $\Leftrightarrow$
$\Pi_1\circ\m_B\circ(\id_B\otimes\Pi_1)=
\Pi_1\circ\m_B\circ(\Pi_1\otimes\Pi_1)$.
\item
$\mu_r=\id_{B_2}\otimes\varepsilon_1$ is trivial
$\Leftrightarrow$ $\Pi_2\circ\m_B=\Pi_2\circ\m_B\circ(\Pi_2\otimes\Pi_2)$
$\Leftrightarrow$ $\proj_2$ is algebra morphism.
\item
$\nu_l=\eta_2\otimes\id_{B_1}$ is trivial $\Leftrightarrow$
$\Delta_B\circ\Pi_1=(\Pi_1\otimes\Pi_1)\circ\Delta_B\circ\Pi_1$
$\Leftrightarrow$ $\inj_1$ is coalgebra morphism.
\item
$\nu_r=\id_{B_2}\otimes\eta_1$ is trivial $\Leftrightarrow$
$(\Pi_2\otimes\id_B)\circ\Delta_B\circ\Pi_2=
(\Pi_2\otimes\Pi_2)\circ\Delta_B\circ\Pi_2$. 
\item
$\sigma=\eta_1\circ(\varepsilon_2\otimes\varepsilon_2)$ is trivial 
$\Leftrightarrow$ $\hat\sigma=\eta_1\otimes\m_2$ $\Leftrightarrow$
$\m_B\circ(\Pi_2\otimes\Pi_2)=\Pi_2\circ\m_B\circ(\Pi_2\otimes\Pi_2)$
$\Leftrightarrow$ $\inj_2$ is algebra morphism.
\item
$\rho= (\eta_1\otimes\eta_1)\circ\varepsilon_2$ is trivial
$\Leftrightarrow$ $\hat\rho=\Delta_1\otimes\varepsilon_2$ $\Leftrightarrow$
$(\Pi_1\otimes\Pi_1)\circ\Delta_B=(\Pi_1\otimes\Pi_1)\circ\Delta_B\circ\Pi_1$
$\Leftrightarrow$ $\proj_1$ is coalgebra morphism.
\item
$\mu_l$ and $\sigma$ are trivial
$\Leftrightarrow$ $\Pi_1\circ\m_B=\Pi_1\circ\m_B\circ(\Pi_1\otimes\Pi_1)$
$\Leftrightarrow$ $\proj_1$ is algebra morphism.
\item
$\nu_r$ and $\rho$ are trivial $\Leftrightarrow$
$\Delta_B\circ\Pi_2=(\Pi_2\otimes\Pi_2)\circ\Delta_B\circ\Pi_2$
$\Leftrightarrow$ $\inj_2$ is coalgebra morphism. 
\end{enumerate}
\end{proposition}

\begin{proof}
We prove Proposition \ref{alg-coalg-triv}.1 and \ref{alg-coalg-triv}.5.
The remaining statements can be derived in a similar
manner or follow directly by $\pi$-symmetric reasoning.
Without loss of generality we may assume 
$\Pi_1=\id_{B_1}\otimes\eta_2\circ\varepsilon_2$ and
$\Pi_2=\eta_1\circ\varepsilon_1\otimes\id_{B_2}$.

Ad 1.: Suppose that $\mu_l$ is trivial. This means that
$(\id_{B_1}\otimes\varepsilon_2)\circ\varphi_{2,1}=
\varepsilon_2\otimes\id_{B_1}$. Using \eqref{cond2-caat} and the unital
identities of \eqref{cond3-caat} then yields
$\Pi_1\circ\m_B\circ(\id_B\otimes\Pi_1)=
(\m_{B_1}\otimes\eta_2\circ\varepsilon_2)\circ
(\id_{B_1}\otimes\varphi_{2,1}\otimes\varepsilon_2)=
(\m_{B_1}\otimes\eta_2)\circ(\id_{B_1}\otimes\mu_l\otimes\varepsilon_2)=
(\m_{B_1}\otimes\eta_2)\circ(\id_{B_1}\otimes\varepsilon_2\otimes
\id_{B_1}\otimes\varepsilon_2)$. On the other hand from the unital identities of
\eqref{cond3-caat} immediately follows
$\Pi_1\circ\m_B\circ(\Pi_1\otimes\Pi_1)=(\m_{B_1}\otimes\eta_2)\circ
(\id_{B_1}\otimes\varepsilon_2\otimes\id_{B_1}\otimes\varepsilon_2)$.
Conversely if $\Pi_1\circ\m_B\circ(\id_B\otimes\Pi_1)=
\Pi_1\circ\m_B\circ(\Pi_1\otimes\Pi_1)$ holds then analogous calculations as
before yield $(\m_{B_1}\otimes\eta_2\circ\varepsilon_2)\circ
(\id_{B_1}\otimes\varphi_{2,1}\otimes\varepsilon_2)=(\m_{B_1}\otimes\eta_2)\circ
(\id_{B_1}\otimes\varepsilon_2\otimes\id_{B_1}\otimes\varepsilon_2)$. And then
the triviality of $\mu_l=(\id_{B_1}\otimes\varepsilon_2)\circ\varphi_{2,1})$
follows straightforwardly.

Ad 5.: Suppose that $\m_B\circ(\Pi_2\otimes\Pi_2)=\Pi_2\circ
\m_B\circ(\Pi_2\otimes\Pi_2)$. Then $\m_2=\proj_2\circ\m_B\circ
(\inj_2\otimes\inj_2)$ is associative for as $\m_2\circ(\id\otimes\m_2)=
\proj_2\circ\m_B\circ(\inj_2\otimes\Pi_2\circ\m_B\circ(\inj_2\otimes\inj_2))=
\proj_2\circ\m_B\circ(\id_B\otimes\m_B)\circ(\inj_2\otimes\inj_2\otimes\inj_2)
=\proj_2\circ\m_B\circ(\m_B\otimes\id_B)\circ(\inj_2\otimes\inj_2\otimes\inj_2)
=\m_2\circ(\m_2\otimes\id_{B_2})$. Hence $B_2$ is an algebra. And 
$\inj_2\circ\m_2=\Pi_2\circ\m_B\circ(\inj_2\otimes\inj_2)=
\m_B\circ(\inj_2\otimes\inj_2)$. Therefore $\inj_2$ is algebra morphism.
Conversely, if $\inj_2$ is algebra morphism then
$\m_B\circ(\Pi_2\otimes\Pi_2)=\inj_2\circ\m_2\circ(\proj_2\otimes\proj_2)
=\inj_2\circ\proj_2\circ\inj_2\circ\m_2\circ(\proj_2\otimes\proj_2)
=\Pi_2\circ\m_B\circ(\Pi_2\otimes\Pi_2)$.

Let $\inj_2$ be algebra morphism. Then from the last identity of
\eqref{struc-morph} and the second identity of \eqref{cons-proj-cycle1}
we derive $\hat\sigma=(\proj_1\otimes\proj_2)\circ\Delta_B\circ\m_B\circ
(\inj_2\otimes\inj_2)=(\proj_1\otimes\proj_2)\circ\Delta_B\circ\inj_2\circ\m_2
=\eta_1\otimes\m_2$. If on the other hand $\hat\sigma=\eta_1\otimes\m_2$ then
we use \eqref{cond2-caat} and \eqref{cond3-caat} to obtain
$\m_B\circ(\inj_2\otimes\inj_2)=\m_B\circ
(\eta_1\otimes\id_{B_2}\otimes\eta_1\otimes\id_{B_2})=
(\eta_1\otimes\id_{B_2})\circ\m_2=\inj_2\circ\m_2$ which shows that $\inj_2$ is
algebra morphism.

It is an easy exercise to prove that triviality of
$\hat\sigma$ and triviality of $\sigma$ are equivalent.\end{proof}

\begin{remark}{\normalfont
Under the equivalent conditions of Theorem \ref{proj-cycle} 
similar results like in Proposition \ref{alg-coalg-triv} can be shown
for general cross product bialgebras if $\m_2$, $\Delta_1$, $\mu_l$,
$\mu_r$, $\nu_l$ and $\nu_r$ will be defined formally as in
\eqref{proj-struc}.}
\end{remark}

\subsection{Strong Cross Product Bialgebras}

In the following we define strong cross product
bialgebras. They will be studied in more detail in Theorems \ref{hp-cp}
and \ref{central-theory} in connection with
so-called strong Hopf data. It turns out that strong
cross product bialgebras are the central objects in our universal,
(co\n-)modular co-cyclic theory of cross product bialgebras
\footnote{strong cross product bialgebras and
strong Hopf data provide a unifying universal and (co\n-)modular
co-cyclic theory of cross product bialgebras. But
they are probably not the most general setting to meet
the same demands. Therefore our notation ``strong'' should be
specified further. However, in order to avoid terminological
blow-up we will henceforth use our notation, always having
in mind that there might be a weaker definition of ``strong''.
This is certainly an interesting direction for further study.}.

\begin{definition}\Label{strong-type-alpha}
A cocycle cross product bialgebra $B_1\cpcybi[\mathfrak{c}] B_2$
is called strong cross product bialgebra if in addition the identities
\begin{gather}\mathlabel{strong-bi1}
\divide\unitlens by 3
\multiply\unitlens by 2
\begin{tangle}
	\hh\hld\step\id\\
	\hh\id\hstep\x\\
	\hh\hru\step\id\\
\end{tangle}
=
\begin{tangle}
\ld\step\counit
\end{tangle}
\quad ,
\quad
\begin{tangle}
	\hh\id\step\hld\\
	\hh\x\hstep\id\\
	\hh\id\step\hru\\
\end{tangle}
=
\begin{tangle}\unit\step\ru
\end{tangle}
\\
\mathlabel{strong-bi2}
\divide\unitlens by 2
\vstretch 60
\begin{tangle}
\Ld\step\id\\
\rd\step\hx\\
\id\step\hx\step\id\\
\hcu*\step\id\step\id
\end{tangle}
=
\begin{tangle}
\unit\step\unit\step\vstr{40}\id\step\id\\
\vstr{40}\id\step\id\step\id\step\counit
\end{tangle}
\quad ,
\quad
\begin{tangle}
\id\step\id\step\hcd*\\
\id\step\hx\step\id\\
\hx\step\lu\\
\id\step\Ru
\end{tangle}
=
\begin{tangle}
\unit\vstr{40}\step\id\step\id\step\id\\
\vstr{40}\id\step\id\step\counit\step\counit
\end{tangle}
\divide\unitlens by 3
\multiply\unitlens by 4
\vstretch 100
\quad ,
\quad
\vstretch 120
\begin{tangle}
\hh\id\step\hld\\
\hh\cu*\hstep\id
\end{tangle}
=
\vstretch 90
\begin{tangle}
\hh\step\id\\
\hh\counit\step\id\\
\hh\unit\step\id
\end{tangle}
\quad ,
\quad
\vstretch 120
\begin{tangle}
\hh\id\hstep\cd*\\
\hh\hru\step\id
\end{tangle}
=
\vstretch 90
\begin{tangle}
\hh\id\\
\hh\id\step\counit\\
\hh\id\step\unit
\end{tangle}
\\[3pt]
\mathlabel{strong-bi3}
\divide\unitlens by 3
\def\FillCircDiam{6}
\vstretch 50
\begin{tangle}
\vstr{70}\id\Step\id\Step\cd*\\
\id\Step\id\step\cd\step\id\\
\id\Step\hx\Step\id\step\id\\
\id\step\dd\ld\step\dd\step\id\\
\hx\step\id\step\hx\Step\id\\
\id\step\d\lu\step\id\Step\id\\
\id\Step\ru\step\id\Step\id
\end{tangle}
=
\vstretch 100
\begin{tangle}
\vstr{70}\id\step\id\step[1.5]\cd*\\
\hh\id\step\id\step\cd\step[1.5]\id\\
\hh\id\step\x\step\id\step[1.5]\id\\
\hh\x\step\x\step[1.5]\id\\
\hh\id\step{\hstr{200}\hru}\step\id\step[1.5]\id
\end{tangle}
\quad,\quad
\vstretch 50
\begin{tangle}
\id\Step\id\step\ld\Step\id\\
\id\Step\id\step\rd\d\step\id\\
\id\Step\hx\step\id\step\hx\\
\id\step\dd\step\ru\dd\step\id\\
\id\step\id\Step\hx\Step\id\\
\id\step\cu\step\id\Step\id\\
\vstr{70}\cu*\Step\id\Step\id
\end{tangle}
=
\vstretch 100
\begin{tangle}
\hh\id\step[1.5]\id\step{\hstr{200}\hld}\step\id\\
\hh\id\step[1.5]\x\step\x\\
\hh\id\step[1.5]\id\step\x\step\id\\
\hh\id\step[1.5]\cu\step\id\step\id\\
\vstr{70}\cu*\step[1.5]\id\step\id
\end{tangle}
\end{gather}
are fullfilled. We denote strong cross product bialgebras
by$\,$\footnote{The deeper meaning of this notation will become clear 
in Section \ref{hopf-data}.} $B_1\overset{\mathfrak{h}}\bowtie B_2$.
\end{definition}

\begin{remark}[!]\Label{non-except-bialg}
{\normalfont A tedious calculation shows that the ``projection'' relations
$(\mathbf{\gamma_1})$, $(\mathbf{\gamma_2})$, $(\mathbf{\delta_1})$, 
and $(\mathbf{\delta_2})$ in \eqref{techn2-cond} as well as
\eqref{type-alpha2} are redundant for strong cross product
bialgebras. Therefore, if we use the morphisms
$\m_{i,jk}$ and $\Delta_{ij,k}$ from
Corollary \ref{co-act-inv} an alternative
definition of strong cross product bialgebras can be given as follows.
\abs
\emph{A strong cross product bialgebra is a cross product bialgebra
for which the identities \eqref{strong-bi1}, \eqref{strong-bi2} and
\eqref{strong-bi3} formally hold for the corresponding morphisms 
$\m_{i,jk}$ and $\Delta_{ij,k}$.}
\abs
We call a strong cross product bialgebra \textit{regular} if all
defining identities have rank 
$(2,2)$, $(1,3)$ or $(3,1)$.
We call it \textit{pure} if \eqref{strong-bi1}, \eqref{strong-bi2} 
and \eqref{strong-bi3} are redundant.}
\end{remark}
\abs
Any bialgebra $B$ which is isomorphic to a certain cocycle
cross product
bialgebra has injections and projections $\inj_1$, $\inj_2$, $\proj_1$,
and $\proj_2$ which are uniquely determined by the cocycle
cross product bialgebra and the given isomorphism
(see Theorem \ref{bialg-cocycle2} and Remark \ref{alpha-ip}).
In the special case
of strong cross product bialgebras the structure
morphisms obey the additional identities
\eqref{strong-bi1} -- \eqref{strong-bi3}.
Using \eqref{struc-morph} and \eqref{proj-struc} these relations can be
translated easily into equations of the idempotents $\Pi_1$ and $\Pi_2$
of the bialgebra $B$ as follows.
\begin{equation}\mathlabel{pi-strong}
\begin{gathered}
(\Pi_2\otimes\id_B)\circ\Phi_{\id_B}\circ(\Pi_1\otimes\Pi_1)=
(\Pi_2\otimes\id_B)\circ\Delta_B\circ\Pi_1\otimes\varepsilon_B\,,\\[5pt]
((\Pi_1\otimes\Pi_1)\circ\Phi_{\Pi_2}\otimes\id_B)\circ
(\Pi_2\otimes\Psi_{B,B})\circ(\Delta_B\circ\Pi_1\otimes\Pi_2)
=\eta_B\otimes\eta_B\otimes\Pi_1\otimes\varepsilon_B\,,\\[5pt]
(\Pi_1\otimes\id_B)\circ\delta_{(2,2,0)}\circ(\id_B\otimes\Pi_1)=
\eta_B\circ\varepsilon_B\otimes\Pi_1\,,\\[5pt]
\begin{split}
&(\Pi_1\otimes\Pi_2\circ\m_B\circ(\id_B\otimes\Pi_2)\otimes\id_B\otimes\Pi_1)
 \circ(\Psi_{B,B}\otimes\Phi_{\Pi_2}\otimes\id_B)\circ\\
&\circ(\id_B\otimes(\Psi_{B,B}\otimes\id_B\otimes\id_B)\circ
 (\Pi_1\otimes(\Delta_B\circ\Pi_1\otimes\id_B)\circ\Delta_B\circ\Pi_2))\\
 =\,\,&((\Pi_1\otimes\Pi_2\circ\m_B)\circ(\Psi_{B,B}\otimes\id_B)\circ
 (\id_B\otimes\Delta_B)\otimes\Pi_1\otimes\Pi_1)\circ\\
&\circ(\id_B\otimes(\Psi_{B,B}\otimes\id_B)\circ(\id_B\otimes(\Pi_1\otimes\id_B)
 \circ\Delta_B\circ\Pi_2))
\end{split}\\[5pt]
\text{and the corresponding $\pi$-symmetric counterparts.}
\end{gathered}
\end{equation}
We used $\Phi_f:=(\m_B\otimes\id_B)\circ(f\otimes\Psi_{B,B})\circ
(\Delta_B\otimes\id_B)$ for $f:B\to B$, and $\delta_{(i,j,k)}$ of
Theorem \ref{bialg-cocycle2}.
\abs
Recall Remark \ref{pi-sym} and observe that our construction of
(strong/cocycle) cross product bialgebras is invariant under
$\pi$-symmetry. However, cross product bialgebras in general fail to be
invariant if duality or rotational symmetry along a vertical axis
(in the plane of the paper) will be transformed separately. Without
problems the dual versions of cross product bialgebras can be defined using
our original definition. The corresponding results follow immediately.
If the category $\C$ has (right) duality then the following proposition shows
how dual cross product bialgebras can be constructed explicitely.

\begin{proposition}
Suppose that $\C$ is a braided category with duality functor
$(\_)^\vee:\C\to\C^{\mathrm{op}}_{\mathrm{op}}$ where
$\C^{\mathrm{op}}_{\mathrm{op}}$ is the opposite category with opposite
tensor product. If $(B;B_1,B_2,\inj_1,\inj_2,\proj_1,\proj_2)$ is a
(strong/cocycle) cross product bialgebra in $\C$ then
the dual tuple $(B^\vee;B_1^\vee,B_2^\vee,\inj_1^\vee,\inj_2^\vee,
\proj_1^\vee,\proj_2^\vee)$ is a dual (strong/cocycle)
cross product bialgebra in $\C^{\mathrm{op}}_{\mathrm{op}}$.
\phantom{xxxxx}\endproof
\end{proposition}

\section{Hopf Data}\Label{hopf-data}

In Section \ref{sec-ccpb} we studied cross product bialgebras
from a universal point of view. We answered the question under which
conditions a bialgebra is isomorphic to a cross product bialgebra.
Now we present an explicit (co\n-)-modular co-cyclic
construction method in terms of Hopf data. 
A Hopf datum consists of two objects with
certain interrelated (co\n-)modular co-cyclic identities. In Theorem
\ref{hp-cp} we will show that so-called strong Hopf data and
strong cross product bialgebras are different descriptions of
the same objects. Universality of our ``strong'' construction will be
demonstrated in Theorem \ref{central-theory}.
We postpone the lengthy proof of Theorem
\ref{hp-cp} to Section \ref{proof-hp-cp}.
In the subsequent definition of Hopf data occur two objects
$B_1$ and $B_2$, and morphisms
\begin{equation}\mathlabel{hp-morph}
\begin{split}
\m_1 &:B_1\otimes B_1\to B_1\,,\\
\Delta_1 &:B_1\to B_1\otimes B_1\,,\\
\eta_1 &:\E\to B_1\,,\\
\varepsilon_1 &:B_1\to\E\,,\\
\mu_l &:B_2\otimes B_1\to B_1\,,\\
\mu_r &:B_2\otimes B_1\to B_2\,,\\
\sigma &:B_2\otimes B_2\to B_1\,,
\end{split}
\qquad
\begin{split}
\m_2 &:B_2\otimes B_2\to B_2\,,\\
\Delta_2 &:B_2\to B_2\otimes B_2\,,\\
\eta_2 &:\E\to B_2\,,\\
\varepsilon_2 &:B_2\to\E\,,\\
\nu_l &:B_1\to B_2\otimes B_1\,,\\
\nu_r &:B_2\to B_2\otimes B_1\,,\\
\rho &: B_2\to B_1\otimes B_1\,.
\end{split}
\end{equation}
Again we use the graphical presentation of Figure \ref{fig-conv} for these
morphisms, and we represent $\sigma$ and $\rho$
by $\sigma=\begin{tangle}\hh \cu* \end{tangle}:B_2\otimes B_2\to B_1$ and
$\rho=\begin{tangle}\hh \cd* \end{tangle}:B_2\to B_1\otimes B_1$.
Similarly as in \eqref{phi-sigma-rho} we define morphisms
$\varphi_{1,2}$, $\varphi_{2,1}$, $\hat\sigma$, and $\hat\rho$ by
\begin{equation*}
\begin{split}
\varphi_{1,2}&:=(\m_2\otimes\m_1)\circ(\id_{B_1}\otimes\Psi_{B_1,B_2}\otimes
\id_{B_2})\circ(\nu_l\otimes\nu_r)\,,\\
\varphi_{2,1}&:=(\mu_l\otimes\mu_r)\circ(\id_{B_2}\otimes\Psi_{B_2,B_1}
\otimes\id_{B_1})\circ(\Delta_2\otimes\Delta_1)\,,\\
\hat\sigma&:=(\sigma\otimes\m_2\circ(\mu_r\otimes\id_{B_2}))\circ
(\id_{B_2}\otimes\Psi_{B_2,B_2}\otimes\id_{B_1\otimes B_2})\circ
(\Delta_2\otimes(\nu_r\otimes\id_{B_2})\circ\Delta_2)\,,\\
\hat\rho&:=(\m_1\circ(\id_{B_1}\otimes\mu_l)\otimes\m_1)\circ
(\id_{B_1\otimes B_2}\otimes\Psi_{B_1,B_1}\otimes\id_{B_1})\circ
((\id_{B_1}\otimes\nu_l)\circ\Delta_1\otimes\rho)
\end{split}
\end{equation*}
or graphically
\begin{equation}\mathlabel{phi-sigma-rho2}
\begin{array}{c}
\hstretch 125
\vstretch 85
\varphi_{1,2}
:=
\divide\unitlens by 3
\begin{tangle}
\hh\step\id\step\id\\
\hh\ld\step\rd\step\\
\id\step\hx\step\id\\
\hh\cu\step\cu\\
\hh\hstep\id\Step\id
\end{tangle}
\multiply\unitlens by 3
\quad ,\quad
\varphi_{2,1}
:=
\divide\unitlens by 3
\begin{tangle}
\hh\hstep\id\Step\id\\
\hh\cd\step\cd\\
\id\step\hx\step\id\\
\hh\lu\step\ru\\
\hh\step\id\step\id
\end{tangle}
\multiply\unitlens by 3
\quad ,\quad
\hat\rho
:=
\hstretch 150
\vstretch 85
\divide\unitlens by 3
\begin{tangle}
\hh\hstep\id\Step\id\\
\hh\cd\step\cd*\\
\hh \id\hstep\hld\step\id\step\id\\
\id\hstep\id\hstep\hx\step\id\\
\hh \id\hstep\hlu\step\id\step\id\\
\hh \cu\step\cu\\
\hh\hstep\id\Step\id
\end{tangle}
\multiply\unitlens by 3
\quad ,\quad
\hat\sigma
:=
\divide\unitlens by 3
\begin{tangle}
\hh\hstep\id\Step\id\\
\hh\cd\step\cd \\
\hh\id\step\id\step\hrd\hstep\id \\
\id\step\hx\hstep\id\hstep\id \\
\hh\id\step\id\step\hru\hstep\id \\
\hh\cu*\step\cu\\
\hh\hstep\id\Step\id
\end{tangle}
\multiply\unitlens by 3
\end{array}
\end{equation}
Occasionally we also use the graphical abbreviations
\begin{equation}\mathlabel{phi-sigma-rho3}
\begin{array}{c}
\varphi_{1,2}=
\divide\unitlens by 2
\begin{tangle}\ox{12}\end{tangle}
\multiply\unitlens by 2
\quad ,\quad
\varphi_{2,1}=
\divide\unitlens by 2
\begin{tangle}\ox{21}\end{tangle}
\multiply\unitlens by 2
\quad ,\quad
\hat\rho
=
\divide\unitlens by 2
\begin{tangle}
\hh\hstep\id\step\id\\
\hh\ffbox 2{\hat\rho}\\
\hh\hstep\id\step\id
\end{tangle}
\multiply\unitlens by 2
\quad ,\quad
\hat\sigma=
\divide\unitlens by 2
\begin{tangle}
\hh\hstep\id\step\id\\
\hh\ffbox 2{\hat\sigma}\\
\hh\hstep\id\step\id
\end{tangle}
\multiply\unitlens by 2
\end{array}
\end{equation}

\begin{definition}\Label{hp}
The tuple $\mathfrak{h}=
\big((B_1,\m_1,\eta_1,\Delta_1,\varepsilon_1),(B_2,\m_2,\eta_2,\Delta_2,
\varepsilon_2);\mu_l,\mu_r,\nu_l,\nu_r,\rho,\sigma\big)$ is called Hopf datum if
\begin{enumerate}
\item
$(B_1,\m_1,\eta_1)$ is an algebra and $\varepsilon_1:B_1\to \E$
is an algebra morphism.
\item
$(B_2,\Delta_2,\varepsilon_2)$ is a coalgebra and $\eta_2:\E\to B_2$
is a coalgebra morphism.
\item
$(B_1,\nu_l)$ is left $B_2$-comodule.
\item
$(B_2,\mu_r)$ is right $B_1$-module.
\item The identities
\begin{equation}\mathlabel{hp1}
\begin{gathered}
\Delta_1\circ\eta_1=\eta_1\otimes\eta_2\,,\quad
(\id_{B_1}\otimes\varepsilon_1)\circ\Delta_1=
(\varepsilon_1\otimes\id_{B_1})\circ\Delta_1=\id_{B_1}\,,\\
\varepsilon_2\circ\m_2=\varepsilon_2\otimes\varepsilon_2\,,\quad
\m_2\circ(\id_{B_2}\otimes\eta_2)=\m_2\circ(\eta_2\otimes\id_{B_2})=
\id_{B_2}\,,\\[5pt]
\mu_r\circ(\eta_2\otimes\id_{B_1})=\eta_2\circ\varepsilon_1\,,\quad
\varepsilon_2\circ\mu_r=\varepsilon_2\otimes\varepsilon_1\,,\\
(\id_{B_2}\otimes\varepsilon_1)\circ\nu_l=\eta_2\circ\varepsilon_1\,,
\quad\nu_l\circ\eta_1=\eta_2\otimes\eta_1\,,\\[5pt]
\mu_l\circ(\eta_2\otimes\id_{B_1})=\id_{B_1}\,,\quad
\mu_l\circ(\id_{B_2}\otimes\eta_1)=\eta_1\circ\varepsilon_2\,,\quad
\varepsilon_1\circ\mu_l=\varepsilon_2\otimes\varepsilon_1\,,\\
(\id_{B_2}\otimes\varepsilon_1)\circ\nu_r=\id_{B_2}\,,\quad
(\varepsilon_2\otimes\id_{B_1})\circ\nu_r=\eta_1\circ\varepsilon_2\,,
\quad\nu_r\circ\eta_2=\eta_2\otimes\eta_1\,,\\[5pt]
\sigma\circ(\eta_2\otimes\id_{B_2})=\sigma\circ(\id_{B_2}\otimes\eta_2)=
\eta_1\circ\varepsilon_2\,,\quad
\varepsilon_1\circ\sigma=\varepsilon_2\otimes\varepsilon_2\,,\\
(\id_{B_1}\otimes\varepsilon_1)\circ\rho=
(\varepsilon_1\otimes\id_{B_1})\circ\rho=\eta_1\circ\varepsilon_2\,,\quad
\rho\circ\eta_2=\eta_1\otimes\eta_1
\end{gathered}
\end{equation}
are satisfied.
\item The subsequent compatibility relations hold.
\end{enumerate}
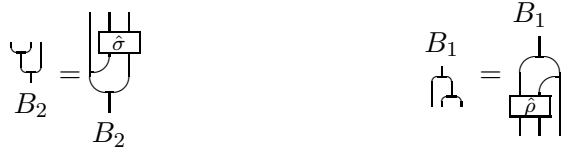
\def\FillCircDiam{5}
\begin{equation*}
{\begin{array}{c}
\divide\unitlens by 2
\begin{tangle}
\step[2]\object{\ }\\
\hh\hcu\hstep\id\\
\hh\hstep\hcu\\
\step\object{B_2}
\end{tangle}
=
\begin{tangle}
        \step[2]\object{\ }\\
	\hh\id\step\id\step\id\\
	\hh\id\hstep\ffbox 2{\hat\sigma}\\
	\hh\ru\step\id\\
	\cu\\
	\step\object{B_2}
\end{tangle}
\multiply\unitlens by 2
\\
\\
\text{Weak associativity of $\m_2$,}
\end{array}}
\qquad
{\begin{array}{c}
\divide\unitlens by 2
\begin{tangle}
\hstep\object{B_1}\\
\hh\hcd\\
\hh\id\hstep\hcd\\
\step[2]\object{\ }
\end{tangle}
=
\begin{tangle}
\step\object{B_1}\\
\hstep\cd\\
\hh\hstep\id\step\ld\\
\hh\ffbox 2{\hat\rho}\hstep\id\\
\hh\hstep\id\step\id\step\id\\
\step[2]\object{\ }\\
\end{tangle}
\multiply\unitlens by 2
\\
\\
\text{Weak coassociativity of $\Delta_1$,}
\end{array}}
\end{equation*}
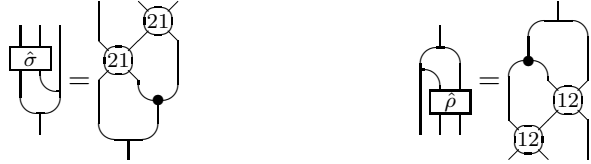
\begin{equation*}
{\begin{array}{c}
\divide\unitlens by 2
\vstretch 110
\begin{tangle}
\hh\hstep\id\step\id\step\id\\
\hh\ffbox2{\hat\sigma}\hstep\id\\
\hh\hstep\id\step\lu\\
\hstep\cu
\end{tangle}
\vstretch 100
\multiply\unitlens by 2
=
\divide\unitlens by 2
\begin{tangle}
\id\Step\ox{21}\\
\ox{21}\Step\id\\
\id\Step\cu*\\
\hstr{150}\cu
\end{tangle}
\multiply\unitlens by 2
\\
\\
\text{Weak associativity of $\mu_l$,}
\end{array}}
\qquad
{\begin{array}{c}
\divide\unitlens by 2
\vstretch 110
\begin{tangle}
	\cd\\
	\hh\rd\step\id\\
	\hh\id\hstep\ffbox2{\hat\rho}\\
	\hh\id\step\id\step\id
\end{tangle}
\vstretch 100
\multiply\unitlens by 2
=
\divide\unitlens by 2
\begin{tangle}
	\step\hstr{150}\cd\\
	\cd*\Step\id\\
	\id\Step\ox{12}\\
	\ox{12}\Step\id
\end{tangle}
\multiply\unitlens by 2
\\
\\
\text{Weak coassociativity of $\nu_r$,}
\end{array}}
\end{equation*}
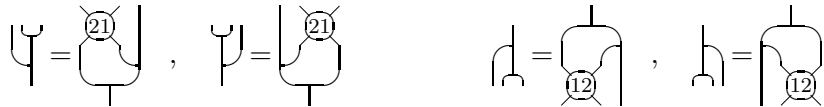
\begin{equation*}
{\begin{array}{c}
\divide\unitlens by 2
\begin{tangle}
\hh\id\hstep\hcu\\
\lu
\end{tangle}
=
\begin{tangle}
	\ox{21}\step\id\\
	\hh\id\Step\lu\\
	\hstr{150}\cu
\end{tangle}
\multiply\unitlens by 2
\quad ,\quad
\divide\unitlens by 2
\begin{tangle}
\hh\hcu\hstep\id\\
\hstep\ru
\end{tangle}
=
\begin{tangle}
	\id\step\ox{21}\\
	\hh\ru\Step\id\\
	\hstr{150}\cu
\end{tangle}
\\
\\
\text{Module-algebra compatibility,}
\end{array}}
\qquad
{\begin{array}{c}
\divide\unitlens by 2
\begin{tangle}
\ld\\
\hh\id\hstep\hcd
\end{tangle}
=
\begin{tangle}
        \hstr{150}\cd\\
        \hh\id\Step\ld\\
        \ox{12}\step\id
\end{tangle}
\multiply\unitlens by 2
\quad ,\quad
\divide\unitlens by 2
\begin{tangle}
\hstep\rd\\
\hh\hcd\hstep\id
\end{tangle}
=
\begin{tangle}
        \hstr{150}\cd\\
        \hh\rd\Step\id\\
        \id\step\ox{12}
\end{tangle}
\multiply\unitlens by 2
\\
\\
\text{Comodule-coalgebra compatibility,}
\end{array}}
\end{equation*}
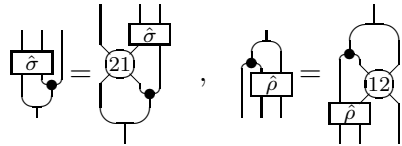
\begin{equation*}
\begin{array}{c}
\divide\unitlens by 2
\vstretch 110
\begin{tangle}
	\hh\hstep\id\step\id\step\id\\
	\hh\ffbox2{\hat\sigma}\hstep\id\\
	\hh\hstep\id\step\cu*\\
	\hh\hstep\hstr{150}\cu
\end{tangle}
\vstretch 100
=
\begin{tangle}
	\hh\id\Step\id\step\id\\
	\hh\id\step[1.5]\ffbox2{\hat\sigma}\\
	\ox{21}\step\id\\
	\hh\id\Step\cu*\\
	\hstr{125}\cu
\end{tangle}
\quad ,\quad
\vstretch 110
\begin{tangle}
	\hh\hstep\hstr{150}\cd\\
	\hh\cd*\step\id\\
	\hh\id\hstep\ffbox2{\hat\rho}\\
	\hh\id\step\id\step\id
\end{tangle}
\vstretch 100
=
\begin{tangle}
	\step\hstr{125}\cd\\
	\hh\hstep\cd*\Step\id\\
	\hstep\id\step\ox{12}\\
	\hh\ffbox2{\hat\rho}\step[1.5]\id\\
	\hh\hstep\id\step\id\Step\id
\end{tangle}
\multiply\unitlens by 2
\\
\\
\text{Cocycle and cycle compatibilities,}
\end{array}
\end{equation*}
\begin{equation*}
\begin{array}{c}
\divide\unitlens by 2
\begin{tangle}
\object{B_1}\step[2]\object{B_1}\\
\cu\\
\cd\\
\object{\ }\step[2]\object{\ }
\end{tangle}
=
 \begin{tangle}
      \hh \cd\step\cd \\
      \hh \id\hstep\hld\step\id\step\id \\
      \id\hstep\id\hstep\hx\step\id \\
      \hh \id\hstep\hlu\step\id\step\id \\
      \hh \cu\step\cu
 \end{tangle}
\multiply\unitlens by 2
\quad ,\quad
\divide\unitlens by 2
\begin{tangle}
\object{B_2}\step[2]\object{B_2}\\
\cu\\
\cd\\
\object{\ }\step[2]\object{\ }
\end{tangle}
=
\begin{tangle}
      \hh \cd\step\cd \\
      \hh \id\step\id\step\hrd\hstep\id \\
      \id\step\hx\hstep\id\hstep\id \\
      \hh \id\step\id\step\hru\hstep\id \\
      \hh \cu\step\cu
 \end{tangle}
\multiply\unitlens by 2
\\
\\
\text{Algebra-coalgebra compatibility,}
\end{array}
\end{equation*}
\def\FillCircDiam{6}
\begin{equation*}
{\begin{array}{c}
\divide\unitlens by 2
\begin{tangle}
\ox{21}\\
\hh\dh\step\ddh\\
\hh\ffbox2{\hat\rho}\\
\hh\hstep\id\step\id
\end{tangle}
=
\multiply\unitlens by 2
\divide\unitlens by 3
\begin{tangle}
\step\Cd\step[2]\Cd\\
\cd*\step[2]\cd\step\id\step[3]\ld\\
\vstr{110}\id\step[2]\ox{\sstyle 12}\step[2]\hx\step[2]\dd\step\id\\
\id\step[2]\id\step[2]\x\step[1]\x\step[2]\id\\
\vstr{110}\id\step[2]\ox{\sstyle 21}\step[2]\hx\step[2]\d\step\id\\
\cu\step[2]\cu*\step\id\step[3]\lu\\
\step\Cu\step[2]\Cu
\end{tangle}
\multiply\unitlens by 3
\quad ,\quad
\divide\unitlens by 2
\begin{tangle}
\hh\hstep\ru\\
\hcd
\end{tangle}
\multiply\unitlens by 2
=
\divide\unitlens by 3
\begin{tangle}
      \hcd\step\cd \\
      \id\step\hx\step\ld \\
      \ru\step\hx\step\id \\
      \cu\step\ru
\end{tangle}
\multiply\unitlens by 3
\\
\\
\text{Module-coalgebra compatibility,}
\end{array}}
\quad\qquad
{\begin{array}{c}
\divide\unitlens by 2
\begin{tangle}
\hstep\hcu\\
\hh\ld
\end{tangle}
\multiply\unitlens by 2
=
\divide\unitlens by 3
 \begin{tangle}
      \ld\step\cd \\
      \id\step\hx\step\ld \\
      \ru\step\hx\step\id \\
      \cu\step\hcu
\end{tangle}
\multiply\unitlens by 3
\quad ,\quad
\divide\unitlens by 2
\begin{tangle}
\hh\hstep\id\step\id\\
\hh\ffbox2{\hat\sigma}\\
\hh\hdd\step\hd\\
\ox{12}
\end{tangle}
\multiply\unitlens by 2
=
\divide\unitlens by 3
\begin{tangle}
\Cd\step[2]\Cd\\
\rd\step[3]\id\step\cd*\step[2]\cd\\
\vstr{110}\id\step\d\step[2]\hx\step[2]\ox{\sstyle 12}\step[2]\id\\
\id\step[2]\x\step\x\step[2]\id\step[2]\id\\
\vstr{110}\id\step\dd\step[2]\hx\step[2]\ox{\sstyle 21}\step[2]\id\\
\ru\step[3]\id\step\cu\step[2]\cu*\\
\Cu\step[2]\Cu
\end{tangle}
\multiply\unitlens by 3
\\
\\
\text{Comodule-algebra compatibility,}
\end{array}}
\end{equation*}
\begin{equation*}
\begin{array}{c}
\divide\unitlens by 2
\begin{tangle}
\ox{21}\\
\hh\id\step[2]\id\\
\ox{12}
\end{tangle}
\multiply\unitlens by 2
=
\divide\unitlens by 3
 \begin{tangle}
      \cd\step\cd \\
      \rd\step\hx\step\ld \\
      \id\step\hx\step\hx\step\id \\
      \ru\step\hx\step\lu \\
      \cu\step\cu
 \end{tangle}
\multiply\unitlens by 3
\\
\\
\text{Module-comodule compatibility,}
\end{array}
\end{equation*}
\begin{equation*}
\begin{array}{c}
\divide\unitlens by 2
\begin{tangle}
\hh\hstep\id\step\id\\
\hh\ffbox2{\hat\sigma}\\
\hh\hstep\id\step\id\\
\hh\ffbox2{\hat\rho}\\
\hh\hstep\id\step\id
\end{tangle}
\multiply\unitlens by 2
=
\divide\unitlens by 3
\begin{tangle}
\step\Cd\step[3]\Cd\\
\cd*\step[2]\cd\step\cd*\step[2]\cd\\
\vstr{110}\id\step[2]\ox{\sstyle 12}\step[2]\hx\step[2]\ox{\sstyle 12}%
 \step[2]\id\\
\id\step[2]\id\step[2]\x\step\x\step[2]\id\step[2]\id\\
\vstr{110}\id\step[2]\ox{\sstyle 21}\step[2]\hx\step[2]\ox{\sstyle 21}%
 \step[2]\id\\
\cu\step[2]\cu*\step\cu\step[2]\cu*\\
\step\Cu\step[3]\Cu
\end{tangle}
\multiply\unitlens by 3
\\
\\
\text{Cycle-cocycle compatibility.}
\\
\text{(End of Definition \ref{hp})}
\end{array}
\end{equation*}
\end{definition}
\def\FillCircDiam{4}
\abs
Observe that Definition \ref{hp} is $\pi$-symmetric invariant.
It will be verified in the next proposition that every cocycle cross
product bialgebra canonically induces a Hopf datum.

\begin{proposition}\Label{cp-hp}
Let $B_1\cpcybi[\mathfrak{c}]B_2$ be a cocycle cross product
bialgebra with corresponding structure morphisms
$\m_1$, $\eta_1$, $\Delta_1$ ,$\varepsilon_1$, $\m_2$, $\eta_2$, $\Delta_2$,
$\varepsilon_2$, $\mu_l$, $\mu_r$, $\nu_l$, $\nu_r$, $\rho$, and $\sigma$.
Then $\big((B_1,\m_1,\eta_1,\Delta_1,\varepsilon_1),
(B_2,\m_2,\eta_2,\Delta_2,\varepsilon_2);
\mu_l,\mu_r,\nu_l,\nu_r,\rho,\sigma\big)$ is a Hopf datum.
\end{proposition}

\begin{proof}
By definition cocycle cross product bialgebras are
cocycle cross product algebras and cycle cross product coalgebras
in particular. Hence $(B_1,\m_1\eta_1)$ is an algebra and
$(B_2,\Delta_2,\varepsilon_2)$
is a coalgebra. Since cocycle cross product bialgebras are
normalized one deduces with \eqref{type-alpha1} that
$\eta_2:\E\to B_2$ is an algebra morphism and
$\varepsilon_1:B_1\to \E$ is a coalgebra morphism. Then the conditions
of Definition \ref{hp}.1 and \ref{hp}.2 hold. The first sixteen
(co\n-)\-unital identities of \eqref{hp1} hold by assumption
(see \eqref{type-alpha1}). Then from \eqref{phi-sigma-rho} one easily derives
\begin{equation}\mathlabel{cp-hp1}
{\begin{split}
\varphi_{1,2}\circ(\eta_1\otimes\id_{B_2})&=\nu_r\,,\\
\varphi_{1,2}\circ(\id_{B_1}\otimes\eta_2)&=\nu_l\,,
\end{split}}
\quad
{\begin{split}
(\varepsilon_1\otimes\id_{B_2})\circ\varphi_{2,1}&=\mu_r\,,\\
(\id_{B_1}\otimes\varepsilon_2)\circ\varphi_{2,1}&=\mu_l\,,
\end{split}}
\quad
{\begin{split}
\hat\rho\circ(\eta_1\otimes\id_{B_1})&=\rho\,\\
(\id_{B_2}\otimes\varepsilon_2)\circ\hat\sigma&=\sigma\,.
\end{split}}
\end{equation}
Since $B_1\cpcybi[\mathfrak{c}] B_2$ is especially cocycle cross product
algebra, the identities \eqref{cond3-caat} are satisfied. Composing
the fifth equation of \eqref{cond3-caat} with
$(\varepsilon_1\otimes\id_{B_2})\circ\vert$ and using
\eqref{cp-hp1} proves that $(B_2,\mu_r)$ is right $B_1$-module.
Similarly the weak associativity of $\mu_l$ is shown by composing
with $(\id_{B_1}\otimes\varepsilon_2)\circ\vert$. From \eqref{cp-hp1}
and \eqref{cond3-caat} it will be concluded straightforwardly that
$\sigma\circ(\eta_2\otimes\id_{B_2}) =(\id_{B_1}\otimes\varepsilon_2)
 \circ\hat\sigma\circ(\eta_2\otimes\id_{B_2})
=(\id_{B_1}\otimes\varepsilon_2)\circ(\eta_1\otimes\id_{B_2})
=\eta_1\circ\varepsilon_2$
and similarly $\sigma\circ(\eta_2\otimes\id_{B_2})=\eta_1\circ\varepsilon_2$.
Since $B_1\cpcybi[\mathfrak{c}] B_2$ is cross product bialgebra,
$\varepsilon_1\otimes\varepsilon_2$ is algebra morphism, and
therefore one easily shows that $\varepsilon_1\circ\sigma=
\varepsilon_2\otimes\varepsilon_2$. Eventually all other (co\n-)unital
identities of \eqref{hp1} follow then by $\pi$-symmetry.
The weak associativity for $\m_2$ follows now from the sixth equation
of \eqref{cond3-caat} by adjoining $(\varepsilon_1\otimes\id_{B_2})\circ\vert$
on both sides. The module-algebra compatibility will be derived
from the fifth equation of \eqref{cond3-caat} through
composition with $(\id_{B_1}\otimes\varepsilon_2)\circ\vert$.
Making use of the bialgebra identity $\Delta_B\circ\m_B=
(\m_B\otimes\m_B)\circ(\id_B\otimes\Psi_{B,B}\otimes\id_B)\circ
(\Delta_B\otimes\Delta_B)$ for $B=B_1\cpcybi[\mathfrak{c}] B_2$ we prove
the algebra-coalgebra compatibility by application of
$\vert\circ(\id_{B_1}\otimes\eta_2\otimes\id_{B_1}\otimes\eta_2)$
and $(\id_{B_1}\otimes\id\varepsilon_2\otimes\id_{B_1}\otimes
\varepsilon_2)\circ\vert$ on both sides of the bialgebra identity.
The first equation of the module-coalgebra compatibility will be shown
similarly by composing with $\vert\circ(\eta_1\otimes\id_{B_2}\otimes
\id_{B_1}\otimes\eta_2)$ and ${(\varepsilon_1\otimes\id_{B_2}\otimes
\varepsilon_1\otimes\id_{B_2})\circ\vert}$. Application of
$\vert\circ(\eta_1\otimes\id_{B_2}\otimes\id_{B_1}\otimes\eta_2)$
and $(\id_{B_1}\otimes\varepsilon_2\otimes\id_{B_1}\otimes
\varepsilon_2)\circ\vert$ yields the second equation of
the module-coalgebra compatibility. The module-comodule compatibility
is derived from the bialgebra identity by composing with
$\vert\circ(\eta_1\otimes\id_{B_2}\otimes\id_{B_1}\otimes\eta_2)$
and $(\varepsilon_1\otimes\id_{B_2}\otimes\id_{B_1}\otimes\varepsilon_2)
\circ\vert$, whereas the cycle-cocycle compatibility
comes from composition with
$\vert\circ(\eta_1\otimes\id_{B_2}\otimes\eta_1\otimes\id_{B_2})$
and $(\id_{B_1}\otimes\varepsilon_2\otimes\id_{B_1}\otimes
\varepsilon_2)\circ\vert$. Application of
$(\id_{B_1}\otimes\varepsilon_2)\circ\vert$ to the sixth equation of
\eqref{cond3-caat} and use of the sixth identity of \eqref{cp-hp1}
yield the cocycle compatibility.
All remaining
compatibility relations of Definition \ref{hp} will be proven by
application of $\pi$-symmetry to the former results.\end{proof}

\subsection{Strong Hopf Data}

The converse of Proposition \ref{cp-hp} is not true in general.
However, in the following we show that so-called strong Hopf
data yield cross product bialgebras. The definition of strong Hopf data
is closely related to the definition of strong cross
product bialgebras.

\begin{definition}\Label{strong-hp}
A Hopf datum $\mathfrak{h}$ is called strong
if in addition the following identities hold.
\begin{gather}\mathlabel{act-coact-triv}
\divide\unitlens by 3
\multiply\unitlens by 2
\begin{tangle}
	\hh\hld\step\id\\
	\hh\id\hstep\x\\
	\hh\hru\step\id\\
\end{tangle}
=
\begin{tangle}
\ld\step\counit
\end{tangle}
\quad ,
\quad
\begin{tangle}
	\hh\id\step\hld\\
	\hh\x\hstep\id\\
	\hh\id\step\hru\\
\end{tangle}
=
\begin{tangle}\unit\step\ru
\end{tangle}
\\
\mathlabel{cocycle-triv1}
\divide\unitlens by 2
\vstretch 60
\begin{tangle}
\Ld\step\id\\
\rd\step\hx\\
\id\step\hx\step\id\\
\hcu*\step\id\step\id
\end{tangle}
=
\begin{tangle}
\unit\step\unit\step\vstr{40}\id\step\id\\
\vstr{40}\id\step\id\step\id\step\counit
\end{tangle}
\quad ,
\quad
\begin{tangle}
\id\step\id\step\hcd*\\
\id\step\hx\step\id\\
\hx\step\lu\\
\id\step\Ru
\end{tangle}
=
\begin{tangle}
\unit\vstr{40}\step\id\step\id\step\id\\
\vstr{40}\id\step\id\step\counit\step\counit
\end{tangle}
\vstretch 100
\divide\unitlens by 3
\multiply\unitlens by 4
\quad ,
\quad
\vstretch 120
\begin{tangle}
\hh\id\step\hld\\
\hh\cu*\hstep\id
\end{tangle}
=
\vstretch 90
\begin{tangle}
\hh\step\id\\
\hh\counit\step\id\\
\hh\unit\step\id
\end{tangle}
\quad ,
\quad
\vstretch 120
\begin{tangle}
\hh\id\hstep\cd*\\
\hh\hru\step\id
\end{tangle}
=
\vstretch 90
\begin{tangle}
\hh\id\\
\hh\id\step\counit\\
\hh\id\step\unit
\end{tangle}
\\[3pt]
\mathlabel{strong-hopf-combo}
\divide\unitlens by 3
\def\FillCircDiam{6}
\vstretch 50
\begin{tangle}
\vstr{70}\id\Step\id\Step\cd*\\
\id\Step\id\step\cd\step\id\\
\id\Step\hx\Step\id\step\id\\
\id\step\dd\ld\step\dd\step\id\\
\hx\step\id\step\hx\Step\id\\
\id\step\d\lu\step\id\Step\id\\
\id\Step\ru\step\id\Step\id
\end{tangle}
=
\vstretch 100
\begin{tangle}
\vstr{70}\id\step\id\step[1.5]\cd*\\
\hh\id\step\id\step\cd\step[1.5]\id\\
\hh\id\step\x\step\id\step[1.5]\id\\
\hh\x\step\x\step[1.5]\id\\
\hh\id\step{\hstr{200}\hru}\step\id\step[1.5]\id
\end{tangle}
\quad ,
\quad
\vstretch 50
\begin{tangle}
\id\Step\id\step\ld\Step\id\\
\id\Step\id\step\rd\d\step\id\\
\id\Step\hx\step\id\step\hx\\
\id\step\dd\step\ru\dd\step\id\\
\id\step\id\Step\hx\Step\id\\
\id\step\cu\step\id\Step\id\\
\vstr{70}\cu*\Step\id\Step\id
\end{tangle}
=
\vstretch 100
\begin{tangle}
\hh\id\step[1.5]\id\step{\hstr{200}\hld}\step\id\\
\hh\id\step[1.5]\x\step\x\\
\hh\id\step[1.5]\id\step\x\step\id\\
\hh\id\step[1.5]\cu\step\id\step\id\\
\vstr{70}\cu*\step[1.5]\id\step\id
\end{tangle}
\end{gather}
\end{definition}
\abs
Strong Hopf data are invariant under $\pi$-symmetry.
This fact will be used extensively in Section \ref{proof-hp-cp} where we prove
Theorem \ref{hp-cp}.

\begin{remark}\Label{non-except-hopf}
{\normalfont Observe that except of \eqref{strong-hopf-combo} 
and the first and the second identity of \eqref{cocycle-triv1} the
defining identities of a strong Hopf datum are either of rank
$(2,2)$, $(1,3)$ or $(3,1)$. Therefore
in the same way as in Remark \ref{non-except-bialg}
we call a strong Hopf datum \textit{regular} if the
defining identities are of either rank $(2,2)$, $(1,3)$ or $(3,1)$.
We call the strong Hopf datum \textit{pure} if \eqref{act-coact-triv},
\eqref{cocycle-triv1} and \eqref{strong-hopf-combo} are redundant.}
\end{remark}

\subsection{Main Results}

This is the central part of the article. In the subsequent
Theorems \ref{hp-cp} and \ref{central-theory} we present
the universal (co\n-)modular co-cyclic theory of (strong) cross
product bialgebras. In Theorem \ref{hp-cp} we describe the (co\n-)modular
co-cyclic construction of strong cross product bialgebras.
Theorem \ref{central-theory} exhibits its universality.

\begin{theorem}\Label{hp-cp}
Let $\mathfrak{h}=
\big((B_1,\m_1,\eta_1,\Delta_1,\varepsilon_1),(B_2,\m_2,\eta_2,\Delta_2,
\varepsilon_2);\mu_l,\mu_r,\nu_l,\nu_r,\rho,\sigma\big)$ be a strong
Hopf datum. Then $B_1\otimes B_2$ with $\varphi_{1,2}$,
$\varphi_{2,1}$, $\hat\sigma$, and $\hat\rho$ defined according to
\eqref{phi-sigma-rho2} is a strong cross product bialgebra
$B_1\overset{\mathfrak{h}}\bowtie B_2$. Strong Hopf data and strong
cross product bialgebras are in one-to-one correspondence.
\endproof
\end{theorem}
\abs
Every strong Hopf datum $\mathfrak{h}$ therefore induces a strong
cross product bialgebra $B_1\overset{\mathfrak{h}}\bowtie B_2$. For any
bialgebra $B$ which is isomorphic to $B_1\overset{\mathfrak{h}}\bowtie B_2$ the
additional relations \eqref{act-coact-triv} -- \eqref{strong-hopf-combo} (or by
one-to-one correspondence the equations \eqref{strong-bi1} --
\eqref{strong-bi3}) imply the additional relations \eqref{pi-strong} for the
idempotents $\Pi_1$ and $\Pi_2$ of $B$.

\begin{theorem}\Label{central-theory}
Let $B$ be a bialgebra in $\C$. Then the following statements are
equivalent.
\begin{entry}
\item[1.]
There is a strong Hopf datum $\mathfrak{h}$ such that the
corresponding strong cross product
bialgebra $B_1\overset{\mathfrak{h}}\bowtie B_2$
is bialgebra isomorphic to $B$.
\item[2.]
There are idempotents $\Pi_1,\Pi_2\in \End(B)$ such that
the conditions of Theorem \ref{proj-cycle}.2 and the
``strong projection'' relations \eqref{pi-strong} hold.
\item[3.] 
There are objects $B_1$ and $B_2$ and
morphisms $B_1\overset{\inj_1}\to B\overset{\proj_1}\to B_1$
and $B_2\overset{\inj_2}\to B\overset{\proj_2}\to B_2$ such that
the conditions of Theorem \ref{proj-cycle}.3 and
the ``strong projection'' relations \eqref{pi-strong} hold
with $\Pi_j:= \inj_j\circ\proj_j$.
\end{entry}
\end{theorem}

\begin{proof}
Theorem \ref{central-theory} can be derived
straightforwardly from Theorem \ref{hp-cp}, Theorem \ref{proj-cycle}, 
(Theorem \ref{bialg-cocycle2}), Remark \ref{alpha-ip}, 
Remark \ref{non-except-bialg} and the relations \eqref{pi-strong}.\end{proof}
\abs
\abs
In the sequel we will consider special versions of strong cross
product bialgebras. They admit equivalent universal and (co\n-)modular
co-cyclic descriptions and can be derived from the most general 
construction given in Theorem \ref{hp-cp}.
In particular all known constructions of cocycle cross product bialgebras
will be recovered, and additionally we find several new types of cocycle
cross product bialgebras. In Proposition \ref{bialg-cocycle6}
we describe comprehensively the universal and (co\n-)modular
co-cyclic properties of one of these special types.

Theorem \ref{hp-cp} and Proposition \ref{alg-coalg-triv}
provide the necessary tools to describe the various subclasses.
We distinguish the different special versions of strong
cross product bialgebras with the help of the boxes
$\divide\unitlens by 3
\multiply\unitlens by 2
\begin{tangle}
 \hh\ffbox1{\normalfont\mu_l}\ffbox1{\normalfont\mu_r}\\
 \hh\ffbox1{\normalfont\sigma}\ffbox1{\normalfont\rho}\\
 \hh\ffbox1{\normalfont\nu_l}\ffbox1{\normalfont\nu_r}
 \end{tangle}
\divide\unitlens by 2
\multiply\unitlens by 3$
where the particular entries will be left blank ``$\ $" or take the values
$\bullet$ for the (weak) (co\n-)actions, $\blacktriangledown$ for $\sigma$ and
$\blacktriangle$ for $\rho$ dependent on the respective morphisms
are trivial or not. For example a strong cross product bialgebra
with trivial cocycle $\sigma=\eta_1\circ(\varepsilon_2\otimes\varepsilon_2)$
and trivial right coaction $\nu_r=\id_{B_2}\otimes\eta_1$ will be
represented by the classification box
$\divide\unitlens by 3
\multiply\unitlens by 2
\begin{tangle}
 \hh\ffbox1{\normalfont\bullet}\ffbox1{\normalfont\bullet}\\
 \hh\ffbox1{\normalfont\ }\ffbox1{\normalfont\blacktriangle}\\
 \hh\ffbox1{\normalfont\bullet}\ffbox1{\normalfont\ }
 \end{tangle}\,$.
\begin{figure}
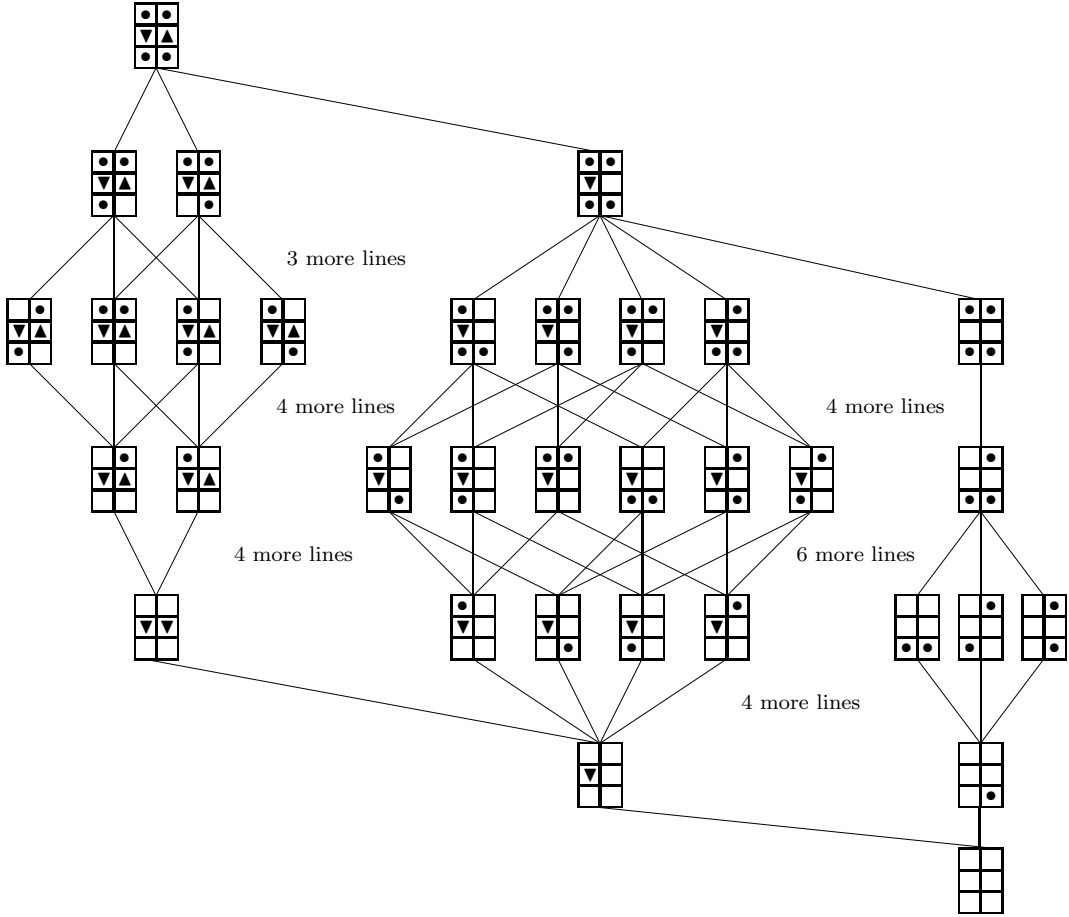

\begin{equation*}
\unitlens 8pt
\providecommand{\0}{\ffbox1{}}
\providecommand{\9}{\ffbox1{\bullet}}
\providecommand{\7}{\ffbox1{\blacktriangledown}}
\providecommand{\8}{\ffbox1{\blacktriangle}}
\providecommand{\4}[1]{\Put(0,5)[lt]{\begin{tangle}\vstr{250}\hstr{1#100}%
 \nw1\end{tangle}}}
\providecommand{\5}{\Put(0,0)[cc]{\vdots}}
\providecommand{\morelines}[3]{\Put(#2,#3)[cc]{\text{#1 more lines}}} 
\begin{tangles}{ccccc}
\HH\9\9 &\step[3]& &\step[3]&\\
\HH\7\8 && &&\\
\HH\9\9 && &&\\
\vstr{200}\step[5]\hstr{200}\dd\d\step[-1]{\hstr{700}\nw3} && &&\\
\HH\9\9\Step\9\9 &&\9\9 &&\\
\HH\7\8\Step\7\8 &&\7\0 &&\\
\HH\9\0\Step\0\9 &&\9\9 &&\\
\vstr{200}\hstr{400}\dd\id\X\id\d &\morelines{3}{5}{20}& %
  \vstr{200}\hstr{200}\step[2]\sw3\step\dd\d\step[-1]\nw3%
  \step[-1]{\hstr{600}\nw3}\Step&&\\
\HH\0\9\Step\9\9\Step\9\0\Step\9\0 && %
  \9\0\Step\9\9\Step\9\9\Step\0\9 &&\9\9\\
\HH\7\8\Step\7\8\Step\7\8\Step\7\8 &&%
   \7\0\Step\7\0\Step\7\0\Step\7\0 && \0\0\\
\HH\9\0\Step\0\0\Step\9\0\Step\0\9 &&%
   \9\9\Step\0\9\Step\9\0\Step\9\9 && \9\9\\
\vstr{200}\hstr{400}\d\id\X\id\dd%
   &\morelines{4}{0}{20}& \vstr{200}\hstr{400}\sw1\step[-1]\sw2\id\sw2%
   \vstr{200}\step[-1]\nw2\id\sw1\step[-1]\nw2\sw1\step[-1]\nw2\id\nw1
   &\morelines{4}{10}{20}&\vstr{200}\id\\
\HH\0\9\Step\9\0 &&\9\0\Step\9\0\Step\9\9\Step\0\0\Step\0\9\Step\0\9 &&\0\9\\
\HH\7\8\Step\7\8 &&\7\0\Step\7\0\Step\7\0\Step\7\0\Step\7\0\Step\7\0 &&\0\0\\
\HH\0\0\Step\0\0 &&\0\9\Step\9\0\Step\0\0\Step\9\9\Step\0\9\Step\9\0 &&\9\9\\
\vstr{200}\hstr{200}\d\dd &\morelines{4}{-20}{20}&\vstr{200}\hstr{400}\nw1\step[-1]\nw2%
   \id\sw1\step[-1]\nw2\XX\step[-2]\sw1\id\sw2\id\sw1 &\morelines{6}{-4}{20}& %
   \vstr{200}\sw3\Step\id\nw3\Step\\
\HH\0\0 && \9\0\Step\0\0\Step\0\0\Step\0\9 && \0\0\step\0\9\step\0\9\\
\HH\7\7 && \hh\7\0\Step\7\0\Step\7\0\Step\7\0 && \hh\0\0\step\0\0\step\0\0\\
\HH\0\0 && \0\0\Step\0\9\Step\9\0\Step\0\0 && \9\9\step\9\0\step\0\9\\%
 \vstr{200}\step[6]\hstr{720}\nw3&&\vstr{200}\hstr{200}\nw3\step\d\dd%
  \step[-1]\sw3\Step &\morelines{4}{-30}{20}&\vstr{200}\nw3\Step\id\sw3\Step\\
\HH && \0\0 && \0\0\\
\HH && \hh\7\0 && \hh\0\0\\
\HH && \0\0 && \0\9\\
 && && \step[-5.5]{\hstr{640}\se3}\step[-1]\id\\
\HH && && \0\0\\
\HH && && \0\0\\
\HH && && \0\0
\end{tangles}
\end{equation*}
\caption{Graph representing the various special versions of strong
cross product bialgebras.}\Label{fig-class}
\end{figure}
\unitlens 12pt
Thereby $2^6=64$ types of special strong cross products can be
obtained; $2^4=16$ of them are co-cycle free. In a similar (dual) way 
the dual strong cross product bialgebras will be denoted
by a classification box
$\divide\unitlens by 3
\multiply\unitlens by 2
\begin{tangle}
 \hh\ffbox1{\normalfont\mu_l}\ffbox1{\normalfont\mu_r}\\
 \hh\ffbox1{\normalfont\sigma}\ffbox1{\normalfont\rho}\\
 \hh\ffbox1{\normalfont\nu_l}\ffbox1{\normalfont\nu_r}
 \end{tangle}
\divide\unitlens by 2
\multiply\unitlens by 3$
with entries $\bullet$ for the (co\-)actions, $\blacktriangle$ for the
cocycle $\sigma$, and $\blacktriangledown$ for the cycle $\rho$. 
Then analogously $64$ special types of dual strong cross 
product bialgebras can be obtained from which $2^4=16$ 
(co\n-cycle free) types coincide with the corresponding 16 types
of strong cross product bialgebras.
Therefore the total number of different
types of strong and dual strong cross product
bialgebras which is classified by our scheme is $2^6+(2^6-2^4)=7\cdot 2^4=112$.
However, henceforth we do not tell between a certain type and its dual or 
$\pi$-symmetric counterparts - it is an easy exercise to obtain one from the 
other\footnote{Observe that a $\pi$-rotation of the classification box in the 
plane of the graphic yields the corresponding $\pi$-symmetric counterpart
of a certain type of cross product bialgebra. And analogously a reflection of
the classification box along a horizontal line yields the dual counterpart.}.
Up to this $\mathbb{Z}_2\times\mathbb{Z}_2$-symmetry we obtain $33$ different
classification boxes. 
In Figure \ref{fig-class} we present a graph where the boxes
(modulo $\mathbb{Z}_2\times\mathbb{Z}_2$-symmetry) are the vertices, and
each box is linked with its descending ``next neighbour" special versions. 
We stratify the whole graph into three layers. The layers consist of boxes 
with 2, 1 or 0 co-cycles respectively.
The most general type of strong cross product bialgebra
(Theorem \ref{hp-cp}) is at the top of the graph in the first row.
In the second row all descendants with one trivial morphism
$\mu_l$, $\mu_r$, $\sigma$, $\rho$, $\nu_l$, or $\nu_r$ are listed, 
in the third row all special types with two trivial morphisms are listed, etc. 
The cross product bialgebra constructions from 
\cite{Rad1:85,Ma1:90,MS:94,Ma6:94,BD1:98} are descendants of the special cases
$\divide\unitlens by 3
\multiply\unitlens by 2
\begin{tangle}
 \hh\ffbox1{\normalfont\bullet}\ffbox1{\normalfont\bullet}\\
 \hh\ffbox1{\normalfont\ }\ffbox1{\normalfont\ }\\
 \hh\ffbox1{\normalfont\bullet}\ffbox1{\normalfont\bullet}
 \end{tangle}
\divide\unitlens by 2
\multiply\unitlens by 3\,$
and
$\divide\unitlens by 3
\multiply\unitlens by 2
\begin{tangle}
 \hh\ffbox1{\normalfont\bullet}\ffbox1{\normalfont\ }\\
 \hh\ffbox1{\normalfont\blacktriangledown}%
 \ffbox1{\normalfont\blacktriangle}\\
 \hh\ffbox1{\normalfont\ }\ffbox1{\normalfont\bullet}
 \end{tangle}
\divide\unitlens by 2
\multiply\unitlens by 3\,$
which are the co-cycle free cross product bialgebras of \cite{BD1:98}
and the braided versions of bicrossed product bialgebras
\cite{Ma6:94,MS:94} respectively\footnote{A special symmetric version of
$\divide\unitlens by 3
\multiply\unitlens by 2
\begin{tangle}
 \hh\ffbox1{\normalfont\bullet}\ffbox1{\normalfont\bullet}\\
 \hh\ffbox1{\normalfont\blacktriangledown}\ffbox1{\normalfont\ }\\
 \hh\ffbox1{\normalfont\ }\ffbox1{\normalfont\bullet}
 \end{tangle}
\divide\unitlens by 2
\multiply\unitlens by 3$ has been studied recently in \cite{Sbg1:98}.}.
In the following tables we describe the different types of cocycle cross
product bialgebras in more detail; we omit the description of the 
various cross product bialgebras studied in 
\cite{Rad1:85,Ma1:90,MS:94,Ma6:94,BD1:98}.
In the second column of the tables the
structure of the multiplication $\m_B$ and the comultiplication $\Delta_B$
is given. In the third column we recall the corresponding equivalent
properties of Proposition \ref{alg-coalg-triv}, and in the last column
the status of the strong cross product bialgebra is listed.
\abs
\unitlens 15pt
\begin{center}
\begin{tabular}{|c|c|c|c|}
\hline Type&Structure morphisms $\m_B$, $\Delta_B$
&Proposition \ref{alg-coalg-triv}&Status\\
\hline &&&\\[-12pt]
\hline $\divide\unitlens by 2\begin{tangle}\classbox 111111\end{tangle}$
&$\mDelta{\mu_l}{\mu_r}{\nu_l}{\nu_r}\sigma\rho
{\hh\id\step\cd\step\cd\step\id\\
\id\step\id\step\hx\step\id\step\id\\
\hh\id\step\lu\step\ru\step\id\\
\hh\id\Step\id\hstep\cd\step\cd\\
\hh\id\Step\id\hstep\id\step\id\step\hrd\hstep\id\\
\cu\hstep\id\step\hx\hstep\id\hstep\id\\
\hh\step\id\step\hstep\id\step\id\step\hru\hstep\id\\
\hh\step\id\step[1.5]\cu*\step\cu\\
\step\cu\Step\id}
{\hstep\id\Step\cd\\
\hh\cd\step\cd*\step[1.5]\id\\
\hh \id\hstep\hld\step\id\step\id\hstep\step\id\\
\id\hstep\id\hstep\hx\step\id\hstep\cd\\
\hh\id\hstep\hlu\step\id\step\id\hstep\id\Step\id\\
\hh\cu\step\cu\hstep\id\Step\id\\
\hh\hstep\id\step\ld\step\rd\step\id\\
\hstep\id\step\id\step\hx\step\id\step\id\\
\hh\hstep\id\step\cu\step\cu\step\id}$
&
&
\\
\hline$\divide\unitlens by 2\begin{tangle}\classbox 111110\end{tangle}$
&$\mDelta{\mu_l}{\mu_r}{\nu_l}{}\sigma\rho
{\hh \id\step\cd\step\cd\step\id \\
\id\step\id\step\hx\step\id\step\id \\
\hh \id\step\lu\step\ru\step\id \\
\hh \id\Step\id\hstep\cd\step\cd \\
\cu\hstep\id\step\hx\step\id \\
\hh \step\id\step[1.5]\cu*\step\cu \\
\step\cu\Step\id}
{\cd\step[1.5]\hstr{225}\hcd\\
\id\step\ld\step\hcd*\step\id\\
\id\step\id\step\hx\step\id\step\id\\
\id\step\lu\step\hcu\step\id\\
\hh{\hstr{300}\cu}\hstep\ld\step\hcd\\
\hh\step\id\step[1.5]\id\step\x\step\id\\
\hh\step\id\step[1.5]\hcu\step\id\step\id}$
&\ref{alg-coalg-triv}.4
&
\\
\hline$\divide\unitlens by 2\begin{tangle}\classbox 111101\end{tangle}$
&$\mDelta{\mu_l}{\mu_r}{}{\nu_r}\sigma\rho
{\hh \id\step\cd\step\cd\step\id \\
\id\step\id\step\hx\step\id\step\id \\
\hh \id\step\lu\step\ru\step\id \\
\hh \id\Step\id\hstep\cd\step\cd \\
\hh \id\Step\id\hstep\id\step\id\step\hrd\hstep\id \\
\cu\hstep\id\step\hx\hstep\id\hstep\id \\
\hh \step\id\step\hstep\id\step\id\step\hru\hstep\id \\
\hh \step\id\step[1.5]\cu*\step\cu \\
\step\cu\Step\id}
{\hstep\id\Step\cd \\
\hh \cd\step\cd*\step[1.5]\id \\
\id\step\hx\step\id\hstep\cd \\
\hh \cu\step\cu\hstep\id\Step\id \\
\hh \hstep\id\Step\id\step\rd\step\id \\
\hstep\id\Step\hx\step\id\step\id \\
\hh \hstep\id\Step\id\step\cu\step\id}$
&\ref{alg-coalg-triv}.3
& 
\\
\hline$\divide\unitlens by 2\begin{tangle}\classbox 111011\end{tangle}$
&$\mDelta{\mu_l}{\mu_r}{\nu_l}{\nu_r}{}\rho
{\hh\id\step\cd\step\cd\step\id\\
\id\step\id\step\hx\step\id\step\id\\
\hh\id\step\lu\step\ru\step\id\\
\hh\id\Step\id\hstep\cd\step\cd\\
\hh\id\Step\id\hstep\id\step\id\step\hrd\hstep\id\\
\cu\hstep\id\step\hx\hstep\id\hstep\id\\
\hh\step\id\step[1.5]\id\step\id\step\hru\hstep\id\\
\hh\step\id\step[1.5]\cu*\step\cu\\
\step\cu\Step\id}
{\cd\step\cd\\
\hh\id\step\ld\step\rd\step\id\\
\id\step\id\step\hx\step\id\step\id\\
\hh\id\step\cu\step\cu\step\id}$
&\ref{alg-coalg-triv}.6
&
\\
\hline$\divide\unitlens by 2\begin{tangle}\classbox 011110\end{tangle}$
&$\mDelta{}{\mu_r}{\nu_l}{}\sigma\rho
{\hh\id\step\id\step\hcd\step\id\\
\hh\id\step\x\step\id\step\id\\
\hh\cu\step\ru\step\id\\
\hh\hstep\id\step\cd\step\cd\\
\hh\hstep\id\step\id\step\x\step\id\\
\hh\hstep\d\hstep\cu*\step\cu\\
\hh\step\cu\step[2]\id}
{\hh\hstep\id\step[2]\cd\\
\hh\cd\step\cd*\hstep\d\\
\hh\id\step\x\step\id\step\id\\
\hh\cu\step\cu\step\id\\
\hh\hstep\id\step\ld\step\cd\\
\hh\hstep\id\step\id\step\x\step\id\\
\hh\hstep\id\step\hcu\step\id\step\id}$
&\ref{alg-coalg-triv}.1, \ref{alg-coalg-triv}.4
& regular
\\
\hline$\divide\unitlens by 2\begin{tangle}\classbox 111100\end{tangle}$
&$\mDelta{\mu_l}{\mu_r}{}{}\sigma\rho
{\hh \id\step\cd\step\cd\step\id \\
\id\step\id\step\hx\step\id\step\id \\
\hh \id\step\lu\step\ru\step\id \\
\hh \id\Step\id\hstep\cd\step\cd \\
\cu\hstep\id\step\hx\step\id \\
\hh \step\id\step[1.5]\cu*\step\cu \\
\step\cu\Step\id}
{\hh\hstep\id\step[2]\cd\\
\hh\cd\step\cd*\hstep\d\\
\hh\id\step\x\step\id\step\id\\
\hh\cu\step\cu\hstep\cd\\
\hh\hstep\id\step[2]\x\step\id}$
&\ref{alg-coalg-triv}.3, \ref{alg-coalg-triv}.4
& 
\\
\hline
\end{tabular}
\end{center}

\begin{center}
\begin{tabular}{|c|c|c|c|}
\hline Type&Structure morphisms $\m_B$, $\Delta_B$
&Proposition \ref{alg-coalg-triv}&Status\\
\hline &&&\\[-12pt]
\hline$\divide\unitlens by 2\begin{tangle}\classbox 101110\end{tangle}$
&$\mDelta{\mu_l}{}{\nu_l}{}\sigma\rho
{\hh \id\step\cd\step\id\Step\id \\
\id\step\id\step\hx\Step\id \\
\hh \id\step\lu\step\id\Step\id \\
\hh \id\Step\id\hstep\cd\step\cd \\
\cu\hstep\id\step\hx\step\id \\
\hh \step\id\step[1.5]\cu*\step\cu \\
\step\cu\Step\id}
{\cd\step[1.5]\hstr{225}\hcd\\
\id\step\ld\step\hcd*\step\id\\
\id\step\id\step\hx\step\id\step\id\\
\id\step\lu\step\hcu\step\id\\
\hh{\hstr{300}\cu}\hstep\ld\step\hcd\\
\hh\step\id\step[1.5]\id\step\x\step\id\\
\hh\step\id\step[1.5]\hcu\step\id\step\id}$
&\ref{alg-coalg-triv}.2, \ref{alg-coalg-triv}.4
& regular
\\
\hline$\divide\unitlens by 2\begin{tangle}\classbox 101101\end{tangle}$
&$\mDelta{\mu_l}{}{}{\nu_r}\sigma\rho
{\hh \id\step\cd\step\id\Step\id \\
\id\step\id\step\hx\Step\id \\
\hh \id\step\lu\step\id\Step\id \\
\hh \id\Step\id\hstep\cd\step\cd \\
\cu\hstep\id\step\hx\step\id \\
\hh \step\id\step[1.5]\cu*\step\cu \\
\step\cu\Step\id}
{\hstep\id\Step\cd \\
\hh \cd\step\cd*\step[1.5]\id \\
\id\step\hx\step\id\hstep\cd \\
\hh \cu\step\cu\hstep\id\Step\id \\
\hh \hstep\id\Step\id\step\rd\step\id \\
\hstep\id\Step\hx\step\id\step\id \\
\hh \hstep\id\Step\id\step\cu\step\id}$
&\ref{alg-coalg-triv}.2, \ref{alg-coalg-triv}.3
& pure
\\
\hline$\divide\unitlens by 2\begin{tangle}\classbox 101011\end{tangle}$
&$\mDelta{\mu_l}{\mu_r}{}{\nu_r}{}\rho
{\hh\id\step\cd\step\id\Step\id\\
\id\step\id\step\hx\Step\id\\
\hh\id\step\lu\step\id\Step\id\\
\hh\id\Step\id\hstep\cd\step\cd\\
\cu\hstep\id\step\hx\step\id\\
\hh\step\id\step[1.5]\cu*\step\cu \\
\step\cu\Step\id}
{\cd\step\cd\\
\hh\id\step\ld\step\rd\step\id\\
\id\step\id\step\hx\step\id\step\id\\
\hh\id\step\cu\step\cu\step\id}$
&\ref{alg-coalg-triv}.2, \ref{alg-coalg-triv}.6
& 
\\
\hline$\divide\unitlens by 2\begin{tangle}\classbox 111001\end{tangle}$
&$\mDelta{\mu_l}{\mu_r}{}{\nu_r}\sigma{}
{\hh \id\step\cd\step\cd\step\id \\
\id\step\id\step\hx\step\id\step\id \\
\hh \id\step\lu\step\ru\step\id \\
\hh \id\Step\id\hstep\cd\step\cd \\
\hh \id\Step\id\hstep\id\step\id\step\hrd\hstep\id \\
\cu\hstep\id\step\hx\hstep\id\hstep\id \\
\hh \step\id\step\hstep\id\step\id\step\hru\hstep\id \\
\hh \step\id\step[1.5]\cu*\step\cu \\
\step\cu\Step\id}
{\cd\step\cd\\
\hh\id\Step\id\step\rd\step\id \\
\id\Step\hx\step\id\step\id \\
\hh\id\Step\id\step\cu\step\id}$
&\ref{alg-coalg-triv}.3, \ref{alg-coalg-triv}.6
& pure
\\
\hline$\divide\unitlens by 2\begin{tangle}\classbox 111010\end{tangle}$
&$\mDelta{}{\mu_r}{\nu_l}{\nu_r}{}\rho
{\hh\id\step\cd\step\cd\step\id\\
\id\step\id\step\hx\step\id\step\id\\
\hh\id\step\lu\step\ru\step\id\\
\hh\id\Step\id\hstep\cd\step\cd\\
\cu\hstep\id\step\hx\step\id\\
\hh\step\id\step[1.5]\cu*\step\cu\\
\step\cu\step[2]\id}
{\cd\step\cd\\
\hh\id\step\ld\step\id\Step\id\\
\id\step\id\step\hx\Step\id\\
\hh\id\step\cu\step\id\Step\id}$
&\ref{alg-coalg-triv}.8
& regular
\\
\hline$\divide\unitlens by 2\begin{tangle}\classbox 011011\end{tangle}$
&$\mDelta{\mu_l}{\mu_r}{\nu_l}{}{}\rho
{\hh\id\step\id\step\hcd\step[1.5]\id\\
\hh\id\step\x\step\id\step[1.5]\id\\
\hh\hcu\step\ru\hstep{\hstr{300}\cd}\\
\hstep\id\step\hcd\step\rd\step\id\\
\hstep\id\step\id\step\hx\step\id\step\id\\
\hstep\id\step\hcu*\step\ru\step\id\\
\hstep{\hstr{225}\hcu}\step[1.5]\cu}
{\cd\step\cd\\
\hh\id\step\ld\step\rd\step\id\\
\id\step\id\step\hx\step\id\step\id\\
\hh\id\step\cu\step\cu\step\id}$
&\ref{alg-coalg-triv}.1, \ref{alg-coalg-triv}.5
&
\\
\hline$\divide\unitlens by 2\begin{tangle}\classbox 110011\end{tangle}$
&
$\mDelta{\mu_l}{\mu_r}{\nu_l}{\nu_r}{}{}{
\hh\id\step\cd\step\cd\step\id \\
\hh\id\step\id\step\x\step\id\step\id \\
\hh\id\step\lu\step\ru\step\id \\
\cu\step\cu}
{\cd\step\cd\\
\hh\id\step\ld\step\rd\step\id \\
\hh\id\step\id\step\x\step\id\step\id \\
\hh\id\step\cu\step\cu\step\id}$
&\ref{alg-coalg-triv}.5 , \ref{alg-coalg-triv}.6
& regular
\\
\hline$\divide\unitlens by 2\begin{tangle}\classbox 010011\end{tangle}$
&$\mDelta{}{\mu_r}{\nu_l}{\nu_r}{}{}
{\hh\id\step[2]\id\step\hcd\step\id\\
\hh\id\step[2]\x\step\id\step\id\\
\hh\id\step[2]\id\step\ru\step\id\\
\cu\step\cu}
{\cd\step\cd\\
\hh\id\step\ld\step\rd\step\id\\
\hh\id\step\id\step\x\step\id\step\id\\
\hh\id\step\hcu\step\hcu\step\id}$
&\ref{alg-coalg-triv}.6, \ref{alg-coalg-triv}.7
& pure
\\
\hline
\end{tabular}
\end{center}

\begin{center}
\begin{tabular}{|c|c|c|c|}
\hline Type&Structure morphisms $\m_B$, $\Delta_B$
&Proposition \ref{alg-coalg-triv}&Status\\
\hline &&&\\[-12pt]
\hline$\divide\unitlens by 2\begin{tangle}\classbox 011100\end{tangle}$
&$\mDelta{}{\mu_r}{}{}\sigma\rho
{\hh\id\step\id\step\hcd\step\id\\
\hh\id\step\x\step\id\step\id\\
\hh\cu\step\ru\step\id\\
\hh\hstep\id\step\cd\step\cd\\
\hh\hstep\id\step\id\step\x\step\id\\
\hh\hstep\d\hstep\cu*\step\cu\\
\hh\step\cu\step[2]\id}
{\hh\hstep\id\step[2]\cd\\
\hh\cd\step\cd*\hstep\d\\
\hh\id\step\x\step\id\step\id\\
\hh\cu\step\cu\hstep\cd\\
\hh\hstep\id\step[2]\x\step\id}$
& \ref{alg-coalg-triv}.1, \ref{alg-coalg-triv}.3, \ref{alg-coalg-triv}.4
& regular
\\
\hline$\divide\unitlens by 2\begin{tangle}\classbox 101010\end{tangle}$
&$\mDelta{}{\mu_r}{}{\nu_r}{}\rho
{\hh\id\step\cd\step\id\Step\id\\
\id\step\id\step\hx\Step\id\\
\hh\id\step\lu\step\id\Step\id\\
\hh\id\Step\id\hstep\cd\step\cd\\
\cu\hstep\id\step\hx\step\id\\
\hh\step\id\step[1.5]\cu*\step\cu\\
\step\cu\Step\id}
{\cd\step\cd\\
\hh\id\step\ld\step\id\Step\id\\
\id\step\id\step\hx\Step\id\\
\hh\id\step\cu\step\id\Step\id}$
& \ref{alg-coalg-triv}.2, \ref{alg-coalg-triv}.8
& regular
\\
\hline$\divide\unitlens by 2\begin{tangle}\classbox 111000\end{tangle}$
&$\mDelta{}{}{\nu_l}{\nu_r}{}\rho
{\hh\id\step\cd\step\cd\step\id\\
\id\step\id\step\hx\step\id\step\id\\
\hh\id\step\lu\step\ru\step\id\\
\hh\id\Step\id\hstep\cd\step\cd\\
\cu\hstep\id\step\hx\step\id\\
\hh\step\id\step[1.5]\cu*\step\cu\\
\step\cu\Step\id}
{\hcd\step\hcd\\
\id\step\hx\step\id}$
& \ref{alg-coalg-triv}.3, \ref{alg-coalg-triv}.8
& pure
\\
\hline$\divide\unitlens by 2\begin{tangle}\classbox 001011\end{tangle}$
&$\mDelta{\mu_l}{\mu_r}{}{}{}{\rho}
{\hh\id\step\id\step\hcd\step\id\\
\hh\id\step\x\step\id\step\id\\
\hh\cu\step\ru\step\id\\
\hh\hstep\id\step\cd\step\cd\\
\hh\hstep\id\step\id\step\x\step\id\\
\hh\hstep\d\hstep\cu*\step\cu\\
\hh\step\cu\step[2]\id}
{\cd\step\cd\\
\hh\id\step\ld\step\rd\step\id\\
\id\step\id\step\hx\step\id\step\id\\
\hh\id\step\cu\step\cu\step\id}$
& \ref{alg-coalg-triv}.1, \ref{alg-coalg-triv}.2, \ref{alg-coalg-triv}.6
& 
\\
\hline$\divide\unitlens by 2\begin{tangle}\classbox 011001\end{tangle}$
&$\mDelta{\mu_l}{}{\nu_l}{}{}{\rho}
{\hh\id\step\id\step\hcd\step[1.5]\id\\
\hh\id\step\x\step\id\step[1.5]\id\\
\hh\hcu\step\ru\hstep{\hstr{300}\cd}\\
\hstep\id\step\hcd\step\rd\step\id\\
\hstep\id\step\id\step\hx\step\id\step\id\\
\hstep\id\step\hcu*\step\ru\step\id\\
\hstep{\hstr{225}\hcu}\step[1.5]\cu}
{\cd\step\cd \\
\hh\id\Step\id\step\rd\step\id \\
\id\Step\hx\step\id\step\id \\
\hh\id\Step\id\step\cu\step\id}$
& \ref{alg-coalg-triv}.1, \ref{alg-coalg-triv}.3, \ref{alg-coalg-triv}.6
& pure
\\
\hline$\divide\unitlens by 2\begin{tangle}\classbox 011010\end{tangle}$
&$\mDelta{}{\mu_r}{\nu_l}{}{}\rho
{\hh\id\step\id\step\hcd\step\id\\
\hh\id\step\x\step\id\step\id\\
\hh\cu\step\ru\step\id\\
\hh\hstep\id\step\cd\step\cd\\
\hh\hstep\id\step\id\step\x\step\id\\
\hh\hstep\d\hstep\cu*\step\cu\\
\hh\step\cu\step[2]\id}
{\cd\step\cd\\
\hh\id\step\ld\step\id\Step\id\\
\id\step\id\step\hx\Step\id\\
\hh\id\step\cu\step\id\Step\id}$
& \ref{alg-coalg-triv}.1, \ref{alg-coalg-triv}.8
& regular
\\
\hline$\divide\unitlens by 2\begin{tangle}\classbox 001010\end{tangle}$
&$\mDelta{}{\mu_r}{}{}{}\rho
{\hh\id\step\x\step[2]\id\\
\hh\cu\hstep\cd\step\cd\\
\hh\hstep\id\step\id\step\x\step\id\\
\hh\hstep\d\hstep\cu*\step\cu\\
\hh\step\cu\step[2]\id}
{\cd\step\cd\\
\hh\id\step\ld\step\id\Step\id\\
\id\step\id\step\hx\Step\id\\
\hh\id\step\cu\step\id\Step\id}$
& \ref{alg-coalg-triv}.1, \ref{alg-coalg-triv}.2, \ref{alg-coalg-triv}.8
& regular
\\
\hline$\divide\unitlens by 2\begin{tangle}\classbox 011000\end{tangle}$
&$\mDelta{}{}{\nu_l}{}{}\rho
{\hh\id\step\id\step\hcd\step\id\\
\hh\id\step\x\step\id\step\id\\
\hh\cu\step\ru\step\id\\
\hh\hstep\id\step\cd\step\cd\\
\hh\hstep\id\step\id\step\x\step\id\\
\hh\hstep\d\hstep\cu*\step\cu\\
\hh\step\cu\step[2]\id}
{\hcd\step\hcd\\
\id\step\hx\step\id}$
& \ref{alg-coalg-triv}.1, \ref{alg-coalg-triv}.3, \ref{alg-coalg-triv}.8
& pure
\\
\hline
\end{tabular}
\end{center}
\unitlens 12pt
\abs
There are various redundant and special relations among the defining
identities of the particular cases listed above. They can be derived
from Theorem \ref{hp-cp}. 

Exemplarily we will describe the explicit structure of
the cross product bialgebras of type
$\divide\unitlens by 3
\multiply\unitlens by 2
\begin{tangle}
 \hh\ffbox1{\normalfont\bullet}\ffbox1{\normalfont\bullet}\\
 \hh\ffbox1{\normalfont\blacktriangledown}\ffbox1{\normalfont\ }\\
 \hh\ffbox1{\normalfont\ }\ffbox1{\normalfont\bullet}
 \end{tangle}
\divide\unitlens by 2
\multiply\unitlens by 3\,$. They consist of two objects $B_1$ and $B_2$
such that
\begin{enumerate}
\item $(B_1,\m_1,\eta_1,\Delta_1,\varepsilon_1)$ is a bialgebra.
\item $(B_2,\Delta_2,\varepsilon_2,\mu_r,\nu_r)$ is 
      $B_1$-module coalgebra and $B_1$-comodule coalgebra and
      $\eta_2:\E\to B_2$ is coalgebra morphism.
\item
\begin{equation}\mathlabel{hp1-epsilon}
\begin{gathered}
\varepsilon_2\circ\m_2=\varepsilon_2\otimes\varepsilon_2\,,\quad
\m_2\circ(\id_{B_2}\otimes\eta_2)=\m_2\circ(\eta_2\otimes\id_{B_2})=
\id_{B_2}\,,\\
\mu_r\circ(\eta_2\otimes\id_{B_1})=
(\id_{B_2}\otimes\varepsilon_1)\circ\nu_l=\eta_2\circ\varepsilon_1\,,
\quad\nu_l\circ\eta_1=\eta_2\otimes\eta_1\,,\\
\mu_l\circ(\eta_2\otimes\id_{B_1})=\id_{B_1}\,,\quad
\mu_l\circ(\id_{B_2}\otimes\eta_1)=\eta_1\circ\varepsilon_2\,,\\
\varepsilon_1\circ\mu_l=\varepsilon_2\otimes\varepsilon_1\,,
\quad\nu_r\circ\eta_2=\eta_2\otimes\eta_1\,,\\
\sigma\circ(\eta_2\otimes\id_{B_2})=\sigma\circ(\id_{B_2}\otimes\eta_2)=
\eta_1\circ\varepsilon_2\,,\quad
\varepsilon_1\circ\sigma=\varepsilon_2\otimes\varepsilon_2\,,\\
(\id_{B_1}\otimes\varepsilon_1)\circ\rho=
(\varepsilon_1\otimes\id_{B_1})\circ\rho=\eta_1\circ\varepsilon_2\,,\quad
\rho\circ\eta_2=\eta_1\otimes\eta_1\,.
\end{gathered}
\end{equation}
\item The weak associativity of $\m_2$ and of $\mu_l$ in Definition
      \ref{hp} hold.
\item The module-algebra compatibility of Definition \ref{hp} is satisfied.
\item The cocycle compatibility of $\hat\sigma$ holds.
\item The algebra-coalgebra compatibility of $B_2$ holds.
\item The module-coalgebra compatibility of $\mu_l$, the
      comodule-algebra compatibility of $\nu_r$, the
      module-comodule compatibility and the cycle-cocycle compatibility
      are respectively given by
\begin{equation*}
\begin{split}
&\Delta_1\circ\mu_l=
(\mu_l\otimes\m_1\circ(\id_{B_1}\otimes\mu_l))\circ
(\id_{B_2}\otimes\Psi_{B_1\otimes B_2,B_1}\otimes\id_{B_1})\circ\\
&\quad\circ((\nu_r\otimes\id_{B_2})\circ\Delta_2\otimes\Delta_1)\,,\\[5pt]
&(\id_{B_2}\otimes\m_1)\circ(\Psi_{B_1,B_2}\otimes\id_{B_1})\circ
(\id_{B_1}\otimes\nu_r)\circ\hat\sigma\\
&=(\m_2\otimes\m_1\circ(\m_1\otimes\sigma)\circ
(\id_{B_1}\otimes\varphi_{2,1}\otimes\id_{B_2}))\circ
(\id_{B_2}\otimes\Psi_{B_1\otimes B_2,B_2}\otimes\id_{B_1\otimes B_2})\circ\\
&\quad\circ
((\nu_r\otimes\id_{B_2})\circ\Delta_2\otimes(\nu_r\otimes\id_{B_2})\circ
\Delta_2)\,,\\[5pt]
&(\id_{B_2}\otimes\m_1)\circ(\Psi_{B_1,B_2}\otimes\id_{B_1})\circ
(\id_{B_1}\otimes\nu_r)\circ\varphi_{2,1}\\
&=(\mu_r\otimes\m_1\circ(\id_{B_1}\otimes\mu_l))\circ
(\id_{B_2}\otimes\Psi_{B_1\otimes B_2,B_1}\otimes\id_{B_1})\circ
((\nu_r\otimes\id_{B_2})\circ\Delta_2\otimes\Delta_1)\,,\\[5pt]
&\Delta_1\circ\sigma=
(\sigma\otimes\m_1\circ(\m_1\otimes\sigma)\circ
(\id_{B_1}\otimes\varphi_{2,1}\otimes\id_{B_2}))\circ
(\id_{B_2}\otimes\Psi_{B_1\otimes B_2,B_2}\otimes\id_{B_1\otimes B_2})\circ\\
&\quad\circ
((\nu_r\otimes\id_{B_2})\circ\Delta_2\otimes(\nu_r\otimes\id_{B_2})\circ
\Delta_2)\,.
\end{split}
\end{equation*}
\end{enumerate}
The set of defining identities of the idempotents $\Pi_1$ and
$\Pi_2$ and of the projections and injections $\proj_1$, $\proj_2$,
$\inj_1$, $\inj_2$ can be determined similarly. Eventually the 
universal properties of the cross product bialgebras of type
$\divide\unitlens by 3
\multiply\unitlens by 2
\begin{tangle}
 \hh\ffbox1{\normalfont\bullet}\ffbox1{\normalfont\bullet}\\
 \hh\ffbox1{\normalfont\blacktriangledown}\ffbox1{\normalfont\ }\\
 \hh\ffbox1{\normalfont\ }\ffbox1{\normalfont\bullet}
 \end{tangle}
\divide\unitlens by 2
\multiply\unitlens by 3$
will be described by the following proposition. We will use the
notations of Corollary \ref{co-act-inv} and Remark \ref{non-except-bialg}.

\begin{proposition}\Label{bialg-cocycle6}
Let $B$ be a bialgebra in $\C$. Then the following equivalent conditions
are satisfied.

\begin{enumerate}
\item
$B$ is isomorphic to a cross product
bialgebra $B_1\bowtie B_2$ where
$\Delta_{21,1}=\eta_2\otimes\id_{B_1}$
and $\Delta_{11,2}= (\eta_1\otimes\eta_1)\circ\varepsilon_2$ are trivial.
\item
There are idempotents $\Pi_1,\Pi_2\in\End(B)$ such that
\begin{enumerate}
\item
$\m_B\circ(\Pi_1\otimes\Pi_1)=\Pi_1\circ\m_B\circ(\Pi_1\otimes\Pi_1)$,
\item
$\Pi_1$ is coalgebra morphism,
\item
$(\Pi_2\otimes\Pi_2)\circ\Delta_B=
(\Pi_2\otimes\Pi_2)\circ\Delta_B\circ\Pi_2$,
\item
$\Pi_1\circ\eta_B=\eta_B$ and $\varepsilon_B\circ\Pi_2=\varepsilon_B$,
\item
$\m_B\circ(\Pi_1\otimes\Pi_2)$ and $(\Pi_1\otimes\Pi_2)\circ\Delta_B$
split the idempotent $\Pi_1\otimes\Pi_2$.
\end{enumerate}
\item
There are objects $B_1$ and $B_2$ and
morphisms $B_1\overset{\inj_1}\longrightarrow A\overset{\proj_1}
\longrightarrow B_1$ and $B_2\overset{\inj_2}\longrightarrow
A\overset{\proj_2}\longrightarrow B_2$ such that
\begin{enumerate}
\item
$\inj_1$ is algebra and coalgebra morphism,
\item
$\proj_1$ is coalgebra morphism,
\item
$\proj_2$ is coalgebra morphism,
\item
$\proj_j\circ\inj_j=\id_{B_j}$ for $j\in\{1,2\}$.
\item
$\m_A\circ(\inj_1\otimes\inj_2):B_1\otimes B_2\to A$ is isomorphism
with inverse
$(\proj_1\otimes\proj_2)\circ\Delta_A$.\endproof
\end{enumerate}
\end{enumerate}
\end{proposition}

\section*{Concluding Remarks}

We found a universal theory of strong cross product
bialgebras with an equivalent (co\n-)modular co-cyclic
characterization in terms of strong Hopf data.
The theory unites all known cross product bialgebras
\cite{Rad1:85,Ma1:90,MS:94,Ma6:94,Ma1:95,BD1:98} in a
single construction. Furthermore various new types of
cross product bialgebras arise out of the most general
construction. The (co\n-)modular co-cyclic structure of strong
cross product bialgebras corresponds canonically to a strong Hopf
datum which in turn completely determines the bialgebra structure.

Thus Hopf data basically provide a pattern for the
realization of cross product bialgebras in terms of
explicit examples - a task to be done in future
investigations. 

There is no conceptual explanation yet 
for the understanding of the strong conditions in Definition
\ref{strong-hp}. A canonical origin of these conditions may be found
in higher dimensional categorical constructions of cross product bialgebras
where the two tensor factors of the strong cross product bialgebras
will be considered as object of two different monoidal categories. 

Since we are working throughout
in braided categories, the results of the article may now be used to 
investigate cross product bialgebras in various types of braided categories
(see \cite{BD1:98,BD3:98} for applications in Hopf bimodule categories).

Our directions of study of cross product bialgebras are particularly 
concerned with these questions as well as with extension theory
and cohomological considerations (see \cite{Hof1:90,Sin1:72,Swe1:68}).

\section{Proof of Theorem \ref{hp-cp}}\Label{proof-hp-cp}

This section is exclusively devoted to the proof of
Theorem \ref{hp-cp}. In what follows we denote by
$\mathfrak{h}=\big((B_1,\m_1,\eta_1,\Delta_1,\varepsilon_1),
(B_2,\m_2,\eta_2,\Delta_2,\varepsilon_2);\mu_l,\mu_r,\nu_l,\nu_r,
\rho,\sigma\big)$ a strong Hopf datum. Although many of the subsequent results
hold for more general Hopf data, we do not explicitely point out this fact
in the particular lemmas and propositions. The reader will easily
verify which of the following results hold under more general assumptions.
Before we start proving the theorem we would like to explain some
useful notations and results on Hopf data which will be used subsequently.

\subsection{Basic Properties of Strong Hopf Data}

For a finite set $I=\{i_1,\dots,i_r\}$ of indices
$i_k\in\{1,2\}$, $k\in\{1,\dots, r\}$ we denote
$B_I:=B_{i_1}\otimes\cdots\otimes B_{i_r}$.
Suppose now $r=\left\vert I\right\vert$ is the lenght of the sequence $I$.
Given morphisms $f:B_I\otimes B_2\otimes B_J\to B_K$ and
$g:B_K\to B_I\otimes B_2\otimes B_J$
we define their \emph{relativization} $f^{[r+1]}$ and $g_{[r+1]}$
in the $(r+1)$st domaine index and codomaine index respectively by
\begin{equation}\mathlabel{rel-morph1}
f^{[r+1]}:=\,
\divide\unitlens by 3
\multiply\unitlens by 2
\begin{tangle}
        \hstep\obj{B_I}\step[2]\obj{B_1}\step[1.5]\obj{B_J}\\
        \hh\hstep\id\step{\hstr{200}\hld}\step\id\\
        \hh\hstep\id\step\id\step\x\\
        \hh\ffbox3{f}\hstep\id\\
        \hh\step[1.5]\id\Step\id\\
        \step[1.5]\obj{B_K}\Step\obj{B_1}
\end{tangle}
\quad\text{and}\quad
g_{[r+1]}:=\,
\begin{tangle}
        \obj{B_2}\Step\obj{B_K}\\
        \hh\id\Step\id\\
        \hh\id\hstep\ffbox3{g}\\
        \hh\x\step\id\step\id\\
        \hh\id\step{\hstr{200}\hru}\step\id\\
        \obj{B_I}\step[1.5]\obj{B_2}\step[1.5]\obj{B_J}
\end{tangle}
\end{equation}
We say that two morphisms $f$ and $f'$
coincide relatively if $f^{[r+1]}=f'{}^{[r+1]}$. Similarly
$g_{[r+1]}=g'_{[r+1]}$ means that $g$ and $g'$ coincide relatively.
Instead of \eqref{rel-morph1} we will often use the obvious shorthand notation
\begin{equation}\mathlabel{rel-morph2}
\divide\unitlens by 3
\multiply\unitlens by 2
\begin{tangle}
        \hh\hstep\id\step\hld\hstep\id\\
        \hh\hstep\id\step\id\step\id\\
        \hh\ffbox3{f}\\
        \hh\step[1.5]\id
\end{tangle}
\,=\,
\begin{tangle}
        \hh\hstep\id\step\hld\hstep\id\\
        \hh\hstep\id\step\id\step\id\\
        \hh\ffbox3{f^\prime}\\
        \hh\step[1.5]\id
\end{tangle}
\quad\text{and}\quad
\begin{tangle}
        \hh\step[1.5]\id\\
        \hh\ffbox3{g}\\
        \hh\hstep\id\step\id\step\id\\
        \hh\hstep\id\hstep\hru\step\id
\end{tangle}
\,=\,
\begin{tangle}
        \hh\step[1.5]\id\\
        \hh\ffbox3{g^\prime}\\
        \hh\hstep\id\step\id\step\id\\
        \hh\hstep\id\hstep\hru\step\id
\end{tangle}
\end{equation}
for $f^{[r+1]}=f'{}^{[r+1]}$ and $g_{[r+1]}=g'_{[r+1]}$ respectively.
This means that the identities \eqref{rel-morph1} result from \eqref{rel-morph2}
by braiding the threads with endpoints in the middle of the graphics with
all neighbouring strings on the right, respectively on the left and  
then completing vertically these threads to the bottom, respectively to the
top of the graphics.

Besides \eqref{phi-sigma-rho2} and \eqref{phi-sigma-rho3} we will use the definitions
\unitlens 15pt                                              %
\begin{equation}\mathlabel{phi-sigma-rho4}
\begin{array}{c}
\varphi_{1,1}
:=
\divide\unitlens by 2
\begin{tangle}\ox{11}\end{tangle}
\multiply\unitlens by 2
:=
\divide\unitlens by 3
\begin{tangle}
\ld\step[2]\id\\\id\step\x\\\lu\step[2]\id
\end{tangle}
\multiply\unitlens by 3
\quad ,\quad
\varphi_{2,2}
:=
\divide\unitlens by 2
\begin{tangle}\ox{22}\end{tangle}
\multiply\unitlens by 2
:=
\divide\unitlens by 3
\begin{tangle}
\id\step[2]\rd\\
\x\step\id\\
\id\step[2]\ru
\end{tangle}
\end{array}
\end{equation}
All subsequent lemmas can be proven straighforwardly with the help of the
definition of (strong) Hopf data. We will therefore only sketch the main
steps of the derivation of the proofs.

\begin{lemma}\Label{weak-strong}
Let $\mathfrak{h}$ be a strong Hopf datum. Then the identities
\begin{equation}\mathlabel{cocycle-triv2}
\divide\unitlens by 3
\multiply\unitlens by 2
\begin{tangle}
\hh\hstep\hld\hstep\id\\
\hh\hstep\id\step\id\\
\hh\ffbox 2{\hat\sigma}\\
\hh\hstep\id\step\id
\end{tangle}
=
\begin{tangle}
\hh\step\hld\hstep\id\\
\hh\unit\step\cu
\end{tangle}
\quad , \quad
\begin{tangle}
\hstep\id\step\ld\\
\hh\ffbox 2{\hat\sigma}\hstep\id\\
\hh\hstep\id\step\id\step\id
\end{tangle}
=
\begin{tangle}
\unit\step\id\step\ld\\
\id\step\hcu\step\id
\end{tangle}
\quad , \quad
\begin{tangle}
\hh\hstep\id\step\id\\
\hh\ffbox 2{\hat\rho}\\
\hh\hstep\id\step\id\\
\hh\hstep\id\hstep\hru
\end{tangle}
=
\begin{tangle}
\hh\cd\\
\hh\id\hstep\hru\step\counit
\end{tangle}
\quad ,
\quad
\begin{tangle}
\hh\id\step\id\step\id\hstep\\
\hh\id\hstep\ffbox 2{\hat\rho}\\
\ru\step\id
\end{tangle}
=
\begin{tangle}
\id\step\hcd\step\id\\
\ru\step\id\step\counit
\end{tangle}
\end{equation}
are satisifed.
\end{lemma}
\begin{proof}
The first identity in \eqref{cocycle-triv2} has been obtained
from \eqref{cocycle-triv1} and \eqref{act-coact-triv}.
With the help of \eqref{cocycle-triv1} and \eqref{strong-hopf-combo}
the second identity will be derived. Application of $\pi$-symmetry
completes the proof.\end{proof}
\abs
The relative associativity of $\m_2$ and $\mu_l$, and by $\pi$-symmetry
the relative coassociativity of $\Delta_2$ and $\nu_r$ will be shown
in the following lemma.

\begin{lemma}\Label{strong-hd3}
For a strong Hopf datum $\mathfrak{h}$ the identities
\begin{equation*}
\begin{split}
\vstretch 60
\divide\unitlens by 2
\begin{tangle}
\ld\step\id\step\id\\
\id\step[2]\hcu\\
\hstr{125}\cu
\end{tangle}
=
\begin{tangle}
\ld\step\id\step\id\\
\cu\step\id\\
\step\cu
\end{tangle}
\quad ,\quad
\begin{tangle}
\id\step\ld\step\id\\
\id\step\cu\\
\cu
\end{tangle}
&=
\vstretch 60
\divide\unitlens by 2
\begin{tangle}
\id\step\ld\step\id\\
\hcu\step[2]\id\\
\hstep\hstr{125}\cu
\end{tangle}
\quad ,\quad
\begin{tangle}
\id\step\id\step\ld\\
\id\step\hcu\\
\hstr{150}\hcu
\end{tangle}
=
\begin{tangle}
\id\step\id\step\ld\\
\hcu\step\id\\
\hstep\hstr{150}\hcu
\end{tangle}
\\
\vstretch 60
\divide\unitlens by 2
\begin{tangle}
\ld\step\id\step\id\\
\cu\step\id\\
\step\lu[2]
\end{tangle}
=
\begin{tangle}
\ld\step\id\step\id\\
\nw1\step\lu\\
\step\lu[2]
\end{tangle}
\quad &,\quad
\vstretch 60
\divide\unitlens by 2
\begin{tangle}
\id\step\ld\step\id\\
\hcu\step[2]\id\\
\hstep\hstr{125}\lu[2]
\end{tangle}
=
\begin{tangle}
\id\step\hld\step\id\\
\id\step\hstr{150}\lu\\
\hstr{125}\lu[2]
\end{tangle}
\end{split}
\end{equation*}
and the corresponding $\pi$-symmetric versions of the relative coassociativity
of $\Delta_1$ and $\nu_r$ hold.
\end{lemma}

\begin{proof}
The relative associativity of $\m_2$ and $\mu_l$ will be derived from the
respective weak associativity of Definition \ref{hp} taking into account
\eqref{cocycle-triv1} and \eqref{cocycle-triv2} of Lemma
\ref{weak-strong}.\end{proof}

\begin{remark}
{\normalfont For a Hopf data satisfying
$\divide\unitlens by 2\begin{tangle}
	\vstr{70}\hh\id\hstep\cu*\\
	\vstr{35}\ru
\end{tangle}
=
\begin{tangle}
\vstr{70}\id\hstep\counit\hstep\counit
\end{tangle}$\
and
$\divide\unitlens by 2\begin{tangle}
	\vstr{35}\hstep\ld\\
	\vstr{70}\hh\cd*\hstep\id
\end{tangle}
=
\begin{tangle}
\vstr{70}\unit\hstep\unit\hstep\id
\end{tangle}$ 
the multiplication $\m_2$ and the comultiplication $\Delta_1$
are (co\n-)associative.}
\end{remark}

\begin{lemma}\Label{strong-hd2}
The relative versions of the morphisms $\varphi_{1,2}$ and $\varphi_{2,1}$
are given by
\begin{equation}\mathlabel{rel-phi}
\begin{array}{c}
\divide\unitlens by 2
\begin{tangle}
\ox{12}\\
\vstr{50}\id\step[2]\id\\
\vstr{60}\id\step\ru
\end{tangle}
=
\vstretch 60
\begin{tangle}
\id\step[2]\rd\\
\x\step\id\\
\id\step[2]\hcu\\
\id\step[1.5]\ru
\end{tangle}
\quad\text{and}\quad
\vstretch 100
\begin{tangle}
\vstr{60}\ld\step\id\\
\vstr{50}\id\step[2]\id\\
\ox{21}
\end{tangle}
=
\vstretch 60
\begin{tangle}
\hstep\ld\step[1.5]\id\\
\hcd\step[2]\id\\
\id\step\x\\
\lu\step[2]\id
\end{tangle}
\end{array}
\end{equation}
\end{lemma}

\begin{proof}
Using that $(B_2,\mu_r)$ is a right module and applying Lemma \ref{weak-strong}
yields the first identity of \eqref{rel-phi}.
\end{proof}
\abs
Then the next lemma follows from the module-algebra and comodule-coalgebra
compatibilities of Definition \ref{hp}.

\begin{lemma}\Label{strong-hd4}
The identities
\begin{gather}\mathlabel{mod-alg-rel}
\divide\unitlens by 2
\begin{tangle}
        \hh\hld\step[1.5]\id\step\id\\
        \hh{\hstr{200}\cu}\step\id\\
        \hh\step\hstr{400}\hru
\end{tangle}
=
\begin{tangle}
        \hh\hld\step[1.5]\hru\\
        \cu
\end{tangle}
\quad ,
\quad
\begin{tangle}
        \hh\hstr{400}\hld\\
        \hh\id\step\hstr{200}\cd\\
        \hh\id\step\id\step[1.5]\hru
\end{tangle}
=
\begin{tangle}
        \hstep\cd\\
        \hh\hld\step[1.5]\hru
\end{tangle}
\\
\mathlabel{mod-alg-rel2}
\divide\unitlens by 2
\begin{tangle}
        \hh\id\step\hld\step\id\\
        \hh\cu\step[1.5]\id\\
        \hstep\Ru
\end{tangle}
\;=\;
\begin{tangle}
        \hh\id\step[1.5]\hld\step\id\\
        \hh\id\step\cd\step\id\\
        \hh\id\step\id\step\x\\
        \hh\id\step{\hstr{200}\hlu}\step\id\\
        \hh\Ru\step\id\\
        \hh\hstr{300}\cu
\end{tangle}
\quad ,
\quad
\begin{tangle}
        \Ld\\
        \hh\id\step[1.5]\cd\\
        \hh\id\step\hru\step\id
\end{tangle}
\;=\;
\begin{tangle}
        \hh\hstr{300}\cd\\
        \hh\id\step\hstr{200}\ld\\
        \hh\id\step{\hstr{200}\hrd}\step\id\\
        \hh\x\step\id\step\id\\
        \hh\id\step\cu\step\id\\
        \hh\id\step\hru\step[1.5]\id
\end{tangle}
\end{gather}
are satisfied for a strong Hopf datum $\mathfrak{h}$.
\end{lemma}

\begin{lemma}\Label{strong-hd5} Let $\mathfrak{h}$ be strong Hopf datum. Then
$(B_1,\m_1,\eta_1,\nu_l)$ is a left $B_2$-comodule algebra and
$(B_2,\Delta_2,\varepsilon_2,\mu_r)$ is a right $B_1$-module coalgebra.
\end{lemma}

\begin{proof}
We use Lemma \ref{weak-strong} and the second identity of the comodule-algebra
compatibility of Definition \ref{hp} to show that
$(B_1,m_1,\eta_1,\nu_l)$ is a left $B_2$-comodule algebra. In $\pi$-symmetric
manner it will be proven that $(B_2,\Delta_2,\varepsilon_2,\mu_r)$ is
a right $B_1$-module coalgebra.\end{proof}

\begin{lemma}\Label{strong-hd8}
For a strong Hopf datum $\mathfrak{h}$
the relativizations of the algebra-coalgebra compatibilities,
the (left) module-coalgebra compatibility, and the (right)
comodule-algebra compatibility are respectively given by
\begin{gather}\mathlabel{alg-coalg-rel}
\divide\unitlens by 2
\begin{tangle}
        \cu\\
        \hh\hstr{200}\cd\\
        \hh\id\step[1.5]\hru
\end{tangle}
=
\begin{tangle}
        \hh\cd\step\cd\\
        \hh\id\step\x\step\id\\
        \hh\cu\step\cu\\
        \hh\hstep\id\step[1.5]\hru
\end{tangle}
\quad ,\quad
\begin{tangle}
        \hh\hld\step[1.5]\id\\
        \hh\hstr{200}\cu\\
        \cd
\end{tangle}
=
\begin{tangle}
        \hh\hstep\hld\step[1.5]\id\\
        \hh\cd\step\cd\\
        \hh\id\step\x\step\id\\
        \hh\cu\step\cu\\
\end{tangle}
\\
\mathlabel{rel-mod-coalg}
\divide\unitlens by 2
\begin{tangle}
\lu\\
\cd\\
\hh\id\step[1.5]\hru
\end{tangle}
=
\begin{tangle}
\cd\step\hcd\\
\hh\rd\step\x\step\id\\
\hh\id\step\x\step\lu\\
\lu\step\cu\\
\hh\step\id\step[1.5]\hru
\end{tangle}
\quad ,\quad
\begin{tangle}
\hstep\lu\\
\hstep\cd\\
\hh\hru\step[2]\id
\end{tangle}
=
\begin{tangle}
\cd\step\hcd\\
\hh\rd\step\x\step\id\\
\hh\id\step\x\step\lu\\
\lu\step\cu\\
\hh\hstep\hru\step[2]\id
\end{tangle}
\\
\mathlabel{rel-comod-alg}
\divide\unitlens by 2
\begin{tangle}
\hh\hld\\
\cu\\
\step\rd
\end{tangle}
=
\begin{tangle}
\hh\step\hld\\
\cd\step\rd\\
\hh\rd\step\x\step\id\\
\hh\id\step\x\step\lu\\
\hcu\step\cu
\end{tangle}
\quad ,\quad
\begin{tangle}
\hh\step[2]\hld\\
\cu\\
\step\rd
\end{tangle}
=
\begin{tangle}
\hh\step[3]\hld\\
\cd\step\rd\\
\hh\rd\step\x\step\id\\
\hh\id\step\x\step\lu\\
\hcu\step\cu
\end{tangle}
\end{gather}
\end{lemma}

\begin{proof}
Since $(B_1,\m_1,\eta_1,\nu_l)$ is a is a left comodule algebra by Lemma
\ref{strong-hd5} the first identity in \eqref{phi-prod} follows from the
relative asociativity of $\m_2$ according to Lemma \ref{strong-hd3}.

The verification of \eqref{rel-mod-coalg} needs a little more
calculation. We prove the first identity of \eqref{rel-mod-coalg}.
The second one can be derived with similar techniques.
We start with the second module-coalgebra compatibility in Definition \ref{hp}.
We apply the relativization (with  $\mu_r$) to the second tensor factor
on both sides of the graphic. The left hand side of this relative
module-coalgebra compatibility yields the left hand side of the
first identity in \eqref{rel-mod-coalg} if we consecutively apply
modularity of $(B_2,\mu_r)$, the second relation of \eqref{cocycle-triv1} and
the second relation of \eqref{act-coact-triv}.
To obtain the right hand side of the first identity of \eqref{rel-mod-coalg}
we transform the right hand side of the relative module-coalgebra compatibility
using successively the third relation of \eqref{cocycle-triv1}, modularity
of $(B_2,\mu_r)$, the first equation of \eqref{rel-phi}, the modularity of $B_2$,
the second identity of \eqref{cocycle-triv1}, and eventually again modularity of
$B_2$.

All other identities can be derived similarly, in particular because of
$\pi$-symmetric reasons.
\end{proof}
\abs
In the subsequent Lemmas \ref{strong-hd6} and \ref{strong-hd7} the
entwining properties of the morphisms $\varphi_{1,1}$, $\varphi_{2,2}$,
$\varphi_{1,2}$ and $\varphi_{2,1}$ will be investigated.

\begin{lemma}\Label{strong-hd6}
The morphism $\varphi_{1,1}$ entwines with the
multiplication $\m_1$, and $\varphi_{2,2}$ entwines with the
comultiplication $\Delta_2$ according to
\begin{equation}\mathlabel{entw1}
\begin{split}
\divide\unitlens by 2
\begin{tangle}
\cu\step\id\\
\step\ox{11}
\end{tangle}
=
\begin{tangle}
\id\step[2]\ox{11}\\
\ox{11}\step[2]\id\\
\id\step[2]\cu
\end{tangle}
\quad &,\quad
\divide\unitlens by 2
\begin{tangle}
\id\step\cu\\
\ox{11}
\end{tangle}
=
\begin{tangle}
\ox{11}\step[2]\id\\
\id\step[2]\ox{11}\\
\cu\step[2]\id
\end{tangle}
\\
\divide\unitlens by 2
\begin{tangle}
\ox{22}\\
\id\step\cd
\end{tangle}
=
\begin{tangle}
\cd\step[2]\id\\
\id\step[2]\ox{22}\\
\ox{22}\step[2]\id
\end{tangle}
\quad &,\quad
\divide\unitlens by 2
\begin{tangle}
\step\ox{22}\\
\cd\step\id
\end{tangle}
=
\begin{tangle}
\id\step[2]\cd\\
\ox{22}\step[2]\id\\
\id\step[2]\ox{22}
\end{tangle}
\\
\divide\unitlens by 2
\begin{tangle}
\cu\step\id\\
\step\ox{22}
\end{tangle}
=
\begin{tangle}
\vstr{50}\id\step[2]\cd\step[2]\id\\
\id\step[2]\id\step[2]\ox{22}\\
\vstr{50}\id\step[2]\id\step[2]\rd\step\id\\
\vstr{50}\id\step[2]\x\step\id\step\id\\
\vstr{50}\x\step[2]\lu\step\id\\
\vstr{50}\id\step[2]{\hstr{300}\ru}\step\id\\
\vstr{50}\id\step[2]\hstr{200}\cu
\end{tangle}
=
\begin{tangle}
\id\step\id\step\Rd\\
\hh\id\step\x\Step\id\\
\Put(0,0)[lb]{\begin{tangle}\hh\x\\
 \hh\id\step\id\end{tangle}}\Step\ox{21}\\
\hh\id\step{\hstr{200}\hru}\Step\id\\
\id\step\hstr{150}\cu
\end{tangle}
\quad &,\quad
\divide\unitlens by 2
\begin{tangle}
\ox{11}\\
\id\step\cd
\end{tangle}
=
\begin{tangle}
\vstr{50}{\hstr{200}\cd}\step[2]\id\\
\vstr{50}\id\step{\hstr{300}\ld}\step[2]\id\\
\vstr{50}\id\step\rd\step[2]\x\\
\vstr{50}\id\step\id\step\x\step[2]\id\\
\vstr{50}\id\step\lu\step[2]\id\step[2]\id\\
\ox{11}\step[2]\id\step[2]\id\\
\vstr{50}\id\step[2]\cu\step[2]\id
\end{tangle}
=
\begin{tangle}
{\hstr{150}\cd}\step\id\\
\hh\id\Step{\hstr{200}\hld}\step\id\\
\ox{12}\step\Put(0,0)[lb]{\begin{tangle}\hh\id\step\id\\
 \hh\x\end{tangle}}\\
\hh\id\Step\x\step\id\\
\Lu\step\id\step\id
\end{tangle}
\end{split}
\end{equation}
\end{lemma}

\begin{proof}
For the proof of the first identity of Lemma \ref{strong-hd6} we use Lemma
\ref{strong-hd5} and the fourth equation of Lemma \ref{strong-hd3}. To verify
the second identity of \eqref{entw1}, the module-algebra compatibility for
Hopf data, Lemma \ref{strong-hd2} and the comodularity of $(B_1,\nu_l)$
have to be applied successively. Using the module-algebra compatibility
and the $\pi$-symmetric version of the fourth identity of Lemma \ref{strong-hd3}
yields the fifth identity of \eqref{entw1}. Simple calculations yield the sixth
identity of \eqref{entw1}. The remaining relations of the lemma
follow by $\pi$-symmetric reasoning.\end{proof}

\begin{lemma}\Label{strong-hd7}
The (relative) entwining identities for $\varphi_{1,2}$ and $\varphi_{2,1}$
are given by
\begin{gather}
\divide\unitlens by 3
\multiply\unitlens by 2
\hstretch 80 \vstretch 80
\mathlabel{phi-prod}
\begin{tangle}
	\ox{21}\\
	\hh\id\step[1.5]\cd
\end{tangle}
=
\begin{tangle}
	\hh{\hstr{160}\cd}\step\cd\\
	\hh\id\Step\x\step\id\\
	\ox{21}\step\ru
\end{tangle}
\quad ,
\quad
\begin{tangle}
	\hh\cu\step[1.5]\id\\
	\hstep\ox{12}
\end{tangle}
=
\begin{tangle}
	\ld\step\ox{12}\\
	\hh\id\step\x\Step\id\\
	\hh\cu\step\hstr{160}\cu
\end{tangle}
\\[5pt]
\mathlabel{phi-prod2}
\divide\unitlens by 3
\multiply\unitlens by 2
\hstretch 80 \vstretch 80
\begin{tangle}
        \hh\id\step[1.5]\cu\\
        \ox{21}
\end{tangle}
=
\begin{tangle}
        \ox{21}\Step\id\\
        \id\Step\ox{21}\\
        \hh\hstr{160}\cu\step\id
\end{tangle}
\quad ,\quad
\begin{tangle}
        \hstep\ox{12}\\
        \hh\cd\step[1.5]\id
\end{tangle}
=
\begin{tangle}
        \hh\id\Step\hstr{160}\cd\\
        \ox{12}\Step\id\\
        \id\Step\ox{12}
\end{tangle}
\\[5pt]
\mathlabel{phi-prod-rel}
\divide\unitlens by 3
\multiply\unitlens by 2
\hstretch 80 \vstretch 70
\begin{tangle}
        \hh\hld\step[1.5]\id\step\id\\
        \hh{\hstr{160}\cu}\step\id\\
        \step\ox{21}
\end{tangle}
=
\begin{tangle}
        \hh\hld\step[1.5]\id\Step\id\\
        \id\step[2]\ox{21}\\
        \ox{21}\step[2]\id\\
        \id\step[2]\cu
\end{tangle}
\quad ,
\quad
\begin{tangle}
        \ox{12}\\
        \hh\id\step\hstr{160}\cd\\
        \hh\id\step\id\step[1.5]\hru
\end{tangle}
=
\begin{tangle}
 \cd\step[2]\id\\
 \id\step[2]\ox{12}\\
 \ox{12}\step[2]\id\\
 \hh\id\Step\id\step[1.5]\hru
\end{tangle}
\\[5pt]
\mathlabel{phi-prod-rel2}
\divide\unitlens by 3
\multiply\unitlens by 2
\hstretch 80 \vstretch 80
\begin{tangle}
        \step\ox{21}\\
        \cd\step\id\\
        \hh\id\step[1.5]\hru\step\id
\end{tangle}
\enspace=\enspace
\begin{tangle}
        \hh\hstr{160}\cd\hstep\cd\\
        \hh{\hstr{160}\hrd}\step\x\Step\id\\
        \id\step\hx\step\ox{21}\\
        \hh\lu\step\cu\Step\id\\
        \hh\step\id\step\hru\step[2.5]\id
\end{tangle}
\quad ,
\quad
\begin{tangle}
        \hh\id\step\hld\step[1.5]\id\\
        \id\step\cu\\
        \ox{12}
\end{tangle}
\enspace=\enspace
\begin{tangle}
        \hh\id\step[2.5]\hld\step\id\\
        \hh\id\Step\cd\step\rd\\
        \ox{12}\step\hx\step\id\\
        \hh\id\Step\x\step{\hstr{160}\hlu}\\
        \hh\hstr{160}\cu\hstep\cu
\end{tangle}
\end{gather}
\end{lemma}

\begin{proof}
We use Lemma \ref{strong-hd5} for $B_1$ and the relative associativity of $\m_2$
(Lemma \ref{strong-hd3}) to obtain the first identity of \eqref{phi-prod}.

The first identity of \eqref{phi-prod2}
will be used explicitely in the proof of Proposition \ref{co-assoc}.
Below we will give its detailed derivation.
\unitlens 10pt                                              %
\begin{equation*}
\hstretch 80 \vstretch 80
\varphi_{2,1}\circ(B_2\otimes\m_1)=
\begin{tangle}
\hh\hstep\id\step\cd\step\cd\\
\hh\hstep\id\step\id\step\x\step\id\\
\hh\cd\hstep\cu\step\cu\\
\hh\id\step\x\Step\id\\
\lu\step\Ru\\
\end{tangle}
=
\begin{tangle}
\hh\hstep{\hstr{160}\cd}\step\cd\step\cd\\
\hh\hstep\id\Step\x\step\x\step\id\\
\hh\cd\step\cd\hstep\x\step\id\step\id\\
\hh\id\step\x\step\id\hstep\id\step\ru\step\id\\
\hh\lu\step\ru\hstep\id\step\Ru\\
\hh\step\id\step{\hstr{120}\lu}\step\id\\
\hh\step{\hstr{200}\cu}\step\id
\end{tangle}
=
\begin{tangle}
\hh\cd\step\cd\Step\id\\
\hh\id\step\x\step\nw 1\step\id\\
\hh\lu\hstep\cd\step\cd\hstep\id\\
\hh\step\id\hstep\id\step\x\step\id\hstep\id\\
\hh\step\id\hstep\ru\step\ru\hstep\id\\
\hh\step\id\hstep\id\Step\id\step\cd\\
\hh\step\id\hstep\nw 1\step\x\step\id\\
\hh\step\d\step\lu\step\ru\\
\step[1.5]\cu\step\id
\end{tangle}
=(\m_1\otimes B_2)\circ\varphi_{2,11}
\end{equation*}
\unitlens 15pt                                              %
where the first equation comes from the relative bialgebra property
of $B_1$ according to \eqref{alg-coalg-rel}.
The second equation can be verified with the help of the
module-algebra compatibility of Definition \ref{hp}, and using
that $(B_2,\mu_r)$ is a right module. In the third identity we
use the relative associativity of $\Delta_1$ proven in Lemma
\ref{strong-hd3}. Finally the result follows because
$(B_2,\Delta_2,\varepsilon_2,\mu_r)$ is a module coalgebra by Lemma
\ref{strong-hd5}.

The left hand side of the first identity of \eqref{phi-prod-rel}
will be transformed to the right hand side of the
identity by using successively the algebra-coalgebra compatibility for
$\Delta_2$ and $\m_2$, the comodularity of $(B_1,\nu_l)$, the first
relation of \eqref{act-coact-triv}, the first identity of \eqref{mod-alg-rel},
the relative associativity of $\mu_l$ according to Lemma
\ref{strong-hd3} and finally the second relation of \eqref{rel-phi}.

The first identity of \eqref{phi-prod-rel2} immediately follows from
the relative module coalgebra properties \eqref{rel-mod-coalg}.

All other relations of the lemma can be derived easily by
applying $\pi$-symmetry.\end{proof}

\begin{lemma}
For the strong Hopf datum $\mathfrak{h}$ the following identities are
satisfied.
\unitlens 15pt                                              %
\begin{gather}
\mathlabel{rel-cocycle-assoc1}
\divide\unitlens by 2
\def\FillCircDiam{2}
\begin{tangle}
        \hh\hld\step[1.5]\id\step\id\\
        \hh\hstr{200}\cu\hstep\id\\
        \hh\step\hstr{200}\cu*
\end{tangle}
=
\def\FillCircDiam{4}
\begin{tangle}
        \hh\hld\step\cu*\\
        \Lu
\end{tangle}
\quad ,\quad
\def\FillCircDiam{2}
\begin{tangle}
        \hh\hstr{200}\cd*\\
        \hh\id\step\hstr{200}\cd\\
        \hh\id\step\id\step[1.5]\hru
\end{tangle}
=
\def\FillCircDiam{4}
\begin{tangle}
        \hstep\Rd\\
        \hh\cd*\step\hru
\end{tangle}
\\
\mathlabel{rel-cocycle-assoc2}
\divide\unitlens by 2
\def\FillCircDiam{2}
\begin{tangle}
        \hh\id\step\hld\step\id\\
        \hh\cu\step[1.5]\id\\
        \hh\hstep\hstr{200}\cu*
\end{tangle}
=
\begin{tangle}
        \hh\id\step\hld\step[1.5]\id\\
        \hh\id\step\hstr{200}\cu\\
        \hh\hstr{200}\cu*
\end{tangle}
\quad ,\quad
\begin{tangle}
        \hh\hstr{200}\cd*\\
        \hh\id\step[1.5]\cd\\
        \hh\id\step\hru\step\id
\end{tangle}
=
\begin{tangle}
        \hh\step\hstr{200}\cd*\\
        \hh{\hstr{200}\cd}\step\id\\
        \hh\id\step[1.5]\hru\step\id
\end{tangle}
\\
\mathlabel{rel-cocycle-assoc3}
\divide\unitlens by 2
\def\FillCircDiam{4}
\begin{tangle}
        \hh\id\hstep\id\step\hld\\
        \hh\id\hstep\cu\hstep\id\\
        \hh\cu*\step\id
\end{tangle}
=
\begin{tangle}
       \cu*\step\id
\end{tangle}
\quad ,\quad
\begin{tangle}
        \hh\id\step\cd*\\
        \hh\id\hstep\cd\hstep\id\\
        \hh\hru\step\id\hstep\id
\end{tangle}
=
\begin{tangle}
        \id\step\cd*
\end{tangle}
\end{gather}
\end{lemma}
\begin{proof}
The proof follows straightforwardly from the cocycle and cycle
compatibilities of Definition \ref{hp} and the identities
\eqref{act-coact-triv} and \eqref{cocycle-triv1}.\end{proof}

\begin{lemma}
The subsequent relations involving $\hat\rho$ and $\hat\sigma$ hold
in $\mathfrak{h}$.
\begin{gather}
\mathlabel{hat-prod}
\divide\unitlens by 2
\vstretch 110
\begin{tangle}
	\hh\hstep\id\step\id\\
	\hh\ffbox2{\hat\sigma}\\
	\hh\hstep\id\hstep\cd
\end{tangle}
=
\begin{tangle}
	\hh\hstep\cd\Step\cd\\
	\hstep\id\step\ox{22}\step\id\\
	\hh\ffbox2{\hat\sigma}\step[1.5]\cu\\
	\hh\hstep\id\step\id\step[2.5]\id
\end{tangle}
\quad ,\quad
\begin{tangle}
	\hh\cu\hstep\id\\
	\hh\ffbox2{\hat\rho}\\
	\hh\hstep\id\step\id
\end{tangle}
=
\begin{tangle}
	\hh\hstep\id\step[2.5]\id\step\id\\
	\hh\cd\step[1.5]\ffbox2{\hat\rho}\\
	\id\step\ox{11}\step\id\\
	\hh\cu\Step\cu
\end{tangle}
\\
\mathlabel{hat-prod-rel}
\divide\unitlens by 2
\vstretch 110
\begin{tangle}
        \hh\hld\hstep\id\hstep\id\\
        \hh\cu\hstep\id\\
        \hh\ffbox2{\hat\sigma}\\
        \hh\hstep\id\step\id
\end{tangle}
=
\begin{tangle}
        \hh\hstep\hld\step\id\step\id\\
        \hh\cd\hstep\ffbox2{\hat\sigma}\\
        \hh\id\step\x\step\id\\
        \hh{\hstr{200}\hlu}\step\cu
\end{tangle}
\quad ,\quad
\begin{tangle}
        \hh\hstep\id\step\id\\
        \hh\ffbox2{\hat\rho}\\
        \hh\hstep\id\hstep\cd\\
        \hh\hstep\id\hstep\id\hstep\hru
\end{tangle}
=
\begin{tangle}
        \hh\hstep\cd\step\hstr{200}\hrd\\
        \hh\hstep\id\step\x\step\id\\
        \hh\ffbox2{\hat\rho}\hstep\cu\\
        \hh\hstep\id\step\id\step\hru
\end{tangle}
\end{gather}
\end{lemma}
\begin{proof} Because of $\pi$-symmetry we will only
demonstrate the first identities of \eqref{hat-prod} and \eqref{hat-prod-rel}.
Applying to the left hand side of \eqref{hat-prod}
the algebra-coalgebra compatibility (for $\Delta_2$ and $\m_2$)
and the entwining property of $\varphi_{2,2}$ (Lemma \ref{strong-hd6})
yields the result.

To obtain the first identity of \eqref{hat-prod-rel} we transform its left hand
side consecutively using the relative bialgebra property \eqref{alg-coalg-rel}
of $\Delta_2$ and $\m_2$, the comodularity of $(B_1,\nu_l)$, the first relation
of \eqref{mod-alg-rel}, the first identity of \eqref{rel-cocycle-assoc1}, the
relative associativity of $\m_2$, and again the fact that $(B_1,\nu_l)$ is left
comodule.\end{proof}

\begin{remark}{\normalfont Note that \eqref{cocycle-triv2} and
\eqref{rel-cocycle-assoc1} are special cases of \eqref{hat-prod-rel}.
Furthermore \eqref{rel-cocycle-assoc3} implies \eqref{cocycle-triv1}.}
\end{remark}

\subsection{Special Properties of Strong Hopf Data}

Before we will prove Theorem \ref{hp-cp} we have to provide several 
auxiliary definitions and specific properties of strong Hopf Data.
Similarly as in Remark \ref{alpha-ip} we define 
$\proj_1=\id_{B_1}\otimes\varepsilon_2$,
$\proj_2=\varepsilon_1\otimes\id_{B_2}$,
$\proj_0=\id_B$, $\inj_1=\id_{B_1}\otimes\eta_2$,
$\inj_2=\eta_1\otimes\id_{B_2}$, and $\inj_0=\id_B$. In this context we
occasionally use the notation $B_0:=B:=B_0\otimes B_1$. 
Given a strong Hopf datum $\mathfrak{h}$ we can build the morphisms
$\m_B:B\otimes B\to B$ and $\Delta_B:B\to B\otimes B$ like in
\eqref{type-alpha3} as
\begin{equation}\mathlabel{strong-mult-comult}
\begin{array}{c}
\m_B :=
\hstretch 150
\divide\unitlens by 3
\begin{tangle}
\hh \id\step\cd\step\cd\step\id \\
\id\step\id\step\hx\step\id\step\id \\
\hh \id\step\lu\step\ru\step\id \\
\hh \id\Step\id\hstep\cd\step\cd \\
\hh \id\Step\id\hstep\id\step\id\step\hrd\hstep\id \\
\cu\hstep\id\step\hx\hstep\id\hstep\id \\
\hh \step\id\step\hstep\id\step\id\step\hru\hstep\id \\
\hh \step\id\step[1.5]\cu*\step\cu \\
\step\cu\Step\id
\end{tangle}
\multiply\unitlens by 3
\qquad\text{and}\qquad
\Delta_B:=
\divide\unitlens by 3
\begin{tangle}
\hstep\id\Step\cd\\
\hh \cd\step\cd*\step[1.5]\id\\
\hh \id\hstep\hld\step\id\step\id\hstep\step\id\\
\id\hstep\id\hstep\hx\step\id\hstep\cd\\
\hh \id\hstep\hlu\step\id\step\id\hstep\id\Step\id\\
\hh \cu\step\cu\hstep\id\Step\id\\
\hh \hstep\id\step\ld\step\rd\step\id\\
\hstep\id\step\id\step\hx\step\id\step\id\\
\hh \hstep\id\step\cu\step\cu\step\id \end{tangle}
\multiply\unitlens by 3
\end{array}
\end{equation}
Similarly as in \eqref{m-delta-ijk} and \eqref{m-delta-ijk1} we define
\unitlens 12pt                                              %
\begin{equation}\mathlabel{aux-m-delta}
\begin{gathered}
\m_{i,jk} :=\proj_i\circ\m_B\circ(\inj_j\otimes\inj_k)=
\vstretch 80
\begin{tangle}
\tu{\sstyle i\vert jk}\\
\end{tangle}
\quad
\text{for all $i,j,k\in\{0,1,2\}$}\,,
\\
\Delta_{ij,k}:=(\proj_i\otimes\proj_j)\circ\Delta_B\circ\inj_k=
\vstretch 80
\begin{tangle}
\td{\sstyle ij\vert k}\\
\end{tangle}
\quad
\text{for all $i,j,k\in\{0,1,2\}$}\,.
\end{gathered}
\end{equation}
and
\begin{equation}\mathlabel{aux-m-delta1}
\begin{gathered}
\m^*_{0,20}:=\hat\sigma\circ(\mu_r\otimes\id_{B_2})\,,\quad
\m^*_{l,2m}:=\proj_l\circ\m^*_{0,20}\circ(\id_{B_2}\otimes\inj_m)
\\
\Delta^*_{01,0}:=(\id_{B_1}\otimes\nu_l)\circ\hat\rho\,,\quad
\Delta^*_{l1,m}:=(\proj_l\otimes\id_{B_1})\circ\Delta^*_{01,0}\circ\inj_m
\end{gathered}
\end{equation}
for $l,m\in\{0,1,2\}$. In particular \eqref{aux-m-delta} implies
\begin{equation*}
\begin{split}
\m_1&=\m_{1,11}\,,\\
\m_2&=\m_{2,22}\,,\\
\m_B&=\m_{0,00}\,,
\end{split}
\quad
\begin{split}
\Delta_1&=\Delta_{11,1}\,,\\
\Delta_2&=\Delta_{22,2}\,,\\
\Delta_B&=\Delta_{00,0}\,,
\end{split}
\quad
\begin{split}
\mu_l&=\m_{1,21}\,,\\
\mu_r&=\m_{2,21}\,,\\
\sigma&=\m_{1,22}\,,
\end{split}
\quad
\begin{split}
\nu_l&=\Delta_{21,1}\,,\\
\nu_r&=\Delta_{21,2}\,,\\
\rho&=\Delta_{11,2}\,.
\end{split}
\end{equation*}
\abs
The following morphisms
$\ixi rsklijmn,\iy rskli,\yi rsklj, \iy rskl2^*,
\yi rskl1^* :B_k\otimes B_l\to B_r\otimes B_s$
turn out to be useful in the sequel. Let $i,j,m,n,r,s,k,l\in\{0,1,2\}$, then
\begin{equation}\mathlabel{x-def}
\begin{split}
\ixi rsklijmn &:= (\m_{r,im}\otimes\m_{s,jn})\circ
                   (\id_{B_i}\otimes\Psi_{B_j,B_m}\otimes\id_{B_n})\circ
                   (\Delta_{ij,k}\otimes\Delta_{mn,l})\,,\\[3pt]
\iy rskli &:=(\m_{r,il}\otimes\id_{B_s})\circ
	       (\id_{B_i}\otimes\Psi_{B_s,B_l})\circ
	       (\Delta_{is,k}\otimes\id_{B_l})\,,\\
\yi rsklj &:=(\id_{B_r}\otimes\m_{s,kj})\circ
	       (\Psi_{B_k,B_r}\otimes\id_{B_j})\circ
	       (\id_{B_k}\otimes\Delta_{rj,l})\,,\\[3pt]
\iy rskl2^* &:=(\m^*_{r,2l}\otimes\id_{B_s})\circ
	       (\id_{B_2}\otimes\Psi_{B_s,B_l})\circ
	       (\Delta_{2s,k}\otimes\id_{B_l})\,,\\
\yi rskl1^* &:=(\id_{B_r}\otimes\m_{s,k1})\circ
	       (\Psi_{B_k,B_r}\otimes\id_{B_1})\circ
	       (\id_{B_k}\otimes\Delta^*_{r1,l})\,.
\end{split}
\end{equation}
For a strong Hopf datum $\mathfrak{h}$ we obtain the following bra-ket
decomposition of morphisms \eqref{x-def}.

\begin{lemma}\Label{bra-ket}
\begin{equation}\mathlabel{bra-ket-decomp}
\begin{gathered}
\ixi rskli0mn = (\id_{B_r}\otimes\m_{s,1s})\circ
		(\iy r1kmi\otimes\id_{B_s}) \circ
                (\id_{B_k}\otimes\yi ms2ln )\circ
		(\Delta_{k2,k}\otimes\id_{B_l})
\end{gathered}
\end{equation}
for any $i,k,l,m,n,r,s\in\{0,1,2\}$.
\end{lemma}
\begin{proof}
The lemma follows straightforwardly from the
(co\n-)associativity of $\m_B$ and $\Delta_B$ and the identities
$\m_{0,12}=\id_{B_0}=\Delta_{12,0}$.\end{proof}
\abs
Below we define another set of auxiliary morphisms which will be used
in the course of the proof of Theorem \ref{hp-cp}.

\begin{gather}\mathlabel{aux-tau-t}
\divide\unitlens by 2
\tau^0_t:=
\begin{tangle}
	\hh\cd\step\cd\step\id\\
	\hh\id\step\x\step\id\step\id\\
	\hh\id\step\id\step{\hstr{200}\hru}\step\id\\
	\id\step\id\step\ox{\sstyle 22}
\end{tangle}
\enspace ,
\enspace
\tau^1_t:=
\begin{tangle}
	\hh\cd\step\cd\step[1.5]\id\\
	\hh\id\step\x\step\id\step[1.5]\id\\
	\hh\id\step\id\step{\hstr{200}\hru}\step\cd\\
	\id\step\id\step\ox{\sstyle 22}\step\id\\
	\hh\id\step\id\step\id\Step\cu\\
\end{tangle}
\enspace ,
\enspace
\tau^2_t:=
\def\FillCircDiam{5}
\begin{tangle}
	\hh\cd\step\cd\step[1.5]\id\\
	\hh\id\step\x\step\id\step[1.5]\id\\
	\hh\id\step\id\step{\hstr{200}\hru}\step\cd\\
	\id\step\id\step\ox{\sstyle 22}\step\id\\
	\hh\id\step\id\step\id\Step\cu*\\
\end{tangle}
\enspace ,
\enspace
\tau^{3}_t:=
\begin{tangle}
	\hh\cd\step\cd\step[1.5]\id\\
	\hh\id\step\x\step\id\step[1.5]\id\\
        \hh\id\step\id\step{\hstr{200}\hru}\step\cd\\
        \id\step\id\step\ox{\sstyle 22}\step\rd\\
        \hh\id\step\id\step\id\Step\x\step\id\\
        \hh\id\step\id\step\id\step[1.5]\cd\hstep\hstr{200}\hlu\\
\end{tangle}
\\ 
\mathlabel{aux-tau-b}
\divide\unitlens by 2
\tau^0_b:=
\begin{tangle}
	\ox{\sstyle 11}\step\id\step\id\\
	\hh\id\step{\hstr{200}\hld}\step\id\step\id\\
	\hh\id\step\id\step\x\step\id\\
	\hh\id\step\cu\step\cu\\
\end{tangle}
\enspace ,
\enspace
\tau^1_b:=
\begin{tangle}
	\hh\cd\Step\id\step\id\step\id\\
	\id\step\ox{\sstyle 11}\step\id\step\id\\
	\hh\cu\step{\hstr{200}\hld}\step\id\step\id\\
	\hh\hstep\id\step[1.5]\id\step\x\step\id\\
	\hh\hstep\id\step[1.5]\cu\step\cu\\
\end{tangle}
\enspace ,
\enspace
\tau^2_b:=
\def\FillCircDiam{5}
\begin{tangle}
	\hh\cd*\Step\id\step\id\step\id\\
	\id\step\ox{\sstyle 11}\step\id\step\id\\
	\hh\cu\step{\hstr{200}\hld}\step\id\step\id\\
	\hh\hstep\id\step[1.5]\id\step\x\step\id\\
	\hh\hstep\id\step[1.5]\cu\step\cu\\
\end{tangle}
\enspace ,
\enspace
\tau^{3}_b:=
\begin{tangle}
        \hh{\hstr{200}\hrd}\hstep\cu\step[1.5]\id\step\id\step\id\\
        \hh\id\step\x\Step\id\step\id\step\id\\
        \lu\step\ox{\sstyle 11}\step\id\step\id\\
        \hh\step\cu\step{\hstr{200}\hld}\step\id\step\id\\
        \hh\step[1.5]\id\step[1.5]\id\step\x\step\id\\
        \hh\step[1.5]\id\step[1.5]\cu\step\cu\\
\end{tangle}
\end{gather}
Besides we use the following definitions.
\begin{equation}\mathlabel{aux-tr-red}
\begin{gathered}
\sigma^{\mathrm{tr}}:=
\divide\unitlens by 2
\def\FillCircDiam{5}
\begin{tangle}
\hh\cd\Step\cd\\
\id\step\ox{\sstyle 22}\step\id\\
\hh\id\step\id\Step\cu*
\end{tangle}
\enspace ,
\enspace
\rho^{\mathrm{tr}}:=
\begin{tangle}
\hh\cd*\Step\id\step\id\\
\id\step\ox{\sstyle 11}\step\id\\
\hh\cu\Step\cu
\end{tangle}
\enspace ,
\enspace
(\tau^3_t)^{\mathrm{red}_1}:=
\begin{tangle}
	\hh\cd\step\cd\step[1.5]\id\\
	\hh\id\step\x\step\id\step[1.5]\id\\
        \hh\id\step\id\step{\hstr{200}\hru}\step\cd\\
        \id\step\id\step\ox{\sstyle 22}\step\rd\\
        \hh\id\step\id\step\id\Step\x\step\id\\
        \hh\id\step\id\step\id\Step\id\step\hstr{200}\hlu\\
\end{tangle}
\enspace ,
\enspace
(\tau^3_b)^{\mathrm{red}_1}:=
\begin{tangle}
        \hh{\hstr{200}\hrd}\step\id\Step\id\step\id\step\id\\
        \hh\id\step\x\Step\id\step\id\step\id\\
        \lu\step\ox{\sstyle 11}\step\id\step\id\\
        \hh\step\cu\step{\hstr{200}\hld}\step\id\step\id\\
        \hh\step[1.5]\id\step[1.5]\id\step\x\step\id\\
        \hh\step[1.5]\id\step[1.5]\cu\step\cu\\
\end{tangle}
\\
\divide\unitlens by 2
(\tau^3_t)^{\mathrm{red}_2}:=
\begin{tangle}
\hh\cd\Step\cd\\
\id\step\ox{\sstyle 22}\step\rd\\
\hh\id\step\id\Step\x\step\id\\
\hh\id\step\id\Step\id\step\hstr{200}\hlu
\end{tangle}
\enspace ,
\enspace
(\tau^3_b)^{\mathrm{red}_2}:=
\begin{tangle}
\hh{\hstr{200}\hrd}\step\id\Step\id\step\id\\
\hh\id\step\x\Step\id\step\id\\
\lu\step\ox{\sstyle 11}\step\id\\
\hh\step\cu\Step\cu
\end{tangle}
\enspace ,
\enspace
(\tau^3_t)^{\mathrm{red}_3}:=
\begin{tangle}
\hh\cd\step\cd\step\id\\
\hh\id\step\x\step\id\step\id\\
\hh\id\step\id\step{\hstr{200}\hru}\step\rd\\
\vstr{50}\id\step\id\step\x\step\id\\
\hh\id\step\id\step\id\Step\hstr{200}\hlu
\end{tangle}
\enspace ,
\enspace
(\tau^3_b)^{\mathrm{red}_3}:=
\begin{tangle}
\hh{\hstr{200}\hrd}\Step\id\step\id\step\id\\
\vstr{50}\id\step\x\step\id\step\id\\
\hh\lu\step{\hstr{200}\hld}\step\id\step\id\\
\hh\step\id\step\id\step\x\step\id\\
\hh\step\id\step\cu\step\cu
\end{tangle}
\end{gathered}
\end{equation}
Observe that the morphisms $(\tau_t)^{\mathrm{red}}$ and
$(\tau_b)^{\mathrm{red}}$ in \eqref{aux-tr-red}
can be expressed with the help $\tau^3_t$ and $\tau^3_b$
respectively. For instance
$(\tau^3_t)^{\mathrm{red}_1}=(\id_{B_2}\otimes\id_{B_1}\otimes
\id_{B_2}\otimes\varepsilon_2\otimes\id_{B_2}\otimes\id_{B_1})\circ\tau^3_t$.

\begin{lemma}\Label{tau-delta}
For a strong Hopf datum $\mathfrak{h}$ the following reduction
identities are satisfied.
\begin{gather}
\mathlabel{tau-R0-aux}
\yi i0201=(\yi i1201\otimes\id_{B_2})\circ\tau^0_t\,,
\quad
\iy 010i1=\tau^0_b\circ(\id_{B_1}\otimes\iy 012i1)\,,
\\[5pt]
\mathlabel{tau-R-aux}
\yi i0200=(\yi i1200\otimes\id_{B_2})\circ\tau^1_t\,,
\quad
\iy 010i0=\tau^1_b\circ(\id_{B_1}\otimes\iy 012i0)\,,
\\[5pt]
\mathlabel{tau-R2-aux}
\yi i1200
=(\id_B\otimes\m_1)\circ(\yi i1201\otimes\id_{B_1})\circ\tau^2_t\,,
\quad
\iy 012i0
=\tau^2_b\circ(\id_{B_2}\otimes\iy 012i2)\circ(\Delta_2\otimes\id_B)\,,
\\[5pt]
\mathlabel{tau-check-aux}
\yi i1220
=(\id_B\otimes\m_1)\circ(\yi i1221\otimes\id_{B_1})\circ\sigma^{\mathrm{tr}}\,,
\enspace
\iy 112i0
=\rho^{\mathrm{tr}}\circ(\id_{B_2}\otimes\iy 112i2)\circ(\Delta_2\otimes\id_B)\,.
\end{gather}
\end{lemma}

\begin{proof}
We only consider the first identities of
\eqref{tau-R0-aux} and \eqref{tau-R-aux}.
The second identities in each row are $\pi$-symmetric analogues.
The identities \eqref{tau-R2-aux} follow from 
\eqref{tau-R0-aux} by application of the mapping
$f\mapsto(\id\otimes\m_{1,02})\circ(f\otimes\id_{B_2})
\circ(\id\otimes\Delta_2)$. Relations \eqref{tau-check-aux} are
special cases of \eqref{tau-R2-aux} through composition with $\inj_2$.
In a first step we derive \eqref{tau-R0-aux} and \eqref{tau-R-aux}
for the case i=2. For \eqref{tau-R0-aux} we obtain
\begin{equation}\mathlabel{phi-phi}
\yi 20201 =
\divide\unitlens by 3
\multiply\unitlens by 2
\begin{tangle}
	\id\step\ox{\sstyle 12}\\
	\hh\x\Step\id\\
	\id\step\ox{\sstyle 21}
\end{tangle}
=
\begin{tangle}
\hh\step\id\step{\hstr{200}\cd\step\id}\\
\step\id\step\id\Step\ox{\sstyle 12}\\
\step\id\step\ox{\sstyle 12}\Step\id\\
\hh\step\x\Step\id\Step\id\\
\hh\sw1{\hstr{200}\cd}\step\id\Step\id\\
\hh\id\step\id\Step\x\Step\id\\
\hh\id\step{\hstr{400}\hlu}\step{\hstr{400}\hru}
\end{tangle}
=
\begin{tangle}
\vstretch 110
\vstr{55}\hcd\step[2]\hcd\step\id\\
\vstr{55}\id\step\x\step\id\step\id\\
\vstr{55}\id\step\id\step[2]\ru\step\id\\
\id\step\id\Step\ox{\sstyle 22}\\
\id\step\ox{\sstyle 12}\Step\id\\
	\hh\x\Step\id\Step\id\\
	\vstr{44}\id\step\Lu\Step\id
\end{tangle}
\end{equation}
where the second equation is an immediate consequence of
the relative entwining property \eqref{phi-prod-rel} of $\varphi_{1,2}$
and the third equation follows from \eqref{rel-phi}
of Lemma~\ref{strong-hd2} and the right module property of 
$(B_2,\mu_r)$.

To derive the first identiy of \eqref{tau-R-aux} for i=2
we use \eqref{phi-prod}, the relative entwining property 
\eqref{phi-prod-rel} of $\varphi_{1,2}$, \eqref{act-coact-triv} 
of Definition \ref{strong-hp}, the coassociativity of $\Delta_2$ and the
entwining property \eqref{entw1} of $\varphi_{2,2}$.

Using \eqref{hat-prod-rel}, coassociativity of $\Delta_2$ and
the entwining property \eqref{entw1} of $\varphi_{2,2}$ we obtain
the relations $(\id_{B_1}\otimes\tau^j_t)\circ f=
(f\otimes\id_{B_2})\circ\tau^j_t\,$ for $j\in\{0,1\}$. We set
$f:=\left((\Psi_{B_2,B_1}\otimes\id_{B_1})\circ(\id_{B_2}\otimes\hat\rho)
\otimes\id_{B_2}\right)\circ(\id_{B_2\otimes B_1}\otimes\Delta_2)$.
These relations allow us to perform easliy the step from i=2 to i=0.
Finally, the case i=1 will be obtained from the corresponding identities
for i=0 by composition with $\varepsilon_2$.
\end{proof}

\begin{lemma}
For a strong Hopf datum $\mathfrak{h}$ the second module-coalgebra
compatibility, the first comodule-algebra compatibility and the
cycle-cocycle compatibility of Definition \ref{hp} can be respectively
converted to
\begin{gather}
\hat\rho\circ\varphi_{2,1}
 =\rho^{\mathrm{tr}}\circ\left(\id_{B_2}\otimes\ixi 11212011\right)
  \circ(\Delta_2\otimes\id_{B_1})\,,
\notag
\\
\mathlabel{phi-hat-new}
\varphi_{1,2}\circ\hat\sigma
 =(\id_{B_2}\otimes\m_1)\circ\left(\ixi 21222021\otimes\id_{B_1}\right)
  \circ\sigma^{\mathrm{tr}}\,,
\\
\hat\rho\circ\hat\sigma
 =(\id_{B_1}\otimes\m_1)\circ(\rho^{\mathrm{tr}}\otimes\id_{B_1})\circ
  \left(\id_{B_2}\otimes\ixi 11222001\otimes\id_{B_1}\right)\circ
  (\id_{B_2}\otimes\sigma^{\mathrm{tr}})\circ(\Delta_2\otimes\id_{B_2})\,.
\notag
\end{gather}
\end{lemma}

\begin{proof} 
The morphisms on right hand side of the corresponding
compatibility relations can be decomposed with the help of
\eqref{bra-ket-decomp}. Then the reduction formulas 
\eqref{tau-check-aux} will be applied to the particular tensor factors.
And finally we use again \eqref{bra-ket-decomp} to obtain the 
result.\end{proof}
\abs
The previous lemma leads to

\begin{lemma}
\begin{equation}\mathlabel{yi-red-adv}
\divide\unitlens by 2
\def\FillCircDiam{3}
\yi 01221 =
\begin{tangle}
\hh\step\id\step[3]\id\\
\hhstep\ffbox6{(\tau^3_t)^{\mathrm{red}_2}}\\
\vstr{110}\hh\id\step{\hstr{200}\cd*}\step\id\step\id\\
\hh\x\step{\hstr{200}\hld}\step\id\step\id\\
\hh\id\step\x\step\id\step\id\step\id\\
\hh\id\step\id\step{\hstr{200}\hlu}\step\id\step\id\\
\hh\id\step\id\Step\x\step\id\\
\hh\id\step{\hstr{200}\cu}\step\cu
\end{tangle}
\quad,\quad
\iy 11202 =
\begin{tangle}
\hh\cd\step{\hstr{200}\cd}\step\id\\
\hh\id\step\x\Step\id\step\id\\
\hh\id\step\id\step{\hstr{200}\hrd}\step\id\step\id\\
\hh\id\step\id\step\id\step\x\step\id\\
\hh\id\step\id\step{\hstr{200}\hru}\step\x\\
\vstr{110}\hh\id\step\id\step{\hstr{200}\cu*}\step\id\\
\hhstep\ffbox6{(\tau^3_b)^{\mathrm{red}_2}}\\
\hh\step\id\step[3]\id
\end{tangle}
\quad,\quad
\ixi 11222001 =
\begin{tangle}
\hh\step\id\step[4]\id\\
\hhstep\ffbox7{(\tau^3_t)^{\mathrm{red}_2}}\\
\vstr{110}\hh\cd\step{\hstr{200}\cd*}\step\id\step\id\\
\hh\id\step\x\Step\id\step\id\step\id\\
\id\step\id\step\ox{\sstyle 21}\step\id\step\id\\
\id\step\id\step\ox{\sstyle 12}\step\id\step\id\\
\hh\id\step\id\step\id\Step\x\step\id\\
\vstr{110}\hh\id\step\id\step{\hstr{200}\cu*}\step\cu\\
\hhstep\ffbox7{(\tau^3_b)^{\mathrm{red}_2}}\\
\hh\step\id\step[4]\id
\end{tangle}
\end{equation}
\begin{equation}\mathlabel{yi-red-adv2}
\divide\unitlens by 2
\yi 21201 =
\begin{tangle}
\hh\step[1.5]\id\step[1.5]\id\step[1.5]\id\\
\ffbox6{(\tau^3_t)^{\mathrm{red}_3}}\\
\hh\step\id\Step\id\step[1.5]\id\step\id\\
\ffbox4{\yi 21211}\hstep\id\step\id\\
\hh\hstep\id\step[3]\x\step\id\\
\hh\hstep{\hstr{300}\cu}\step\cu
\end{tangle}
\quad,\quad
\iy 10212 =
\begin{tangle}
\hh\cd\step{\hstr{300}\cd}\\
\hh\id\step\x\step[3]\id\\
\id\step\id\hstep\ffbox4{\iy 21212}\\
\hh\id\step\id\step[1.5]\id\Step\id\\
\hhstep\ffbox6{(\tau^3_b)^{\mathrm{red}_3}}\\
\hh\step\id\step[1.5]\id\step[1.5]\id
\end{tangle}
\end{equation}
\begin{equation}\mathlabel{yi-red-adv3}
\divide\unitlens by 2
\yi 01201 =
\begin{tangle}
\hh\step[1.5]\id\step[1.5]\id\step[1.5]\id\\
\ffbox6{(\tau^3_t)^{\mathrm{red}_1}}\\
\hh\hstep\id\step[1.5]\id\step[1.5]\id\step\id\step\id\\
\ffbox4{\yi 01201^*}\hstep\id\step\id\\
\hh\hstep\id\step[1.5]\id\step[1.5]\x\step\id\\
\hh\hstep\id\step[1.5]{\hstr{150}\cu}\step\cu
\end{tangle}
\quad,\quad
\iy 01202 =
\begin{tangle}
\hh\cd\step{\hstr{150}\cd}\step[1.5]\id\\
\hh\id\step\x\step[1.5]\id\step[1.5]\id\\
\id\step\id\hstep\ffbox4{\iy 01202^*}\\
\hh\id\step\id\step\id\step[1.5]\id\step[1.5]\id\\
\hhstep\ffbox6{(\tau^3_b)^{\mathrm{red}_1}}\\
\hh\step\id\step[1.5]\id\step[1.5]\id
\end{tangle}
\end{equation}
\end{lemma}

\begin{proof}
The first identity of \eqref{yi-red-adv} has been derived
by successive application of the left module-algebra
compatibility of Definition \ref{hp} and the relations
\eqref{mod-alg-rel} and \eqref{rel-cocycle-assoc1}.
The second identity is its $\pi$-symmetric counterpart.
To get the third identity we use \eqref{bra-ket-decomp},
the first and the second equation of \eqref{yi-red-adv},
the relations \eqref{rel-cocycle-assoc2} and the module-comodule 
compatibility.

The first identity in \eqref{yi-red-adv2} is obtained from the left
module-algebra compatibility and \eqref{mod-alg-rel}.
The identities in \eqref{yi-red-adv3} are derived from
\eqref{yi-red-adv2} with the help of \eqref{hat-prod-rel}.\end{proof}
\abs
The next lemma is a straightforward implication of the previous result.

\begin{lemma}
\begin{equation}\mathlabel{ixi-rel}
\divide\unitlens by 4
\multiply\unitlens by 3
\hstretch 120
\begin{tangle}
\hh\hstep\id\step\id\\
\ffbox2{\ixi 11222001}\\
\hh\hstep\id\step\id\\
\hh\hstep\id\hstep\hru
\end{tangle}
=
\begin{tangle}
\hh\hstep\id\step\id\\
\ffbox2{\ixi 11222221}\\
\hh\hstep\id\step\id\\
\hh\hstep\id\hstep\hru
\end{tangle}
\quad,\quad
\begin{tangle}
\hh\hstep\hld\hstep\id\\
\hh\hstep\id\step\id\\
\ffbox2{\ixi 11222001}\\
\hh\hstep\id\step\id
\end{tangle}
=
\begin{tangle}
\hh\hstep\hld\hstep\id\\
\hh\hstep\id\step\id\\
\ffbox2{\ixi 11222111}\\
\hh\hstep\id\step\id
\end{tangle}
\end{equation}
\end{lemma}
\begin{proof}
Both identities follow from the third identity of \eqref{yi-red-adv}
using \eqref{act-coact-triv}, \eqref{rel-cocycle-assoc1} and
\eqref{cocycle-triv1}.\end{proof}

\begin{lemma}\Label{sigma-Delta}
Given a strong Hopf datum $\mathfrak{h}$, then it holds
\begin{gather}
\mathlabel{sigma-Delta-rel}
\divide\unitlens by 2
\def\FillCircDiam{6}
\begin{tangle}
        \cu*\\
        \cd\\
        \hh\id\step[1.5]\hru
\end{tangle}
=
\multiply\unitlens by 2
\begin{tangle}
\vstr{60}\hh\hstep\id\step\id\\
\vstr{80}\ffbox2{\ixi 11222020}\\
\vstr{60}\hh\hstep\id\step\id\\
\hstr{50}\vstr{50}\hh\step\id\step[1.5]\hru\\
\end{tangle}
\quad,
\quad
\vstretch 100
\divide\unitlens by 2
\begin{tangle}
        \hh\hld\step[1.5]\id\\
        \cu\\
        \cd*
\end{tangle}
=
\multiply\unitlens by 2
\begin{tangle}
\hstr{50}\vstr{50}\hh\step\hld\step[1.5]\id\\
\vstr{60}\hh\hstep\id\step\id\\
\vstr{80}\ffbox2{\ixi 11220011}\\
\vstr{60}\hh\hstep\id\step\id\\
\end{tangle}
\\
\mathlabel{hat-sigma-Delta}
\hstretch 70
\begin{tangle}
\vstr{60}\hh\hstep\id\step\id\\
\hh\ffbox2{\hat\sigma}\\
\vstr{60}\hh\cd\hstep\id\\
\vstr{60}\hh\id\hstep\hru\hstep\id
\end{tangle}
=
\begin{tangle}
\vstr{60}\hh\step\id\Step\id\\
\vstr{150}\hh\ffbox4{\ixi 10222020}\\
\vstr{60}\hh\hstep\id\Step\id\step\id\\
\vstr{60}\hh\hstep\id\step[1.5]\hru\step\id
\end{tangle}
\quad,
\quad
\begin{tangle}
\vstr{60}\hh\hstep\id\hstep\hld\hstep\id\\
\vstr{60}\hh\hstep\id\hstep\cu\\
\hh\ffbox2{\hat\rho}\\
\vstr{60}\hh\hstep\id\step\id
\end{tangle}
=
\begin{tangle}
\vstr{60}\hh\hstep\id\step\hld\step[1.5]\id\\
\vstr{60}\hh\hstep\id\step\id\Step\id\\
\vstr{150}\hh\ffbox4{\ixi 01220011}\\
\vstr{60}\hh\step\id\Step\id
\end{tangle}
\end{gather}
\end{lemma}

\begin{proof}
We prove the first identity of \eqref{sigma-Delta-rel} graphically.
\begin{equation}\mathlabel{sigma-Delta-expl}
\divide\unitlens by 2
\def\FillCircDiam{6}
\begin{tangle}
        \cu*\\
        \cd\\
        \hh\id\step[1.5]\hru
\end{tangle}
=
\multiply\unitlens by 2
\vstretch 60
\begin{tangle}
\hh\hstep\id\step\id\\
\vstr{80}\ffbox2{\ixi 11222000}\\
\hh\hstep\id\step\id\\
\hstr{50}\vstr{50}\hh\step\id\step[1.5]\hru\\
\end{tangle}
=
\begin{tangle}
\hh\step\id\step\id\\
\vstr{40}\ffbox3{\sigma^{\mathrm{tr}}}\\
\hh\hstep\id\step\id\step\id\\
\vstr{80}\ffbox2{\ixi 11222001}\hstep\id\\
\vstr{40}\hh\hstep\id\step\id\step\id\\
\hh\vstr{50}{\hstr{50}\step\id\step[1.5]\hru\hhstep}\hstr{125}\ru\\
\vstr{40}\hh\hstep\id\step[0.75]\id
\end{tangle}
=
\begin{tangle}
\hh\step\id\step\id\\
\vstr{40}\ffbox3{\sigma^{\mathrm{tr}}}\\
\hh\hstep\id\step\id\step\id\\
\vstr{80}\ffbox2{\ixi 11222221}\hstep\id\\
\vstr{40}\hh\hstep\id\step\id\step\id\\
\hh\vstr{50}{\hstr{50}\step\id\step[1.5]\hru\hhstep}\hstr{125}\ru\\
\vstr{40}\hh\hstep\id\step[0.75]\id
\end{tangle}
=
\begin{tangle}
\hh\hstep\id\step\id\\
\vstr{80}\ffbox2{\ixi 11222020}\\
\hh\hstep\id\step\id\\
\hstr{50}\vstr{50}\hh\step\id\step[1.5]\hru\\
\end{tangle}
\end{equation}
The first equation in the graphic has been obtained with the help of the
cycle-cocycle compatibility (Definition \ref{hp}), 
\eqref{act-coact-triv} and \eqref{cocycle-triv1}. In the second
equation we use \eqref{bra-ket-decomp} and \eqref{tau-check-aux}.
The third equation is \eqref{ixi-rel}. Application of 
\eqref{tau-check-aux} to the fourth diagram yields the fifth relation.
The identities \eqref{hat-sigma-Delta} follow from
\eqref{sigma-Delta-rel} with the help of the entwining properties of 
$\varphi_{2,2}$ and $\varphi_{1,1}$ respectively.
\end{proof}

\begin{lemma}\Label{phi-sigma-rel}
Let $\mathfrak{h}$ be a strong Hopf datum. Then
\unitlens 15pt                                              %
\begin{equation}\mathlabel{phi-sigma}
\divide\unitlens by 2
\begin{tangle}
\hh\id\step\id\step\id\\
\hh\id\hstep\ffbox 2{\hat\sigma}\\
\ru\step\id\\
\ox{\sstyle 22}
\end{tangle}
=
\hstretch 120
\begin{tangle}
\hh\id\Step\cd\step\id\\
\ox{22}\hstep\ffbox2{\yi 21220}\\
\hh\id\Step\x\step\id\\
\cu\step\ru
\end{tangle}
\end{equation}
and the corresponding $\pi$-symmetric version for $\hat\rho$ holds.
\end{lemma}

\begin{proof}
The identity follows straighforwardly from \eqref{phi-hat-new}.
\end{proof}

\begin{lemma}
\begin{equation}\mathlabel{yi-red-adv4}
\unitlens = 8pt
\yi 00221^*=
\begin{tangle}
\hh\hstep\hstr{200}\cd\step\id\\
\hstep\id\Step\ox{\sstyle 22}\\
\ffbox3{\yi 01221^*}\step[1.5]\id\\
\hh\hstep\id\step\id\step\id\Step\id
\end{tangle}
\quad,\quad
\iy 11002^*=
\begin{tangle}
\hh\id\Step\id\step\id\step\id\\
\id\step[1.5]\ffbox3{\iy 11202^*}\\
\ox{\sstyle 11}\Step\id\\
\hh\hstr{200}\id\step\cu
\end{tangle}
\end{equation}
\begin{equation}\mathlabel{yi-red-adv5}
\unitlens = 8pt
\yi 00221=
\begin{tangle}
\hh\hstep\hstr{200}\cd\step\id\\
\hstep\id\Step\ox{\sstyle 22}\\
\ffbox3{\yi 01221}\step[1.5]\id\\
\hh\hstep\id\step\id\step\id\Step\id
\end{tangle}
\quad,\quad
\iy 11002=
\begin{tangle}
\hh\id\Step\id\step\id\step\id\\
\id\step[1.5]\ffbox3{\iy 11202}\\
\ox{\sstyle 11}\Step\id\\
\hh\hstr{200}\id\step\cu
\end{tangle}
\end{equation}
\end{lemma}

\begin{proof}
The identities \eqref{yi-red-adv4} follow from
\eqref{rel-cocycle-assoc1}, \eqref{mod-alg-rel}, whereas
\eqref{yi-red-adv5} will be derived from \eqref{yi-red-adv}
and the entwining properties \eqref{entw1}, \eqref{phi-prod2} of
$\varphi_{2,2}$, $\varphi_{2,1}$.
\end{proof}

\unitlens 12pt                                              %

\subsection{Proof of Theorem \ref{hp-cp}}

Given a strong Hopf datum $\mathfrak{h}$ we have to prove that the object
$B=B_1\otimes B_2$ provided with the structure given in Theorem \ref{hp-cp},
is a strong cross product bialgebra. The (co\n-)unital identities
\eqref{type-alpha1} can be verified easily. The ``strong conditions"
of Definition \ref{strong-type-alpha} hold by construction.
Using the first and the second identity of Lemma \ref{cocycle-triv2} it follows
straightforwardly that the ``projection relations" \eqref{type-alpha2} are
fulfilled. It remains to show that $B$ is a bialgebra.

\begin{proposition}\Label{co-assoc}
Let $\mathfrak{h}$ be a strong Hopf datum and $\varphi_{1,2}$, $\varphi_{2,1}$,
$\hat\rho$, and $\hat\sigma$ be the morphisms defined in \eqref{phi-sigma-rho2}.
Then $B=B_1\otimes B_2$ is a cocycle cross product algebra
$B_1\cpcyalg B_2$ and a cycle cross product coalgebra $B_1\cpcyco B_2$.
\end{proposition}

\begin{proof}
According to Proposition \ref{cross-equiv} the identities \eqref{cond3-caat}
have to be verified in order to prove that $B$ is a cross product
algebra. Similar $\pi$-symmetric procedures are needed to demonstrate that
$B$ is a cross product coalgebra.

Without difficulties the unital identities of \eqref{cond3-caat} can be
verified. The fifth relation of \eqref{cond3-caat} has been proven
in Lemma \ref{strong-hd7} in the first identity of \eqref{phi-prod2}.
The seventh identity of \eqref{cond3-caat} will be proven subsequently.
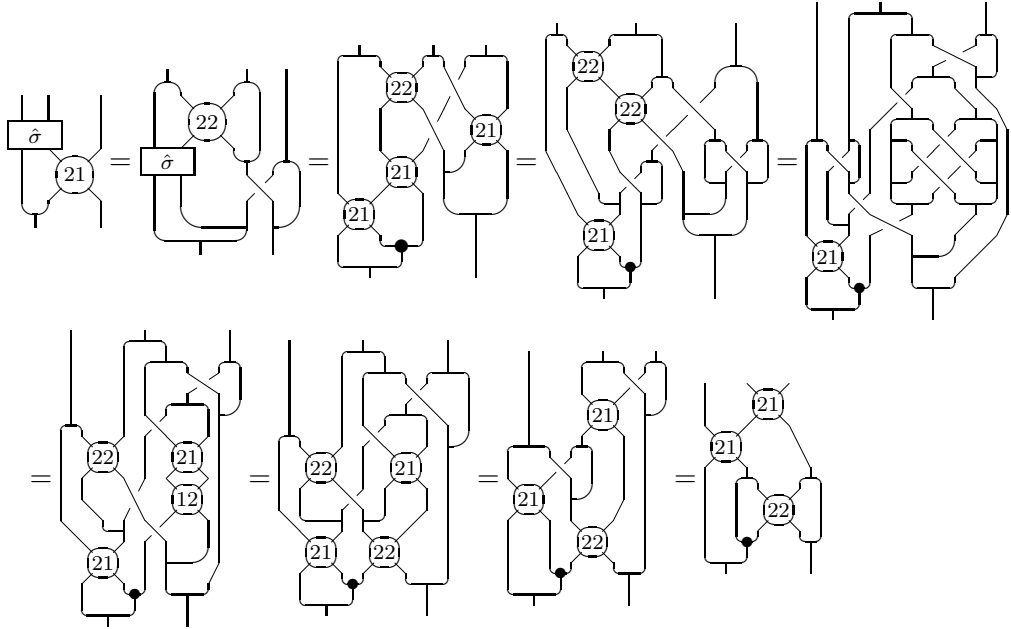
\begin{equation*}
\unitlens 10pt                                              %
\begin{split}
&\begin{tangle}
	\hh\id\step\id\Step\id\\
	\hh\hhstep\ffbox2{\hat\sigma}\step[1.5]\id\\
        \id\step\ox{21}\\
	\hh\cu\step[2]\id
\end{tangle}
=\enspace
\def\FillCircDiam{3}
\begin{tangle}
	\hh\cd\Step\cd\step\id\\
	\id\step\ox{22}\step\id\step\id\\
	\hh\hhstep\ffbox2{\hat\sigma}\step[1.5]\cu\hstep\cd\\
	\hh\id\step\id\step[2.5]\x\step\id\\
	\hh\id\step{\hstr{250}\lu}\step\ru\\
	\hh{\hstr{350}\cu}\step\id
\end{tangle}
=
\hstretch 80 \vstretch 80
\begin{tangle}
	\hh{\hstr{160}\cd}\Step\cd\step\hstr{160}\cd\\
	\id\Step\ox{22}\step\hx\Step\id\\
	\id\Step\id\Step\hx\step\ox{21}\\
	\id\Step\ox{21}\step\ru\Step\id\\
	\ox{21}\Step\id\step\hstr{120}\cu\\
        \hh\id\Step{\hstr{160}\cu*}\step[2.5]\id\\
	\hh{\hstr{240}\cu}\step[3.5]\id
\end{tangle}
=
\def\FillCircDiam{5}
\begin{tangle}
        \hh\cd\Step{\hstr{200}\cd}\step[3.5]\id\\
        \id\step\ox{22}\Step\hcd\Step\cd\\
        \id\step\id\Step\ox{22}\step\x\Step\id\\
        \hh\id\step\d\step[1.5]\id\Step\x\step[1.5]\cd\step\cd\\
        \hh\d\step\d\step\d\step\cd\hstep\d\step\id\step\x\step\id\\
        \hh\hstep\d\step\d\step\x\step\id\step\id\step
                                           \hstr{160}\hlu\hstep\hru\\
        \hh\step\d\step{\hstr{160}\hlu\hstep\hru}\step\Ru\step\id\\
        \step[1.5]\ox{21}\step\id\Step\hstr{120}\cu\\
        \hh\step[1.5]\id\Step\cu*\step[3.5]\id\\
        \hh\step[1.5]{\hstr{200}\cu}\step[4]\id
\end{tangle}
=
\begin{tangle}
\hh\hstep\id\step[1.5]{\hstr{240}\cd}\step[3.5]\id\\
\hh\hstep\id\step[1.5]\id\Step{\hstr{160}\cd}\Step\cd\\
\hh\hstep\id\step[1.5]\id\Step\id\Step{\hstr{160}\x}\step\id\\
\hh\hstep\id\step[1.5]\id\Step\id\step\hstr{160}\cd\hstep\hru\\
\hh\hstep\id\step[1.5]\id\Step\x\Step\nw1\nw1\\
\hh\hstep\id\step[1.5]\id\step\sw1{\hstr{160}\cd\hstep\cd}\d\\
\hh\cd\step\hrd\hstep\id\step{\hstr{160}\hrd}\step\x\step{\hstr{160}\hld}
        \hstep\id\\
\hh\id\step\x\hstep\id\hstep\id\step\id\step\x\step\x\step\id\hstep\id\\
\hh\id\step\id\step\hru\hstep\id\step{\hstr{160}\hru}\step\x\step
        {\hstr{160}\hlu}\hstep\id\\
\hh\id\step\id\step\x\step{\hstr{160}\cu\hstep\cu}\hstep\id\\
\hh\id\step{\hstr{160}\hlu\hstep\x\step\dd}\step\dd\\
\ox{21}\step\id\Step\Ru\sw2\\
\hh\id\Step\cu*\Step\hstr{160}\cu\\
\hh{\hstr{200}\cu}\step[3.5]\id
\end{tangle}
\\
&=
\def\FillCircDiam{5}
\vstretch 80 \hstretch 80
\begin{tangle}
        \hh\hstep\id\step[2.5]{\hstr{160}\cd}\step[3]\id\\
        \hh\hstep\id\step[2.5]\id\step{\hstr{160}\cd}\step[1.5]\cd\\
        \hh\hstep\id\step[2.5]\id\step\id\Step{\hstr{120}\x}\step\id\\
        \hh\hstep\id\step[2.5]\id\step\id\step{\hstr{160}\cd}\hstep\ru\\
        \hh\cd\Step\id\step\x\Step\id\hstep\id\\
        \id\step\ox{22}\step\id\step\ox{21}\hstep\id\\
        \id\step{\hstr{160}\dh}\step\hx\step\ox{12}\hstep\id\\
        \hh\nw1\step{\hstr{160}\hlu}\step\x\Step\id\hstep\id\\
        \step\ox{21}\step\id\step\Ru\ddh\\
        \hh\step\id\Step\cu*\step{\hstr{160}\cu}\\
        \hh\step{\hstr{200}\cu}\step[2.5]\id
\end{tangle}
=
\begin{tangle}
	\hh\hstep\id\step[2.5]{\hstr{160}\cd}\step[3]\id\\
	\hh\hstep\id\step[2.5]\id\step\hstr{160}\cd\hstep\cd\\
	\hh\hstep\id\step[2.5]\id\step\id\Step\x\Step\id\\
	\hh\hstep\id\step[2.5]\id\step\id\step{\hstr{160}\cd}
                                                   \nw1\step\id\\
	\hh\cd\Step\id\step\x\Step\id\step\ru\\
	\id\step\ox{22}\step\id\step\ox{21}\step\id\\
	\hh\id\step\id\Step\x\step\id\Step\id\step\id\\
	\hh\nw1{\hstr{320}\hlu}\step{\hstr{160}\hru}\step\sw1\step\id\\
	\step\ox{21}\step\ox{22}\Step\id\\
	\hh\step\id\Step\cu*\Step\hstr{160}\cu\\
	\hh\step{\hstr{200}\cu}\step[3.5]\id
\end{tangle}
=
\begin{tangle}
	\hh\step\id\step[2.5]{\hstr{160}\cd}\step\cd\\
	\hh\step\id\step[2.5]\id\Step\x\step\id\\
	\step\id\step[2.5]\ox{21}\step\ru\\
	\hh{\hstr{160}\cd}\step\cd\step[1.5]\id\step\id\\
	\hh\id\Step\x\step\id\step[1.5]\id\step\id\\
	\ox{21}\step\ru\step\ddh\step\id\\
	\id\Step\id\step\ox{22}\step[1.5]\id\\
	\hh\id\Step\cu*\Step\hstr{120}\cu\\
	\hh{\hstr{200}\cu}\step[3.25]\id
\end{tangle}
=
\begin{tangle}
	\id\Step\ox{21}\\
	\ox{21}\Step\d\\
	\hh\id\step[1.5]\cd\Step\cd\\
	\id\step[1.5]\id\step\ox{22}\step\id\\
	\hh\id\step[1.5]\cu*\Step\cu\\
	\hh{\hstr{160}\cu}\step[3]\id
\end{tangle}
\end{split}
\end{equation*}
\unitlens 12pt                               %
where the first identity has been obtained from \eqref{hat-prod}.
We use the weak associativity of $\mu_l$ and
the module-algebra compatibility of $\mu_r$ according to Definition \ref{hp}
to get the second equation in the graphic. Then the entwining property
of $\varphi_{2,2}$ (see Lemma \ref{strong-hd6}) yields the third relation.
In the fourth identity we use \eqref{rel-cocycle-assoc3} and the relative
coassociativity of $\Delta_1$ corresponding to Lemma \ref{strong-hd3}.
To derive the fifth identity the module-comodule compatibility of
Definition \ref{hp} has been used. In the sixth equation we apply
\eqref{phi-prod-rel2}, and in the seventh equation we use \eqref{phi-prod}.
Hence the seventh identity of \eqref{cond3-caat} has been verified.
Finally we will prove the sixth identity of \eqref{cond3-caat}. Its left
hand side will be transformed according to
\begin{equation}
\mathlabel{assoc-2l}
\divide\unitlens by 3
\multiply\unitlens by 2
\def\FillCircDiam{5}
\begin{split}
&\hstretch 80 \vstretch 80
\begin{tangle}
        \hh\id\Step\id\step[3]\id\\
        \hh\id\step[1.5]\ffbox4{\hat\sigma}\\
        \ox{\sstyle 21}\step[3]\id\\
        \hh\id\step[1.5]\cd\Step\cd\\
        \id\step[1.5]\id\step\ox{\sstyle 22}\step\id\\
        \hh\id\step[1.5]\cu*\Step\cu\\
        \hh{\hstr{160}\cu}\step[3]\id
\end{tangle}
= 
\begin{tangle}
        \hh\step\id\step[2.5]\cd\Step\cd\\
        \step\id\step[2.5]\id\step\ox{\sstyle 22}\step\id\\
        \hh\step\id\Step\ffbox2{\hat\sigma}\step[1.5]\cu\\
        \hh{\hstr{160}\cd}\step\cd\hstep\d\Step\id\\
        \hh\id\Step\x\step\id\step\id\Step\id\\
        \hh\id\Step\id\step{\hstr{160}\hru}\step\id\Step\id\\
        \ox{\sstyle 21}\step\ox{\sstyle 22}\Step\id\\
        \hh\id\Step\cu*\Step\hstr{160}\cu\\
        \hh{\hstr{200}\cu}\step[3.5]\id
\end{tangle}
= 
\begin{tangle}
\hh\step[1.5]\id\step[4]\cd\Step\cd\\
\step[1.5]\id\Step\sw2\Step\ox{\sstyle 22}\step\id\\
\hh\step[1.5]\id\step{\hstr{160}\cd\hstep\cd}\step\cu\\
\hh\step[1.5]\id\step{\hstr{160}\hrd}\step\id\step{\hstr{320}\hrd}\nw1
    \hstep\d\\
\hh\step[1.5]\id\step\id\step\id\step\x\Step\id\step\id\step\id\step\\
\step[1.5]\id\step\id\step\hx\step\ox{\sstyle 21}\step\id\step\id\\
\hh\step\cd\hstep\cu*\step\cu\step[1.5]\ffbox2{\hat\sigma}\hstep\id\\
\hh\sw1\step\x\Step\id\step[1.5]\sw1\sw1\step\id\\
\hh\id\Step\id\step{\hstr{320}\hru}\hstep\sw1\sw1\Step\id\\
\hh\id\Step\id\step{\hstr{400}\hru}\sw1\step[3]\id\\
\ox{\sstyle 21}\step{\hstr{100}\ox{\sstyle 22}}\Step\sw2\\
\Put(0,0)[lb]{\begin{tangle}\hh\id\Step\cu*\\ \hh{\hstr{200}\cu}
    \end{tangle}}\step[5.5]\cu
\end{tangle}
= 
\begin{tangle}
\hh\hstep\id\hstr{240}\step\cd\hstr{160}\step\cd\\
\hstep\id\step[3]\id\step[3]\ox{\sstyle 22}\Step\id\\
\hh\cd\Step\cd\step{\hstr{240}\cd}\hstep\hstr{160}\cu\\
\id\step\ox{\sstyle 22}\step\Put(0,0)[lb]{\begin{tangle}
    \hh\id\step{\hstr{320}\hrd}\step\nw1\hstep\nw1\\
    \hh\x\Step\id\Step\nw1\hstep\nw1
    \end{tangle}}\\
\id\step\id\Step\Put(0,0)[lb]{\begin{tangle}
    \hh\x\\ \hh\id\step\d\end{tangle}}
    \Step\ox{\sstyle 21}\step[3]\id\step[1.5]\hstr{160}\hd\\
\hh\id\step\nw1\step\id\step[1.5]\hru\step[1.5]\cd\step{\hstr{240}\cd}
    \step\id\\
\Put(0,0)[lb]{\begin{tangle}
    \hh\d\step[1.5]\cu*\\
    \hstep\ox{\sstyle 21}\\
    \hh\hstep\id\Step\nw1
    \end{tangle}}
    \step[4.5]
    \Put(0,0)[lb]{\begin{tangle}
    \ox{\sstyle 22}\\ \id\Step\hx
    \end{tangle}}
    \step[3]
    \Put(0,0)[lb]{\begin{tangle}
    \hh\id\step{\hstr{320}\hrd}\\
    \hh\x\Step\id\\
    \step\ox{\sstyle 21}
    \end{tangle}}
    \vstr{160}\step[4]\id\step\id\\
\hh\hstep\nw1\Step\id\step{\hstr{160}\cu\hstep\hru\step}\cu*\step\id\\
\hh\step[1.5]\nw1\step\def\FillCircDiam{2}{\hstr{160}\cu*\step}
    {\hstr{560}\hru}\step[1.5]\id\\
\hh\step[2.5]{\hstr{160}\cu}\step[3]\hstr{400}\cu
\end{tangle}
\\
&=  
\hstretch 80 \vstretch 80
\begin{tangle}
\hh\hstep\id\step[4.5]{\hstr{240}\cd}\Step{\hstr{160}\cd}\\
\hstep\id\step[3]\Put(0,0)[lb]{\begin{tangle}
    \hh\hstep\sw1\\ \hh\cd\end{tangle}}
    \step[4.5]\ox{\sstyle 22}\Step\id\\
\hh\hstep\id\step[3]\id\hstep\cd\Step{\hstr{160}\cd\hstep\cu}\\
\hcd\Step\ddh\hstep\id\step\ox{\sstyle 22}\Step
    \Put(0,0)[lb]{\begin{tangle}\hh\nw1\step\nw1\\
    \hh\step\d\step[1.5]\nw1
    \end{tangle}}\\
\id\step\ox{\sstyle 22}\step
    \Put(0,0)[lb]{\begin{tangle}
    \hh\id\step{\hstr{320}\hrd}\nw1\\
    \hh\x\Step\id\step\id
    \end{tangle}}
    \step[6.5]\id\step[2.5]\id\\
\id\step
    \Put(0,0)[lb]{\begin{tangle}
    \hh\id\Step\x\\
    \hh\def\FillCircDiam{2}{\hstr{160}\cu*}\step\d
    \end{tangle}}
    \step[4]\ox{\sstyle 21}\step\id\step[2.5]\id\step[2.5]\id\\
\hh\id\Step\id\step[2.5]\hru\Step\id\step\id\step{\hstr{240}\cd}\step\id\\
\ox{\sstyle 21}\step[2.5]{\hstr{100}\ox{\sstyle 22}}\step
    \Put(0,0)[lb]{\begin{tangle}
    \hh\id\step{\hstr{320}\hrd}\\
    \hh\x\Step\id
    \end{tangle}}
    \step[4]\id\step\id\\
{\hstr{160}\dh}\step\nw2\step[1.5]\dh\Step\hx\step\ox{\sstyle 21}%
 \step\id\step\id\\
\hh\step\nw1\Step\id\step{\hstr{160}\cu\hstep\hru\step}\cu*\step\id\\
\hh\Step\nw1\step\def\FillCircDiam{2}{\hstr{160}\cu*\step}
    {\hstr{560}\hru}\step[1.5]\id\\
\hh\step[3]{\hstr{160}\cu}\step[3]{\hstr{400}\cu}
\end{tangle}
=
\begin{tangle}
\hh\hstep\id\step[6]{\hstr{120}\cd}\Step\cd\\
\hstep\id\step[4]\sw2\step[2.5]\ox{\sstyle 22}\step\id\\
\hh\hstep\id\step[3.5]\cd\step[2.5]\cd\step[1.5]\cu\\
\hcd\Step\hddcd\step{\hstr{100}\ox{\sstyle 22}}\step\nw2\step\nw2\\
\id\step\ox{\sstyle 22}\step\Put(0,0)[lb]{
        \begin{tangle}
        \hh\id\step{\hstr{160}\hrd}\step[1.5]\nw1\Step\nw1\\
        \hh\x\step\nw1\step\cd\Step\cd
        \end{tangle}
    }\step[8]\hstep\id\\
\id\step\Put(0,0)[lb]{
        \begin{tangle}
        \hh\id\Step\x\\
        \hh\def\FillCircDiam{2}{\hstr{160}\cu*}\step\id
        \end{tangle}
    }\step[4]\ox{\sstyle 21}\step\id\step\ox{\sstyle 22}\step\id\hstep\id\\
\ox{\sstyle 21}\Step\id\step\ox{\sstyle 12}\step\Put(0,0)[lb]{
        \begin{tangle}
        \hh\id\step{\hstr{160}\hrd}\\
        \hh\x\step\id
        \end{tangle}
    }\step[3]\hcu*\hstep\id\\
\hh\id\Step\nw1\step\id\step\id\Step\x\step{\hstr{160}\hlu}\step[1.5]\id
    \step\id\\
\hh\id\step[3]\d\hstep\nw1{\hstr{160}\cu}\step\id\step\sw1\hstep\sw1\sw1\\
\hh\id\step[3.5]\d\step\x\step\sw1\sw1\hstep\sw1\sw1\\
\hh\id\step[4]\cu*\step{\hstr{160}\hru}\step[-1]\Ru\hstep\sw1\sw1\\
{\hstr{180}\cu}\step[1.5]\Put(0,0)[lb]{
        \begin{tangle}
        \hh{\hstr{400}\hru}\sw1\\
        \hh\hstr{200}\cu
        \end{tangle}
    }
\end{tangle}
\enspace=\enspace
\begin{tangle}
\hh\hstep\id\step[5]{\hstr{240}\cd}\step[2.5]{\hstr{160}\cd}\\
\hstep\id\step[4.5]\hdd\step[3]{\hstr{100}\ox{\sstyle 22}}\Step\nw2\\
\hh\hstep\id\step[3.5]{\hstr{160}\cd}\Step\cd\Step\nw1\step[3]\id\\
\hstep\id\step[3]\hdd\Step\ox{\sstyle 22}\step\nw2\Step\nw2\Step\id\\
\hh\cd\Step\cd\Step{\hstr{160}\hrd}\step\nw1\Step\nw1\Step\id\step\id\\
\id\step\ox{\sstyle 22}\step\x\step\nw2\step\nw2\Step\rd\hstep\ffbox2{\hat\sigma}\\
\hh\id\step\id\Step\x\step{\hstr{160}\cd\hstep\cd}\step\x\step\id
    \step\id\step\id\\
\id\step\cu*\step\id\step\Put(0,0)[lb]{
    \begin{tangle}
    \hh{\hstr{160}\hrd}\step\x\step{\hstr{160}\hld}\\
    \hh\id\step\x\step\x\step\id
    \end{tangle}
    }\step[6]\id\step\lu\step\id\step\id\\
\ox{\sstyle 21}\Step\id\step\Put(0,0)[lb]{
    \begin{tangle}
    \hh{\hstr{160}\hru}\step\x\step{\hstr{160}\hlu}\\
    \hh\hstr{160}\cu\hstep\cu
    \end{tangle}
    }\step[6]\id\Step\id\step\id\step\id\\
\hh\id\Step\nw1\step\nw1\step\nw1\Step{\hstr{160}\x}\Step\id\step\id
    \step\id\\
\hh\nw1\Step\nw1\step\nw1\step{\hstr{160}\cu}\Step\id\step\sw1\step\id
    \step\id\\
\hh\step\nw1\Step\nw1\step{\hstr{160}\x}\Step\sw1\sw1\step\sw1\step\id\\
\hh\Step\nw1\Step\cu*\Step{\hstr{320}\hru}\step[-2]{\hstr{240}\ru}\step
    \sw1\step\sw1\\
\hh\step[3]\nw1\step\dd\step[2.5]{\hstr{640}\hru}\step\sw1\\
\hh\step[4]\cu\step[3]\hstr{400}\cu
\end{tangle}
\end{split}
\end{equation}
where we used \eqref{phi-prod} and \eqref{hat-prod} in the first equation,
and \eqref{hat-sigma-Delta} and the right modularity of $(B_2,\mu_r)$ in
the second equation. To get the third identity in the graphic one has to
apply \eqref{phi-sigma} and again the fact that $(B_2,\mu_r)$ is right module.
In the fourth relation we use \eqref{phi-prod} and the relative comodule
property of $\nu_r$ according to Lemma \ref{strong-hd3}. The fifth 
identity can be
verified with the help of \eqref{rel-phi}. The module-comodule compatibility,
the associativity of $\Delta_2$ and the entwining property \eqref{entw1} of
$\varphi_{2,2}$ yield the final equation in the graphic.

For the right hand side of the sixth identity of \eqref{cond3-caat} we get
\begin{equation}\mathlabel{assoc-2r}
\divide\unitlens by 4
\multiply\unitlens by 3
\begin{split}
&\hstretch 80 \vstretch 80
\begin{tangle}
        \hh\id\step\id\step[3]\id\\
        \hh\hhstep\ffbox2{\hat\sigma}\step[2.5]\id\\
        \hh\id\hstep\cd\Step\cd\\
        \id\hstep\id\step\ox{22}\step\id\\
        \hh\id\hstep\cu*\Step\cu\\
        \hh\cu\step[3]\id
\end{tangle}
=
\begin{tangle}
        \hh\cd\Step\cd\step\id\\
        \id\step\ox{22}\step\id\hstep\hcd\\
        \hh\id\step\id\Step\cu\hstep\hrd\hstep\id\\
        \hh\hhstep\ffbox2{\hat\sigma}\Step\x\hstep\id\hstep\id\\
        \hh\def\FillCircDiam{2}
        \id\step{\hstr{200}\cu*}\step\hru\hstep\id\\
        \hh{\hstr{180}\cu}\step[2.25]\cu
\end{tangle}
=
\begin{tangle}
        \hh\hstep\id\step[3]\id\step[1.5]\hstr{200}\cd\\
        \hh\cd\Step\cd\step{\hstr{320}\hrd}\hstep\id\\
        \id\step\ox{22}\step\hx\Step\id\hstep\id\\
        \id\step\id\Step\hx\step\ox{21}\hstep\id\\
        \hh\hhstep\ffbox2{\hat\sigma}\hstep\sw1\step{\hstr{160}\hru}
                \Step\id\hstep\id\\
        \hh\id\step\cu*\Step{\hstr{240}\cu}\hstep\id\\
        \hh{\hstr{120}\cu}\step[4]{\hstr{160}\cu}
\end{tangle}
=
\begin{tangle}
        \hh\hstep\id\step[3]\id\step[1.5]\hstr{240}\cd\\
        \hh\cd\Step\cd\step{\hstr{320}\hrd}\step\id\\
        \id\step\ox{22}\step\hx\Step\id\step\id\\
        \id\step\d\step\hx\step\ox{21}\step\id\\
        \hh\id\step[1.5]\ffbox2{\hat\sigma}\hstep{\hstr{160}\hru}
                \step[1.5]\ffbox2{\hat\sigma}\\
        \ox{21}\step\id\step{\hstr{120}\Ru}\step\id\\
        \hh\id\Step\cu*\step{\hstr{320}\cu}\\
        \hh{\hstr{200}\cu}\step[3.5]\id
\end{tangle}
=
\begin{tangle}
\hh\step\id\step[4]\id\step[3]\hstr{200}\cd\\
\hh\step\id\step[4]\id\step[1.5]{\hstr{240}\cd}\step\id\\
\hh{\hstr{160}\cd}\step[2.5]\cd\step{\hstr{160}\hrd}\Step\id\step\id\\
\id\Step{\hstr{100}\ox{22}}\step\hx\step\ox{12}\step\id\\
\hh\id\step[1.5]\cd\Step\x\step\x\Step\id\step\id\\
\id\step[1.5]\id\step\ox{22}\step\hx\step\ox{21}\step\id\\
\hh\id\step[1.5]\cu*\Step\cu\step\d\hstep\id\step[1.5]\ffbox2{\hat\sigma}\\
\ox{21}\step\sw2\step[3]\hd\cu\step\id\\
\hh\id\Step\cu*\step[4.5]{\hstr{160}\hru}\Step\id\\
\hh{\hstr{200}\cu}\step[5]\hstr{240}\cu
\end{tangle}
\\
&=\enspace
\hstretch 80 \vstretch 80
\begin{tangle}
\hh\step[1.25]\id\step[3.75]\id\step[4.5]\hstr{240}\cd\\
\hh\step[1.25]\id\step[3.75]\id\step[2.5]{\hstr{320}\cd}\step\id\\
\hh{\hstr{200}\cd}\Step\cd\Step{\hstr{400}\hrd}\step[1.5]\id\step\id\\
\id\step[2.5]\ox{22}\step\x\step[2.5]\id\step[1.5]\id\step\id\\
\hh\id\Step\dd\Step\x\step{\hstr{160}\cd}\step\cd\step\id\step\id\\
\hh\id\Step\id\step[2.5]\id\step\id\step\id\Step\x\step\id\step\id\step\id\\
\hh\id\step[1.5]\cd\Step\id\step\id\step\id\step\sw1\step{\hstr{160}\hru}
        \step\id\step\id\\
\id\step[1.5]\id\step\ox{22}\step\id\step\id\step\id\Step\ox{22}\step\id\\
\id\step[1.5]\id\step\id\Step\id\step\id\step\id\step\ox{12}\Step\id
        \step\id\\
\hh\id\step[1.5]\cu*\Step\id\step\id\step\x\Step\id\step[1.5]
        \ffbox2{\hat\sigma}\\
\ox{21}\step[2.5]\id\step\hx\step\Lu\Step\id\step\id\\
\hh\id\Step\id\step[2.5]\cu\step\id\step[3]{\hstr{160}\cu}\step\id\\
\hh\def\FillCircDiam{2}\id\Step{\hstr{240}\cu*}\step[1.5]{\hstr{640}\hru}
        \Step\id\\
\hh{\hstr{280}\cu}\step[3]\hstr{480}\cu
\end{tangle}
=
\begin{tangle}
\hh\step\id\step[4]{\hstr{200}\cd}\Step\cd\\
\step\id\step[4]\id\step[2.5]\ox{22}\step{\hstr{160}\hd}\\
\hh\step\id\step[4]\id\step[1.5]{\hstr{160}\cd}\step\nw1\step\id\\
\cd\step[2.5]\hcd\step\Put(0,0)[lb]{
        \begin{tangle}
        \hh{\hstr{240}\hrd}\hstep\d\step\ffbox2{\hat\sigma}\\
        \hh\id\hstep{\hstr{160}\hld\hstep\hrd}\hstep\id\step\id
        \end{tangle}
    }\\
\id\Step{\hstr{100}\ox{22}}\step\Put(0,0)[lb]{
        \begin{tangle}
        \hh\id\step\id\hstep\id\step\x\step\id\hstep\id\step\id\\
        \hh\x\hstep\cu\step\d\hstep\id\hstep\id\step\id
        \end{tangle}
    }\\
\hh\id\step[1.5]\cd\Step\x\step\x\Step\id\hstep\id\hstep\id
        \step\id\\
\id\step[1.5]\id\step\ox{22}\step\hx\step\ox{21}\hstep\id\hstep\id
        \step\id\\
\hh\id\step[1.5]\cu*\Step\cu\step\id\hstep\dd\Step\hlu\hstep\id\step\id\\
\ox{21}\step[3]\id\step[1.5]\Put(0,0)[lb]{
        \begin{tangle}
        \hh\id\hstep{\hstr{240}\cu}\hstep\id\step\id\\
        \hh{\hstr{320}\hru}\step[-2]{\hstr{160}\Ru}\step\id
        \end{tangle}
    }\\
\Put(0,0)[lb]{
        \begin{tangle}
        \hh\def\FillCircDiam{2}\id\Step{\hstr{240}\cu*}\\
        \hh{\hstr{280}\cu}
        \end{tangle}
    }\step[6.5]\hstr{200}\cu
\end{tangle}
\end{split}
\end{equation}
In the first equation of \eqref{assoc-2r} we used \eqref{hat-prod}.
Then we applied \eqref{entw1} to get the second equation.
To derive the third identity we make use of the cocycle compatibility 
and the weak associativity of $\m_2$ according to Definition \ref{hp}.
In the forth equation we apply the right comodule-coalgebra compatibility,
and the fifth identity follows from \eqref{phi-phi}.
The sixth equation in the graphical calculation \eqref{assoc-2r}
has been obtained from the entwining property of $\varphi_{2,2}$
and the left module-algebra compatibility.

Eventually, the seventh graphic in \eqref{assoc-2l} and in 
\eqref{assoc-2r} coincide. This can be verified by applying
Lemma \ref{strong-hd3} for $\m_2$ and $\nu_r$, coassociativity of 
$\Delta_2$, the right modularity of $(B_2,\mu_r)$ and
\eqref{mod-alg-rel} of Lemma \ref{strong-hd4}.\end{proof}
\abs
\begin{remark}\Label{dress}
{\normalfont To finish the proof of Theorem \ref{hp-cp} we have to show that
$\Delta_B$ is an algebra morphism. For this purpose we will use a sort of
\emph{dressing transformation} technique. This means we define two sequences
of morphisms $(\beta_\ell^i)_{i=0}^4$ and $(\beta_r^i)_{i=0}^4$
and a sequence of ``dressing transformations" $\left(T^{(i)}\right)_{i=1}^4$ 
where $\beta_r^0=(\m_B\otimes\m_B)\circ(\id\otimes\Psi_{B,B}\otimes\id)
\circ(\m_B\otimes\m_B)$, $\beta_\ell^0=\Delta_B\circ\m_B$, 
$\beta_\ell^4=\beta_r^4$, and $T^{(i)}(\beta_\ell^i)=\beta_\ell^{i-1}$,
$T^{(i)}(\beta_r^i)=\beta_r^{i-1}$ for all $i\in\{1,2,3,4\}$.
This implies $\beta_\ell^0=\beta_r^0$ and therefore the statement.}
\end{remark}

\begin{remark}{\normalfont Once it is proven that $B$ is an algebra,
Proposition \ref{alg-morph-equiv} provides an alternative method to show that
$\Delta_B$ is an algebra morphism. The present proof using dressing
transformations is invariant under $\pi$-symmetry, however.}
\end{remark}
\abs
The ``dressing transformations'' will be defined as
\begin{gather}\mathlabel{dressing}
\divide\unitlens by 3
\multiply\unitlens by 2
\hstretch 80
\vstretch 120
T^{(1)}(f):=
\begin{tangle}
        \hh\id\step[1.5]\id\step\id\step\id\\
        \hh\vstr{144}\id\hstep\ffbox4{\tau^1_t}\\
        \hh\id\step\id\step\id\step\id\step\id\\
        \hh\vstr{144}\id\hstep\ffbox3{f}\hstep\id\\
        \hh\id\step\id\step\id\step\id\step\id\\
        \hh\vstr{144}\hhstep\ffbox4{\tau^1_b}\hstep\id\\
        \hh\hstep\id\step\id\step\id\step[1.5]\id
\end{tangle}
\,\,,\,\,
T^{(2)}(f):=
\begin{tangle}
        \hh{\hstr{120}\cd}\step\id\step\id\\
        \hh\vstr{144}\id\hstep\ffbox4{\tau^2_t}\\
        \hh\id\step\id\step\id\step\id\step\id\\
        \hh\vstr{144}\id\hstep\ffbox3{f}\hstep\id\\
        \hh\id\step\id\step\id\step\id\step\id\\
        \hh\vstr{144}\hhstep\ffbox4{\tau^2_b}\hstep\id\\
        \hh\hstep\id\step\id\step\hstr{120}\cu
\end{tangle}
\,\, ,\,\,
\hstretch 70
T^{(3)}(g):=
\begin{tangle}
        \hh\cd\step[1.5]{\hstr{140}\cd}\step[2]\id\\
        \hh\id\step{\hstr{105}\x}\step[2]\id\step[2]\id\\
        \hh\vstr{144}\id\step\id\hstep\ffbox6{\tau^{3}_t}\\
        \hh\id\step\id\step\id\step\id\step\id\step\id\step\id\step\id\\
        \hh\vstr{144}\id\step\id\hstep\ffbox4{g}\hstep\id\step\id\\
        \hh\id\step\id\step\id\step\id\step\id\step\id\step\id\step\id\\
        \hh\vstr{144}\hhstep\ffbox6{\tau^{3}_b}\hstep\id\step\id\\
        \hh\hstep\id\step[2]\id\step[2]{\hstr{105}\x}\step\id\\
        \hh\hstep\id\step[2]{\hstr{140}\cu}\step[1.5]\cu
\end{tangle}
\,\, ,\,\,
T^{(4)}(h):=
\begin{tangle}
\hh\id\step[1.5]{\hstr{210}\hld}\step\id\step\id\\
\hh{\hstr{105}\x}\step[1.5]\id\step\id\step\id\\
\hh\id\step\cd\step\id\step\id\step\hstr{140}\hrd\\
\hh\id\hstep\ffbox5{h}\hstep\id\\
\hh{\hstr{140}\hlu}\step\id\step\id\step\cu\step\id\\
\hh\step\id\step\id\step\id\step[1.5]\hstr{105}\x\\
\hh\step\id\step\id\step{\hstr{210}\hru}\step[1.5]\id\\
\end{tangle}
\end{gather}
for any $f\in\Hom(B_{212},B_{121})$, $g\in\Hom(B_{2122},B_{1121})$
and $h\in\Hom(B_{22122},B_{11211})$.
\abs
The remainder of this section is devoted to the proof of the 
following proposition.

\begin{proposition}\Label{b-bialg}
Given a strong Hopf datum $\mathfrak{h}$,
then there exist two sequence of morphisms $(\beta^i_\ell)_{i=0}^4$,
$(\beta^i_r)_{i=0}^4$ with $\beta^0_\ell:=\Delta_B\circ\m_B$ and
$\beta^0_r:=(\m_B\otimes\m_B)\circ(\id_{B}\otimes\Psi_{B,B}\otimes\id_{B})\circ
(\Delta_{B}\otimes\Delta_{B})$ such that
\begin{equation}\mathlabel{beta-lr}
\begin{diagram}[height=15pt]
\beta^0_\ell&\lMapsto{T^{(1)}}&
\beta^1_\ell&\lMapsto{T^{(2)}}&
\beta^2_\ell&\lMapsto{T^{(3)}}&
\beta^3_\ell&\lMapsto{T^{(4)}}&
\beta^4_\ell\\
 & & & & & & & &\dEqualsto\\
\beta^0_r&\lMapsto{T^{(1)}}&
\beta^1_r&\lMapsto{T^{(2)}}&
\beta^2_r&\lMapsto{T^{(3)}}&
\beta^3_r&\lMapsto{T^{(4)}}&
\beta^4_r
\end{diagram}
\end{equation}
Therefore $\beta^0_\ell=\beta^0_r$ and $B=B_1\otimes B_2$ is a
bialgebra which proves Theorem \ref{hp-cp}.
\end{proposition}

\begin{proof}
The proof will be splitted into several parts. At first we verify the following
diagram.
\begin{equation}
\begin{diagram}[height=18pt]
\beta^0_r&&\beta^1_r&&\beta^2_r&&\beta^3_r\\
\dCongruent&&\dCongruent&&\dCongruent&&\dCongruent\\
\ixi 00000000&&\ixi 01200000&&\ixi 01202001&&\\
\dEqualsto&&\dEqualsto&&\dEqualsto&&\\
\hstretch 60
\begin{tangle}
        \hh\id\step\cd\step\id\step\id\\
        \hh\vstr{120}\id\step\id\hstep\ffbox3{\yi 00200}\\
        \hh\id\step\id\step[.75]\id\hstep\id\step[.75]\id\step\id\\
        \hh\vstr{120}\hhstep\ffbox3{\iy 01000}\hstep\id\step\id\\
        \hh\id\step\id\step\cu\step\id
\end{tangle}
&\lMapsto{T^{(1)}}&
\hstretch 60
\begin{tangle}
        \hh\cd\step\id\step\id\\
        \hh\vstr{120}\id\hstep\ffbox3{\yi 01200}\\
        \hh\id\step\id\step\id\step\id\\
        \hh\vstr{120}\hhstep\ffbox3{\iy 01200}\hstep\id\\
        \hh\id\step\id\step\cu
\end{tangle}
&\lMapsto{T^{(2)}}&
\hstretch 60
\begin{tangle}
        \hh\cd\step\id\step\id\\
        \hh\vstr{120}\id\hstep\ffbox3{\yi 01201}\\
        \hh\id\step\id\step\id\step\id\\
        \hh\vstr{120}\hhstep\ffbox3{\iy 01202}\hstep\id\\
        \hh\id\step\id\step\cu
\end{tangle}
&\lMapsto{T^{(3)}}&
\hstretch 60 \vstretch 60
\begin{tangle}
\hh{\hstr{150}\cd}\step\cd\Step\cd*\hstep\id\\
{\hstr{90}\rd}\step\hx\step\ox{11}\step\id\hstep\id\\
\hh\id\step[1.5]\x\step\x\Step\cu\hstep\id\\
\hh\id\step\dd\step\x\step\id\step{\hstr{90}\ld}\step\id\\
\hh\id\step\id\step\cd\hstep\id\step\id\hstep\cd\step\x\\
\hh\id\step\x\step\id\hstep\id\step\id\hstep\id\step\x\step\id\\
\hh\x\step\cu\hstep\id\step\id\hstep\cu\step\id\step\id\\
\hh\id\step{\hstr{90}\ru}\step\id\step\x\step\dd\step\id\\
\hh\id\hstep\cd\Step\x\step\x\step[1.5]\id\\
\id\hstep\id\step\ox{22}\step\hx\step{\hstr{90}\lu}\\
\hh\id\hstep\cu*\Step\cu\step\hstr{150}\cu
\end{tangle}
\end{diagram}
\end{equation}
where the identities ($\equiv$) in the first row are definitions and
the equalities ($=$) in the second row are special cases
of \eqref{bra-ket-decomp}. The morphisms in the third row will be obtained by
applying $T^{(1)}$, $T^{(2)}$ and $T^{(3)}$ using the identities
\eqref{tau-R-aux}, \eqref{tau-R2-aux} and \eqref{yi-red-adv3} respectively.

A similar diagram can be set up for the $(\beta^i_\ell)$. 
For that we use the following auxilarity definitions
\begin{equation}\mathlabel{aux-gamma0}
\gamma^0_t:=(\id_{B_1}\otimes\Delta_1)\circ\m_{0,20}
\quad,\quad
\gamma^0_b:=\Delta_{01,0}\circ(\m_1\otimes\id_{B_2})
\end{equation}
\begin{equation}
\mathlabel{aux-gamma}
\begin{gathered}
\unitlens = 8pt
\gamma^1_t:=
\begin{tangle}
\hh\step\id\step\cd\step[2.5]\id\\
\hh\step\x\step\id\step[2.5]\id\\
\hh\sw1\step{\hstr{200}\hru}\step[2.5]\id\\
\hh\id\step{\hstr{200}\cd\step}\cd\\
\id\step\id\Step\ox{\sstyle 22}\step\id\\
\id\hstep\ffbox3{\yi 01220}\step[1.5]\hcu\\
\id\step\id\step\id\step\hstr{125}\ox{\sstyle 12}
\end{tangle}
\quad,\quad
\gamma^2_t:=
\begin{tangle}
\hh\step\id\step\cd\step[3]\id\\
\hh\step\x\step\id\step[3]\id\\
\hh\sw1\step\ru\step[3]\id\\
\hh\id\step\hstr{200}\cd\step\cd\\
\id\step\id\Step\ox{\sstyle 22}\Step\id\\
\vstr{115}\id\hstep\ffbox3{\yi 01221}\step\ffbox3{\ixi 21222021}\\
\hh\id\step\id\hstep\id\hstep\hstr{200}\hld\step\id\step\id\\
\hh\id\step\id\hstep\id\hstep\id\step\hstr{200}\x\hstep\dd\\
\hh\id\step\id\hstep\id\hstep\cu\Step\cu
\end{tangle}
\quad,\quad
\unitlens = 9pt
\def\FillCircDiam{2}
\gamma^3_t:=
\begin{tangle}
\hh\id\step\cd\hstep\hstr{200}\cd*\hstep\id\\
\hh\x\step\id\hstep\id\step{\hstr{200}\hld}\step\rd\\
\hh\id\step{\hstr{200}\hru}\hstep\id\step\id\step\x\step\id\\
\hh\id\step{\hstr{150}\x}\step\cu\step\id\step\id\\
\hh\id\step\id{\hstr{150}\step\x\step}\id\step\id\\
\id\step\id\step[1.5]\id\step\hstr{50}\ffbox{5}{\ixi 21212001}\step\id\\
\hh\id\step\id\step[1.5]\id\step[1.5]\id\step[1.5]\x\\
\hh\id\step\id\step[1.5]\id\step[1.5]{\hstr{300}\hru}\step\id
\end{tangle}
\\
\unitlens = 8pt
\gamma^1_b:=
\begin{tangle}
\hstep{\hstr{125}\ox{\sstyle 21}}\step\id\step\id\step\id\\
\hcd\step[1.5]\ffbox3{\iy 11200}\hstep\id\\
\id\step\ox{\sstyle 11}\Step\id\step\id\\
\hh\cu\Step{\hstr{200}\cu}\step\id\\
\hh\hstep\id\step[2.5]\hstr{200}\hld\hstep\dd\\
\hh\hstep\id\step[2.5]\id\step\x\\
\hh\hstep\id\step[2.5]\cu\step\id
\end{tangle}
\quad,\quad
\gamma^2_b:=
\begin{tangle}
\hh\step\cd\Step\cd\hstep\id\hstep\id\step\id\\
\hh{\hstr{200}\dd\hstep\x}\step\id\hstep\id\hstep\id\step\id\\
\hh{\hstr{200}\id\step\id\step\hru}\hstep\id\hstep\id\step\id\\
\vstr{115}\hhstep\ffbox3{\ixi 11212011}\step\ffbox3{\iy 11202}\hstep\id\\
\id\Step\ox{\sstyle 11}\Step\id\step\id\\
\hh\hstr{200}\cu\step\cu\hstep\id\\
\hh\step\id\step[3]\ld\step\sw1\\
\hh\step\id\step[3]\id\step\x\\
\hh\step\id\step[3]\cu\step\id
\end{tangle}
\quad,\quad
\unitlens = 9pt
\def\FillCircDiam{2}
\gamma^3_b:=
\begin{tangle}
\hh\id\step{\hstr{300}\hld\hstep\id\hstep}\id\step\id\\
\hh\x{\hstr{150}\step\id\step\id\step\id}\step\id\\
\id\hstep{\hstr{50}\ffbox{5}{\ixi 21212001}}\step\id\step[1.5]\id\step\id\\
\hh\id\step\id{\hstr{150}\step\x\step}\id\step\id\\
\hh\id\step\id\step\cd\step{\hstr{150}\x}\step\id\\
\hh\id\step\x\step\id\step\id\hstep{\hstr{200}\hld}\step\id\\
\hh\lu\step{\hstr{200}\hru}\step\id\hstep\id\step\x\\
\hh\step\id\step{\hstr{200}\cu*}\hstep\cu\step\id
\end{tangle}
\end{gathered}
\end{equation}
Then
\begin{equation}
\begin{diagram}[height=18pt]
\beta^0_\ell && \beta^1_\ell && \beta^2_\ell && \beta^3_\ell\\
\dCongruent && \dCongruent && \dCongruent && \dCongruent\\
\Delta_0\circ\m_0 && \Delta_{01,0}\circ\m_{0,20} && &&\\
\dEqualsto&\ldMapsto^{\sstyle T^{(1)}}&\dEqualsto&& &&\\
\hstretch 60
\begin{tangle}
        \hh\id\step\cd\step\id\step\id\\
        \hh\vstr{120}\id\step\id\hstep\ffbox3{\gamma^0_t}\\
        \hh\id\step\id\step[.75]\id\hstep\id\step[.75]\id\step\id\\
        \hh\vstr{120}\hhstep\ffbox3{\gamma^0_b}\hstep\id\step\id\\
        \hh\id\step\id\step\cu\step\id
\end{tangle}
&&
\hstretch 60
\begin{tangle}
        \hh\cd\step\id\step\id\\
        \hh\vstr{120}\id\hstep\ffbox3{\gamma^1_t}\\
        \hh\id\step\id\step\id\step\id\\
        \hh\vstr{120}\hhstep\ffbox3{\gamma^1_b}\hstep\id\\
        \hh\id\step\id\step\cu
\end{tangle}
&\lMapsto{T^{(2)}}&
\hstretch 60
\begin{tangle}
        \hh\cd\step\id\step\id\\
        \hh\vstr{120}\id\hstep\ffbox3{\gamma^2_t}\\
        \hh\id\step\id\step\id\step\id\\
        \hh\vstr{120}\hhstep\ffbox3{\gamma^2_b}\hstep\id\\
        \hh\id\step\id\step\cu
\end{tangle}
&\lMapsto{T^{(3)}}&
\hstretch 60
\begin{tangle}
        \hh\cd\step\id\step\id\\
        \hh\vstr{120}\id\hstep\ffbox3{\gamma^3_t}\\
        \hh\id\step\id\step[0.33]\id\step[0.34]\id\step[0.33]\id\step\id\\
        \hh\vstr{120}\hhstep\ffbox3{\gamma^3_b}\hstep\id\\
        \hh\id\step\id\step\cu
\end{tangle}
\end{diagram}
\end{equation}
Again the identities ($\equiv$) in the first row are definitions.
The remaining identities will be proven in the subsequent lemmas. 

\begin{lemma}
\begin{equation}\mathlabel{tau-L-aux}
\gamma^0_t=(\m_{0,20}\otimes\id_{B_2})\circ\tau^1_t\,,\qquad
\gamma^0_b=\tau^1_b\circ(\id_{B_1}\otimes\Delta_{01,0})
\end{equation}
and hence $\beta^0_\ell=T^{(1)}\left(\beta^1_\ell\right)$.
\end{lemma}
\begin{proof}
The identities are obtained using \eqref{hat-prod} and \eqref{phi-prod}.
\end{proof}

\begin{lemma}
\begin{gather}\mathlabel{beta-1l}
\beta^1_\ell=(\id\otimes\m_1)\circ(\gamma^1_b\otimes\id_{B_1})\circ
 (\id_{B_2}\otimes\gamma^1_t)\circ(\Delta_2\otimes\id)
\\
\mathlabel{gamma-12}
\gamma^1_t=(\id\otimes\m_1)\circ(\gamma^2_t\otimes\id_{B_1})\circ\tau^2_t\,,
\quad
\gamma^1_b=\tau^2_b\circ(\id_{B_2}\otimes\gamma^2_b)\circ(\Delta_2\otimes\id)
\end{gather}
and therefore $\beta^1_\ell=T^{(2)}\left(\beta^2_\ell\right)$.
\end{lemma}
\begin{proof}
The proof of the identity  \eqref{beta-1l} will be given
in the following graphical calculation.
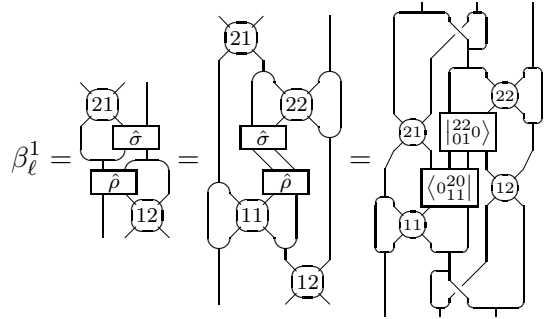
\begin{equation*}
\beta^1_\ell=
\vstretch 70 \hstretch 70
\begin{tangle}
	\ox{21}\step\id\\
	\hh\id\step[1.5]\ffbox2{\hat\sigma}\\
	\hh\hstr{140}\cu\cd\\
	\hh\hstep\ffbox2{\hat\rho}\step[1.5]\id\\
	\step\id\step\ox{12}\\
\end{tangle}
=\begin{tangle}
	\hstep\ox{21}\step[3]\id\\
	\hh\hstep\id\step[1.5]\cd\Step\cd\\
	\hstep\id\step[1.5]\id\step\ox{22}\step\id\\
	\hh\hstep\id\step\ffbox2{\hat\sigma}\step[1.5]\cu\\
	\hh\hstep\id\step[1.5]\nw1\nw1\step[1.5]\id\\
	\hh\cd\step[1.5]\ffbox2{\hat\rho}\step\id\\
	\id\step\ox{11}\step\id\step[1.5]\id\\
	\hh\cu\Step\cu\step[1.5]\id\\
	\hstep\id\step[3]\ox{12}
\end{tangle}
=\enspace
\unitlens=7pt
\vstretch 100 \hstretch 100
\begin{tangle}
\hh{\hstr{300}\cd}\step\cd\step[2.5]\id\\
\hh\id\step[3]\x\step\id\step[2.5]\id\\
\hh\id\Step\sw1\step{\hstr{200}\hru}\step[2.5]\id\\
\hh\id\Step\id\step{\hstr{200}\cd\step}\cd\\
\id\Step\id\step\id\Step\ox{\sstyle 22}\step\id\\
\ox{\sstyle 21}\hstep\ffbox3{\yi 01220}\step[1.5]\hcu\\
\hh\hhstep\dd\Step\id\step\id\step\id\step\id\Step\dd\\
\step[-1]\hcd\step[1.5]\ffbox3{\iy 11200}\hstep\ox{\sstyle 12}\\
\step[-1]\id\step\ox{\sstyle 11}\Step\id\step\id\Step\id\\
\hh\step[-1]\cu\Step\hstr{200}\cu\hstep\id\step\id\\
\hh\hhstep\id\step[2.5]\hstr{200}\hld\hstep\dd\step\id\\
\hh\hhstep\id\step[2.5]\id\step\x\step[3]\id\\
\hh\hhstep\id\step[2.5]\cu\step\hstr{300}\cu
\end{tangle}
\end{equation*}
where the first equation holds by definition and the second identity is an
application of \eqref{hat-prod}. In the third equation we use the cycle-cocycle
compatibility, \eqref{bra-ket-decomp}, the (co-)associativity of $\m_1$ and
$\Delta_2$, and \eqref{phi-prod}.

To prove the first identity of \eqref{gamma-12} we transform
in the next calculation the morphism $\gamma^1_t$ 
with the help of \eqref{tau-check-aux}, the entwining
properties \eqref{entw1} of $\varphi_{2,2}$, and \eqref{phi-prod} of
$\varphi_{1,2}$.
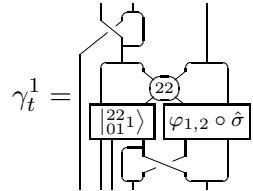
\begin{equation*}
\divide\unitlens by 3
\multiply\unitlens by 2
\gamma^1_t=
\vstretch 80
\begin{tangle}
\hh\step\id\step\cd\step[3]\id\\
\hh\step\x\step\id\step[3]\id\\
\hh\ne1\step{\hstr{200}\hru}\step[3]\id\\
\hh\id\step\hstr{200}\cd\step\cd\\
\id\step\id\step[2]\ox{\sstyle 22}\step[2]\id\\
\id\hstep\ffbox3{\yi 01221}\hstep\ffbox4{\varphi_{1,2}\circ\hat\sigma}\\
\hh\id\step\id\hstep\id\hstep{\hstr{200}\hld\step\id\step\id}\\
\hh\id\step\id\hstep\id\hstep\id\step\hstr{200}\x\step\id\\
\hh\id\step\id\hstep\id\hstep\cu\hstr{200}\step\cu
\end{tangle}
\end{equation*}
From the right hand side of this equation we get the right hand side
of \eqref{gamma-12} by applying \eqref{phi-hat-new} and again the 
entwining property \eqref{entw1} of $\varphi_{2,2}$.\end{proof}

\begin{lemma}
\begin{equation}
\unitlens = 8pt
\mathlabel{gamma-34}
\gamma^2_t=
\begin{tangle}
\hh\step\id\step[1.5]\id\step[1.5]\id\\
\hstr{50}\ffbox{11}{\tau^{3}_t}\\
\hh\step\id\hstep\id\hstep\id\hstep\id\step[1.5]\id\step\id\\
{\hstr{50}\ffbox{7}{\gamma^3_t}}\hstep\id\step\id\\
\hh\hstep\id\hstep\id\hstep\id\hstep\id\step\x\step\id\\
\hh\hstep\id\hstep\id\hstep\id\hstep\cu\step\cu
\end{tangle}
\quad,\quad
\gamma^2_b=
\begin{tangle}
\hh\cd\step\cd\hstep\id\hstep\id\hstep\id\\
\hh\id\step\x\step\id\hstep\id\hstep\id\hstep\id\\
\id\step\id\hstep\hstr{50}\ffbox{7}{\gamma^3_b}\\
\hh\id\step\id\step[1.5]\id\hstep\id\hstep\id\hstep\id\\
\hhstep\hstr{50}\ffbox{11}{\tau^{3}_t}\\
\hh\step\id\step[1.5]\id\step[1.5]\id
\end{tangle}
\end{equation}
and therefore $\beta^2_\ell = T^{(3)}(\beta_\ell^3)$.
\end{lemma}

\begin{proof}
The following identities hold.
\begin{equation}\mathlabel{gamma-derive}
\unitlens =7pt
\begin{tangle}
\hh\hstr{200}\cd\step\id\step\id\\
\id\Step\ox{\sstyle 22}\Step\id\\
\hhstep\ffbox3{\yi 01221}\step\ffbox3{\iy 21222}\\
\hh\id\hstep\id\hstep\hstr{200}\hld\hstep\dd\step\id\\
\hh\id\hstep\id\hstep\id\step\x\Step\sw1\\
\hh\id\hstep\id\hstep\cu\step\hstr{200}\cu\\
\object{\sstyle B}\step[1.5]\object{\sstyle B_2}\step[2.5]\object{\sstyle B_1}
\end{tangle}
\;=
\unitlens =12pt
\vstretch 80
\begin{tangle}
\id\step\td{\sstyle 02\vert 2}\step\id\\
\hh\x\Step\id\step\id\\
\id\hstep\ffbox3{\ixi 21212001}\hstep\id\\
\hh\id\step\id\Step\x\\
\hh\id\step\hstr{200}\cu\hstep\id\\
\object{\sstyle B}\step[2]\object{\sstyle B_2}\step[2]\object{\sstyle B_1}
\end{tangle}
\quad ,
\quad
\gamma^2_t=
\unitlens =12pt
\vstretch 80
\begin{tangle}
\vstr{60}\hh\step\id\Step\id\Step\id\\
\vstr{60}\hhstep\ffbox7{(\tau^3_t)^{\mathrm{red}_1}}\\
\hh\id\step\cd\step\id\Step\id\step\id\\
\Put(0,0)[lb]{\begin{tangle}\hh\x\step\id\\ 
\hh\id\step\hstr{200}\hru\end{tangle}}
    \Step\td{\sstyle 02\vert 2}\step\id\step\id\\
\hh\id\step\x\Step\id\step\id\step\id\\
\id\step\id\hstep\ffbox3{\ixi 21212001}\hstep\id\step\id\\
\hh\id\step\id\step\id\Step\x\step\id\\
\hh\id\step\id\step{\hstr{200}\cu}\step\cu\\
\object{\sstyle B_1}\step\object{\sstyle B}\step[2]\object{\sstyle B_2}%
  \step[2.5]\object{\sstyle B_1}
\end{tangle}
\end{equation}
The first identity of \eqref{gamma-derive} is a consequence of
\eqref{yi-red-adv5} and the module-comodule compatibility. Using 
this result we obtain the second relation in
\eqref{gamma-derive} with the help of \eqref{bra-ket-decomp} 
and the entwining property \eqref{entw1} of $\varphi_{2,2}$.

We denote
\begin{equation}
\tau^4_t:=
\unitlens =7pt
\vstretch 100
\begin{tangle}
\hh\cd\step\cd\step\cd\\
\hh\id\step\x\step\id\step\id\step\rd\\
\hh\id\step\id\step{\hstr{200}\hru}\step\id\step\id\hstep\hld\\
\id\step\id\step\ox{\sstyle 22}\step\id\hstep\id\hstep\id\\
\hh\id\step\id\step\id\step[2]{\vstr{50}\hx}\hstep\id\hstep\id\\
\hh\id\step\id\step\id\step[2]\id\step{\hstr{50}\vstr{50}\hx}\hstep\id\\
\hh\id\step\id\step\id\step[2]\id\step\id\hstep\hlu
\end{tangle}
\end{equation}
Then we obtain
\begin{equation}\mathlabel{gamma-derive2}
\begin{tangle}
\hh\id\step\cd\step\id\\
\Put(0,0)[lb]{\begin{tangle}
                 \hh\x\step\id\\ \hh\id\step\hstr{200}\hru
              \end{tangle}}
    \Step\td{\sstyle 02\vert 2}\\
\hh\id\step\x\Step\id\\
\id\step\id\hstep\ffbox3{\ixi 21212001}\\
\hh\id\step\id\step\id\Step\id\\
\object{\sstyle B_1}\step\object{\sstyle B}\step\object{\sstyle B_2}%
  \step[2]\object{\sstyle B_1}
\end{tangle}
=
\hstretch 80
\begin{tangle}
\hh\step\id\step[1.5]\id\step[1.5]\id\\
\hstr{40}\ffbox{11}{\tau^{4}_t}\\
\hh\step\id\hstep\id\hstep\id\hstep\id\step[1.5]\id\step\id\\
{\hstr{40}\ffbox{7}{\gamma^3_t}}\hstep\id\step\id\\
\hh\hstep\id\hstep\id\hstep\id\hstep\id\step\x\step\id\\
\hh\hstep\id\hstep\id\hstep\id\hstep\cu\step\cu
\end{tangle}
\quad,\quad
\tau^{3}_t=
\begin{tangle}
\hh\step[1.5]\id\step[1.5]\id\step[1.5]\id\\
\ffbox6{(\tau^3_t)^{\mathrm{red}_1}}\\
\hh\step\id\step\id\step\id\step[1.5]\id\step\id\\
\ffbox4{\tau^4_t}\hstep\id\step\id\\
\hh\hstep\id\hstep\id\hstep\id\hstep\id\hstep\id\step\x\step\id\\
\hh\hstep\id\hstep\id\hstep\id\hstep\id\hstep\cu\step\cu
\end{tangle}
\hstretch 100
\end{equation}
The first equation in \eqref{gamma-derive2} results from the 
algebra-coalgebra compatibility,
the right comodule-algebra compatibility, the left module-algebra
compatibility, \eqref{mod-alg-rel}, \eqref{rel-cocycle-assoc1}
and the right module algebra property of $B_2$.
The second identity can be obtained by subsequent application of
the right comodule-coalgebra compatibility, the left module-algebra
compatibility, \eqref{mod-alg-rel}, the right module-coalgebra property 
of $B_2$, and the entwining property of $\varphi_{2,2}$.

The first relation of \eqref{gamma-34} can now be derived by recursive 
substitution of the identities \eqref{gamma-derive2} into the second
identity of \eqref{gamma-derive}.\end{proof}
\abs
Eventually we define
\begin{equation}\mathlabel{beta-4}
\hstretch 100 \vstretch 100
\divide\unitlens by 2
\def\FillCircDiam{4}
\beta^4_\ell :=
\begin{tangle}
\hh\id\Step\id\step{\hstr{550}\cd}\step[3.25]{\hstr{125}\cd*}\Step\id\\
\hh{\hstr{200}\hrd}\step\x\step[2.5]{\hstr{600}\hld}\step[2.25]\sw1
    \step[1.25]\nw1\step\id\\
\hh\id\step\x\step{\hstr{250}\x}\step[3]\id\step[1.5]{\hstr{150}\cd}\step
    {\hstr{300}\hld}\step\id\\
\hh{\hstr{200}\hru}\step\x\step{\hstr{300}\cd}\step\cd\step\id\hstep
    {\hstr{200}\hld\hstep\hrd}\hstep\id\step\id\\
\hh\id\step\sw1\step\id\step\id\step[3]\x\step\id\step\id\hstep\id\step\x
    \step\id\hstep\id\step\id\\
\hh\id\step\id\Step\id\step\id\Step\sw1\step{\hstr{200}\hru}\step\id\hstep
    \cu\step\cu\hstep\x\\
\hh\id\step\id\Step\id\step\id\Step\id\Step{\hstr{200}\x}\step\id\Step\x
    \step\id\\
\hh\id\step\id\Step\id\step\id\Step{\hstr{200}\x}\Step\nw1{\hstr{200}\cu}
    \step\id\step\id\\
\hh\id\step\id\Step\id\step{\hstr{200}\x}\Step\nw1\Step\x\Step\id\step\id\\
\hh\id\step\id\Step\id\step\id\step{\hstr{200}\cd\hstep\cd}\step\id\step
    \id\Step\id\step\id\\
\hh\id\step\id\Step\id\step\id\step{\hstr{200}\hrd}\step\x\step
    {\hstr{200}\hld}\step\id\step\id\Step\id\step\id\\
\hh\id\step\id\Step\id\step\id\step\id\step\x\step\x\step\id\step\id
    \step\id\Step\id\step\id\\
\hh\id\step\id\Step\id\step\id\step{\hstr{200}\hru}\step\x\step
    {\hstr{200}\hlu}\step\id\step\id\Step\id\step\id\\
\hh\id\step\id\Step\id\step\id\step{\hstr{200}\cu\hstep\cu}\step\id\step
    \id\Step\id\step\id\\
\hh\id\step\id\Step\x\Step\nw1\Step{\hstr{200}\x}\step\id\Step\id\step\id\\
\hh\id\step\id\step{\hstr{200}\cd}\nw1\Step{\hstr{200}\x}\Step\id\step\id
    \Step\id\step\id\\
\hh\id\step\x\Step\id\step{\hstr{200}\x}\Step\id\Step\id\step\id\Step\id
    \step\id\\
\hh\x\hstep\cd\step\cd\hstep\id\step{\hstr{200}\hld}\step\sw1\Step\id\step
    \id\Step\id\step\id\\
\hh\id\step\id\hstep\id\step\x\step\id\hstep\id\step\id\step\x\step[3]\id
    \step\id\step\sw1\step\id\\
\hh\id\step\id\hstep{\hstr{200}\hlu\hstep\hru}\hstep\id\step\cu\step
    {\hstr{300}\cu}\step\x\step\hstr{200}\hld\\
\hh\id\step{\hstr{300}\hru}\step{\hstr{150}\cu\step}\id\step[3]
    {\hstr{250}\x}\step\x\step\id\\
\hh\id\step\nw1\step[1.25]\sw1\step[2.25]{\hstr{600}\hru}\step[2.5]\x
    \step{\hstr{200}\hlu}\\
\hh\id\Step{\hstr{125}\cu*}\step[3.25]{\hstr{550}\cu}\step\id\Step\id
\end{tangle}
\quad\text{and}\quad
\beta^4_r :=
\def\FillCircDiam{5}
\begin{tangle}
\hh\id\step[3.25]\id\step[5.25]{\hstr{300}\cd}\step[3.5]\id\step[3]\id\\
\hh\id\step[3.25]\id\step[4.75]\dd\Step{\hstr{200}\hld}\step[3]\cd*
    \step[2.5]\id\\
\hh\id\step[3.25]\id\step[4.25]\dd\step[1.5]\sw1\step{\hstr{300}\x}\step
    \d\Step\id\\
\hh\id\step[3.25]\id\step[3.75]\dd\step\sw1\step\sw1\step[3]\id
    \step[1.5]\d\step[1.5]\id\\
\hh\id\step[2.5]{\hstr{150}\cd}\step[2.5]\dd\step\dd\step\sw1\hstep
    {\hstr{350}\ld}\step[-1.5]{\hstr{300}\hld}\step{\hstr{200}\hld}
    \step[1.5]\id\\
\hh{\hstr{300}\hrd}\step\id\step[1.5]\id\Step\dd\step[1.5]\id
    \step\dd\step[1.5]\id\Step\id\step[1.5]\x\step\id\step[1.5]\id\\
\hh\id\hstep{\hstr{200}\hld\hstep\hrd}\hstep\id\step[1.5]
    \dd\Step\id\step\id\Step\id\Step\id\step\cd\hstep\id\step
    {\hstr{150}\x}\\
\hh\id\hstep\id\step\x\step\id\hstep\id\step\dd\step[2.5]\id
    \step\id\Step\id\Step\x\step\id\hstep\x\step[1.5]\id\\
\hh\id\hstep\cu\step\cu\hstep\x\step[2.5]\dd\step\id\step
    {\hstr{200}\cd\hstep\hrd}\cu\hstep\id\step\d\step\id\\
\hh\id\step\id\Step\x\step{\hstr{250}\x}\step\cd\hstep
    {\hstr{200}\hrd}\step\x\step\id\hstep\x\step[1.5]\id\step\id\\
\hh\id\step\id\step[1.5]\dd\step\x\step[2.5]\x\step\id\hstep\id
    \step\x\step{\hstr{200}\hlu}\hstep\id\step\id\step[1.5]\id\step\id\\
\hh\id\step\id\step[1.5]\id\step[1.5]\id\step{\hstr{250}\x}
    \step\x\hstep\cu\step{\hstr{200}\cu}\dd\step\id\step[1.5]\id
    \step\id\\
\hh\id\step\id\step[1.5]\id\step[1.5]\x\step[2.5]\x\step\x
    \step[2.5]\x\step[1.5]\id\step[1.5]\id\step\id\\
\hh\id\step\id\step[1.5]\id\step\dd{\hstr{200}\cd}\step\cd
    \hstep\x\step{\hstr{250}\x}\step\id\step[1.5]\id\step[1.5]\id
    \step\id\\
\hh\id\step\id\step[1.5]\id\step\id\hstep{\hstr{200}\hrd}\step
    \x\step\id\hstep\id\step\x\step[2.5]\x\step\dd\step[1.5]\id\step\id\\
\hh\id\step\id\step[1.5]\x\hstep\id\step\x\step{\hstr{200}\hlu}
    \hstep\cu\step{\hstr{250}\x}\step\x\Step\id\step\id\\
\hh\id\step\d\step\id\hstep\cd{\hstr{200}\hlu\hstep\cu}
    \step\id\step\dd\step[2.5]\x\hstep\cd\step\cd\hstep\id\\
\hh\id\step[1.5]\x\hstep\id\step\x\Step\id\Step\id\step\id
    \step[2.5]\dd\step\id\hstep\id\step\x\step\id\hstep\id\\
\hh{\hstr{150}\x}\step\id\hstep\cu\step\id\Step\id\Step\id
    \step\id\Step\dd\step[1.5]\id\hstep{\hstr{200}\hlu\hstep\hru}
    \hstep\id\\
\hh\id\step[1.5]\id\step\x\step[1.5]\id\Step\id\step[1.5]\dd\step\id
    \step[1.5]\dd\Step\id\step[1.5]\id\step{\hstr{300}\hlu}\\
\hh\id\step[1.5]{\hstr{200}\hru}\step{\hstr{300}\hru}\step[-1.5]
    {\hstr{350}\ru}\hstep\sw1\step\dd\step\dd\step[2.5]{\hstr{150}\cu}
    \step[2.5]\id\\
\hh\id\step[1.5]\d\step[1.5]\id\step[3]\sw1\step\sw1\step\dd
    \step[3.75]\id\step[3.25]\id\\
\hh\id\Step\d\step{\hstr{300}\x}\step\sw1\step[1.5]\dd\step[4.25]\id
    \step[3.25]\id\\
\hh\id\step[2.5]\cu*\step[3]{\hstr{200}\hru}\Step\dd\step[4.75]\id
    \step[3.25]\id\\
\hh\id\step[3]\id\step[3.5]{\hstr{300}\cu}\step[5.25]\id\step[3.25]\id\\
\end{tangle}
\end{equation}

In what follows we will tacitly use the various forms of (weak)
(co\n-)associativity of $\m_1$, $\m_2$, $\Delta_1$, $\Delta_2$,
$\mu_l$, $\mu_r$, $\nu_l$ and $\nu_r$ which are given by Lemma 
\ref{strong-hd3} in particular. The subsequent graphical calculation yields 
$\beta^3_\ell=T^{(4)}(\beta^4_\ell)$.
\begin{equation*}
\beta^3_\ell =\;
\divide\unitlens by 2
\hstretch 100 \vstretch 100
\def\FillCircDiam{5}
\begin{tangle}
\hh\step[4.5]\id\step{\hstr{400}\hld}\step[2.5]\id\step[1.5]\id\\
\hh\step[4.5]\x\Step\id\step[2.5]\id\step[1.5]\id\\
\hh\step[3.5]\sw1\hstep\cd\step\cd\step[1.5]\cd*\step\id\\
\hh\step[2.5]\sw1\hstep\sw1\step\x\hstep\dd\step\dd\hstep\hld\step
    {\hstr{200}\hrd}\\
\hh\step[1.5]\sw1\hstep\sw1\step[1.5]\dd\step\hru\hstep\sw1\hstep\dd\hstep
    \x\step\nw1\\
\hh\hstep\sw1\step\dd\Step\dd\step[1.5]\x\step[1.5]\cu\step\nw1\step\nw1\\
\hh\dd\step[1.5]\dd\step[2.5]\id\Step\id\step{\hstr{200}\x}\step[2.5]\d
    \step[1.5]\d\\
\hh\hstr{200}\id\hstep\cd\hstep\cd\hstep\id\hstep\id\hstep\cd\hstep\cd
    \hstep\id\\
\hh\id\step{\hstr{200}\hrd}\step\x\step{\hstr{200}\hld}\step\id\step
    \id\step{\hstr{200}\hrd}\step\x\step{\hstr{200}\hld}\step\id\\
\hh\id\step\id\step\x\step\x\step\id\step\id\step\id\step\id
    \step\x\step\x\step\id\step\id\\
\hh\id\step{\hstr{200}\hru}\step\x\step{\hstr{200}\hlu}\step\id\step
    \id\step{\hstr{200}\hru}\step\x\step{\hstr{200}\hlu}\step\id\\
\hh\hstr{200}\id\hstep\cu\hstep\cu\hstep\id\hstep\id\hstep\cu\hstep\cu
    \hstep\id\\
\hh\d\step[1.5]\d\step[2.5]{\hstr{200}\x}\step\id\Step\id\step[2.5]\dd
    \step[1.5]\dd\\
\hh\hstep\nw1\step\nw1\step\cd\step[1.5]\x\step[1.5]\dd\Step\dd\step
    \sw1\\
\hh\step[1.5]\nw1\step\x\hstep\dd\hstep\sw1\hstep\hld\step\dd\step[1.5]
    \sw1\hstep\sw1\\
\hh\step[2.5]{\hstr{200}\hlu}\step\hru\hstep\dd\step\dd\hstep\x\step
    \sw1\hstep\sw1\\
\hh\step[3.5]\id\step\cu*\step[1.5]\cu\step\cu\hstep\sw1\\
\hh\step[3.5]\id\step[1.5]\id\step[2.5]\id\Step\x\\
\hh\step[3.5]\id\step[1.5]\id\step[2.5]{\hstr{400}\hlu}\step\id
\end{tangle}
\multiply\unitlens by 2
\;=\;
\divide\unitlens by 3
\multiply\unitlens by 2
\hstretch 80 \vstretch 80
\def\FillCircDiam{3}
\begin{tangle}
\hh\id\Step{\hstr{600}\hld}\step[4.25]\id\Step\id\\
\hh{\hstr{160}\x}\step[3.75]\id\step[4.25]\id\Step\id\\
\hh\id\step{\hstr{160}\cd}\step{\hstr{280}\cd}\step[2.5]\id\Step\id\\
\hh\id\step{\hstr{160}\hrd}\step\x\step[1.5]{\hstr{320}\hld}
        \step[2.5]\id\Step\id\\
\hh\id\step\id\step\x\step{\hstr{120}\x}\Step\id\step[2.5]\id\Step\id\\
\hh\id\step{\hstr{160}\hru}\step\x\step\cd\step\cd\Step\id\Step\id\\
\hh\id\step{\hstr{160}\cu}\step\id\step\id\step\x\step\id\step
        {\hstr{160}\cd*}\step\id\\
\hh\id\Step\id\Step\id\step\id\step\id\step{\hstr{160}\hru}\sw1\Step\id
        \step\id\\
\hh\id\Step\id\Step\id\step\id\step\id\step\x\step[3]\id\step\id\\
\hh\id\Step\id\Step\id\step\id\step\x\step\nw1\step{\hstr{160}\hld}
        \step{\hstr{320}\hrd}\\
\hh\id\Step\id\Step\id\step\x\step\nw1\step\id\step\id\step\x\Step\id\\
\id\Step\id\Step\hx\step\ox{\sstyle 21}\step\id\step\hddcu\step\id\Step\id\\
\id\Step\id\step\hddcd\step\id\step\ox{\sstyle 12}\step\hx\Step\id\Step\id\\
\hh\id\Step\x\step\id\step\id\step\nw1\step\x\step\id\Step\id\Step\id\\
\hh{\hstr{320}\hlu}\step{\hstr{160}\hru}\step\nw1\step\x\step\id\step\id
        \Step\id\Step\id\\
\hh\Step\id\step\id\step[3]\x\step\id\step\id\step\id\Step\id\Step\id\\
\hh\Step\id\step\id\Step\sw1{\hstr{160}\hld}\step\id\step\id\step\id
        \Step\id\Step\id\\
\hh\Step\id\step{\hstr{160}\cu*}\step\id\step\x\step\id\step\id\step
        {\hstr{160}\cd}\step\id\\
\hh\Step\id\Step\id\Step\cu\step\cu\step\x\step{\hstr{160}\hld}\step\id\\
\hh\Step\id\Step\id\step[2.5]\id\Step{\hstr{120}\x}\step\x\step\id
        \step\id\\
\hh\Step\id\Step\id\step[2.5]{\hstr{320}\hru}\step[1.5]\x\step
        {\hstr{160}\hlu}\step\id\\
\hh\Step\id\Step\id\step[2.5]{\hstr{280}\cu}\step{\hstr{160}\cu}\step\id\\
\hh\Step\id\Step\id\step[4.25]\id\step[3.75]{\hstr{160}\x}\\
\hh\Step\id\Step\id\step[4.25]{\hstr{600}\hru}\Step\id
\end{tangle}
\;=\;
\def\FillCircDiam{2}
\begin{tangle}
\hh\id\Step{\hstr{680}\hld}\step[6]\id\step[2.75]\id\\
\hh{\hstr{160}\x}\step[4.25]\id\step[6]\id\step[2.75]\id\\
\hh\id\step{\hstr{160}\cd}\step{\hstr{360}\cd}\step[3.75]\id
        \step[2.75]\id\\
\hh\id\step{\hstr{160}\hrd}\step\x\step[2.5]{\hstr{320}\hld}\step[3.75]
        \id\step[2.75]\id\\
\hh\id\step\id\step\x\step{\hstr{200}\x}\Step\id\step[3.75]
        \id\step[2.75]\id\\
\hh\id\step{\hstr{160}\hru}\step\x\Step\cd\step\cd\step[1.5]
        {\hstr{280}\cd*}\step\id\\
\hh\id\step\id\step\sw1\step\id\Step\id\step\x\step\id\hstep
        {\hstr{160}\cd}\step{\hstr{240}\hld}\step{\hstr{160}\hrd}\\
\hh\id\step\id\step\id\Step\id\Step\id\step\id\step{\hstr{160}\hru}\dd
        \step{\hstr{160}\hld\hstep\hrd}\hstep\id\step\id\step\id\\
\hh\id\step\id\step\id\Step\id\Step\id\step\id\step\x\step[1.5]\id
        \step\x\step\id\hstep\id\step\id\step\id\\
\hh\id\step\id\step\id\Step\id\Step\id\step\x\step\nw1\hstep\cu\step
        \cu\hstep\x\step\id\\
\hh\id\step\id\step\id\Step\id\Step\x\step\nw1\step\id\step\id\Step
        \x\step\id\step\id\\
\id\step\id\step\id\Step\x\step\ox{\sstyle 21}\step\id\step\cu\step\id
        \step\id\step\id\\
\id\step\id\step\id\step\cd\step\id\step\ox{\sstyle 12}\step\x\Step\id\step\id
        \step\id\\
\hh\id\step\id\step\x\Step\id\step\id\step\nw1\step\x\Step\id\Step\id
        \step\id\step\id\\
\hh\id\step\x\hstep\cd\step\cd\hstep\nw1\step\x\step\id\Step\id\Step
        \id\step\id\step\id\\
\hh\id\step\id\step\id\hstep\id\step\x\step\id\step[1.5]\x\step\id
        \step\id\Step\id\Step\id\step\id\step\id\\
\hh\id\step\id\step\id\hstep{\hstr{160}\hlu\hstep\hru}\step\dd
        {\hstr{160}\hld}\step\id\step\id\Step\id\Step\id\step
        \id\step\id\\
\hh{\hstr{160}\hlu}\step{\hstr{240}\hru}\step{\hstr{160}\cu}\hstep\id
        \step\x\step\id\Step\id\step\sw1\step\id\step\id\\
\hh\step\id\step{\hstr{280}\cu*}\step[1.5]\cu\step\cu\Step\x\step
        {\hstr{160}\hld}\step\id\\
\hh\step\id\step[2.75]\id\step[3.75]\id\Step{\hstr{200}\x}\step\x
        \step\id\step\id\\
\hh\step\id\step[2.75]\id\step[3.75]{\hstr{320}\hru}\step[2.5]\x\step
        {\hstr{160}\hlu}\step\id\\
\hh\step\id\step[2.75]\id\step[3.75]{\hstr{360}\cu}\step{\hstr{160}\cu}
        \step\id\\
\hh\step\id\step[2.75]\id\step[6]\id\step[4.25]{\hstr{160}\x}\\
\hh\step\id\step[2.75]\id\step[6]{\hstr{680}\hru}\Step\id
\end{tangle}
\;=T^{(4)}(\beta^4_\ell)
\end{equation*}
where we used the module-comodule compatibility
and \eqref{mod-alg-rel} in the first equation. In the second 
identity we use again \eqref{mod-alg-rel}. To get the third equation in the
graphic we applied the right module-algebra and the left comodule-coalgebra
compatibilities of Definition \ref{hp} as well as \eqref{rel-cocycle-assoc2}.
With the help of the module-comodule compatibility the fourth identity
can be derived from the definition of $\beta^4_\ell$ in 
\eqref{beta-4} and the definition of $T^{(4)}$ in \eqref{dressing}.

Next we will prove the identity $\beta^3_r = 
T^{(4)}\left(\beta^4_r\right)$. And therefore all 
relations in both rows of the diagram in Proposition \ref{b-bialg} 
are satisfied.

\begin{equation*}
\beta^3_r=\;
\hstretch 60 \vstretch 60
\def\FillCircDiam{5}
\begin{tangle}
\hh{\hstr{120}\cd}\step{\hstr{270}\cd}\Step\id\step[1.5]\id\\
\hh{\hstr{120}\hrd}\step\x\Step{\hstr{300}\hld}\step[1.5]\cd*\step\id\\
\hh\id\step\x\step\id\Step\id\step[2.5]{\hstr{90}\x}\step\id\step\id\\
\hh\id\step\id\step\id\step\id\step{\hstr{120}\cd}\step\cd\step
        \cu\step\id\\
\hh\id\step\id\step\id\step\id\step{\hstr{120}\hrd}\step\x\step\id
        \step[1.5]\id\step[1.5]\id\\
\hh\id\step\id\step\id\step\id\step\id\step\x\step{\hstr{120}\hlu}
        \hstep\ld\step[1.5]\id\\
\hh\id\step\id\step\id\step\nw1{\hstr{120}\hlu}\step{\hstr{120}\cu}
        \hstep\id\step\d\step\id\\
\hh\id\step\id\step\nw1\step\x\step\sw1\step\cd\step\x\\
\hh\id\step\id\Step\x\step\x\step\sw1\step\x\step\id\\
\hh\id\step\id\step\sw1\step\x\step\x\step\sw1\step\id\step\id\\
\hh\id\step\x\step\sw1\step\x\step\x\Step\id\step\id\\
\hh\x\step\cu\step\sw1\step\x\step\nw1\step\id\step\id\\
\hh\id\step\d\step\id\hstep{\hstr{120}\cd\hstep\hrd\d}\step\id
        \step\id\step\id\\
\hh\id\step[1.5]\ru\hstep{\hstr{120}\hrd}\step\x\step\id
        \step\id\step\id\step\id\step\id\\
\hh\id\step[1.5]\id\step[1.5]\id\step\x\step{\hstr{120}\hlu}
        \step\id\step\id\step\id\step\id\\
\hh\id\step\cd\step\cu\step{\hstr{120}\cu}
        \step\id\step\id\step\id\step\id\\
\hh\id\step\id\step{\hstr{90}\x}\step[2.5]\id\Step\id\step\x\step\id\\
\hh\id\step\cu*\step[1.5]{\hstr{300}\hru}\Step\x\step\hstr{120}\hlu\\
\hh\id\step[1.5]\id\Step{\hstr{270}\cu}\step\hstr{120}\cu
\end{tangle}
\;=\;
\def\FillCircDiam{3}
\begin{tangle}
\hh\id\step[5]{\hstr{180}\hld}\step[5.5]\id\step[1.5]\id\\
\hh\id\step[4]\sw1\hstep{\hstr{120}\cd}\step[4.5]\id\step[1.5]\id\\
\hh\id\step[3]\sw1\step[1.5]\id\hstep{\hstr{180}\hld}\step[3.5]
        {\hstr{120}\cd*}\hstep\id\\
\hh\nw1\step\sw1\step[2.5]\id\hstep\id\step[1.5]{\hstr{210}\x}\Step\id
        \hstep\id\\
\hh\step\x\step[3.5]\id\hstep\id\step\cd\step[1.5]{\hstr{180}\hld}
        \step{\hstr{120}\hld}\hstep\hrd\\
\hh\sw1{\hstr{120}\cd}\Step\dd\hstep\x\step\id\step\dd\step[1.5]
        \x\step\id\hstep\id\hstep\id\\
\hh\id\step{\hstr{120}\hrd}\step{\hstr{120}\x}\step\id\step\lu\hstep
        \cd\step\cd\hstep\cu\hstep\id\hstep\id\\
\hh\id\hstep\dd\step\x\Step\x\Step\id\hstep\id\step\x\step\id\step\x
        \hstep\id\\
\hh\id\hstep\id\step[1.5]\id\step{\hstr{120}\x\hstep\x}\hstep\cu\step\cu
        \hstep\dd\step\id\hstep\id\\
\hh\id\hstep\id\step[1.5]\x\Step\x\Step\x\Step\x\step[1.5]\id\hstep\id\\
\hh\id\hstep\id\step\dd\hstep\cd\step\cd\hstep{\hstr{120}\x\hstep\x}
        \step\id\step[1.5]\id\hstep\id\\
\hh\id\hstep\x\step\id\step\x\step\id\hstep\id\Step\x\Step\x\step\dd
        \hstep\id\\
\hh\id\hstep\id\hstep\cd\hstep\cu\step\cu\hstep\rd\step\id\step
        {\hstr{120}\x\hstep\hlu}\step\id\\
\hh\id\hstep\id\hstep\id\step\x\step[1.5]\dd\step\id\step\x\hstep\dd
        \Step{\hstr{120}\cu}\sw1\\
\hh\hlu\hstep{\hstr{120}\hru}\step{\hstr{180}\hru\hstep}\cu\step\id
        \hstep\id\step[3.5]\x\\
\hh\hstep\id\hstep\id\Step{\hstr{210}\x}\step[1.5]\id\hstep\id
        \step[2.5]\sw1\step\nw1\\
\hh\hstep\id\hstep{\hstr{120}\cu*}\step[3.5]{\hstr{180}\hru}\hstep\id
        \step[1.5]\sw1\step[3]\id\\
\hh\hstep\id\step[1.5]\id\step[4.5]{\hstr{120}\cu}\hstep\sw1\step[4]
        \id\\
\hh\hstep\id\step[1.5]\id\step[5.5]{\hstr{180}\hru}\step[5]\id
\end{tangle}
\;=\;
\hstretch 100 \vstretch 100
\divide\unitlens by 2
\def\FillCircDiam{5}
\begin{tangle}
\hh\step\id\step[6]{\hstr{800}\hld}\step[5]\id\step[3]\id\\
\hh\step\id\step[4]\sw2\step[3.5]{\hstr{300}\cd}\step[3.5]\id\step[3]\id\\
\hh\step\id\Step\sw2\step[5]\dd\Step{\hstr{200}\hld}\step[3]\cd*
    \step[2.5]\id\\
\hh\step{\hstr{200}\x}\step[5.5]\dd\step[1.5]\sw1\step{\hstr{300}\x}\step
    \d\Step\id\\
\hh\sw1\step{\hstr{200}\cd}\step[4]\dd\step\sw1\step\sw1\step[3]\id
    \step[1.5]\d\step[1.5]\id\\
\hh\id\step{\hstr{200}\cd}\step\nw1\step[2.5]\dd\step\dd\step\sw1\hstep
    {\hstr{350}\ld}\step[-1.5]{\hstr{300}\hld}\step{\hstr{200}\hld}
    \step[1.5]\id\\
\hh\id\step{\hstr{300}\hrd}\hstep\d\step[1.5]\id\Step\dd\step[1.5]\id
    \step\dd\step[1.5]\id\Step\id\step[1.5]\x\step\id\step[1.5]\rd\\
\hh\id\step\id\hstep{\hstr{200}\hld\hstep\hrd}\hstep\id\step[1.5]
    \dd\Step\id\step\id\Step\id\Step\id\step\cd\hstep\id\step
    {\hstr{150}\x}\step\id\\
\hh\id\step\id\hstep\id\step\x\step\id\hstep\id\step\dd\step[2.5]\id
    \step\id\Step\id\Step\x\step\id\hstep\x\step[1.5]\id\step\id\\
\hh\id\step\id\hstep\cu\step\cu\hstep\x\step[2.5]\dd\step\id\step
    {\hstr{200}\cd\hstep\hrd}\cu\hstep\id\step\d\step\id\step\id\\
\hh\id\step\id\step\id\Step\x\step{\hstr{250}\x}\step\cd\hstep
    {\hstr{200}\hrd}\step\x\step\id\hstep\x\step[1.5]\id\step\id\step\id\\
\hh\id\step\id\step\id\step[1.5]\dd\step\x\step[2.5]\x\step\id\hstep\id
    \step\x\step{\hstr{200}\hlu}\hstep\id\step\id\step[1.5]\id\step\id
    \step\id\\
\hh\id\step\id\step\id\step[1.5]\id\step[1.5]\id\step{\hstr{250}\x}
    \step\x\hstep\cu\step{\hstr{200}\cu}\dd\step\id\step[1.5]\id
    \step\id\step\id\\
\hh\id\step\id\step\id\step[1.5]\id\step[1.5]\x\step[2.5]\x\step\x
    \step[2.5]\x\step[1.5]\id\step[1.5]\id\step\id\step\id\\
\hh\id\step\id\step\id\step[1.5]\id\step\dd{\hstr{200}\cd}\step\cd
    \hstep\x\step{\hstr{250}\x}\step\id\step[1.5]\id\step[1.5]\id
    \step\id\step\id\\
\hh\id\step\id\step\id\step[1.5]\id\step\id\hstep{\hstr{200}\hrd}\step
    \x\step\id\hstep\id\step\x\step[2.5]\x\step\dd\step[1.5]\id\step\id
    \step\id\\
\hh\id\step\id\step\id\step[1.5]\x\hstep\id\step\x\step{\hstr{200}\hlu}
    \hstep\cu\step{\hstr{250}\x}\step\x\Step\id\step\id\step\id\\
\hh\id\step\id\step\d\step\id\hstep\cd{\hstr{200}\hlu\hstep\cu}
    \step\id\step\dd\step[2.5]\x\hstep\cd\step\cd\hstep\id\step\id\\
\hh\id\step\id\step[1.5]\x\hstep\id\step\x\Step\id\Step\id\step\id
    \step[2.5]\dd\step\id\hstep\id\step\x\step\id\hstep\id\step\id\\
\hh\id\step{\hstr{150}\x}\step\id\hstep\cu\step\id\Step\id\Step\id
    \step\id\Step\dd\step[1.5]\id\hstep{\hstr{200}\hlu\hstep\hru}
    \hstep\id\step\id\\
\hh\lu\step[1.5]\id\step\x\step[1.5]\id\Step\id\step[1.5]\dd\step\id
    \step[1.5]\dd\Step\id\step[1.5]\d\hstep{\hstr{300}\hlu}\step\id\\
\hh\step\id\step[1.5]{\hstr{200}\hru}\step{\hstr{300}\hru}\step[-1.5]
    {\hstr{350}\ru}\hstep\sw1\step\dd\step\dd\step[2.5]\nw1\step
    {\hstr{200}\cu}\step\id\\
\hh\step\id\step[1.5]\d\step[1.5]\id\step[3]\sw1\step\sw1\step\dd
    {\hstr{200}\Step\cu}\step\sw1\\
\hh\step\id\Step\d\step{\hstr{300}\x}\step\sw1\step[1.5]\dd\step[5.5]
    {\hstr{200}\x}\\
\hh\step\id\step[2.5]\cu*\step[3]{\hstr{200}\hru}\Step\dd\step[4]\sw2
    \step[3]\id\\
\hh\step\id\step[3]\id\step[3.5]{\hstr{300}\cu}\step[2.5]\sw2\step[5]\id\\
\hh\step\id\step[3]\id\step[5]{\hstr{800}\hru}\step[6]\id
\end{tangle}
\;=T^{(4)}(\beta^4_r)
\end{equation*}
The first identity in this graphical calculation has been obtained from
\eqref{rel-mod-coalg} and \eqref{rel-comod-alg} of Lemma \ref{strong-hd8}.
In the second equation we use \eqref{mod-alg-rel2}, Lemma
\ref{strong-hd5} and \eqref{alg-coalg-rel}.
With the help of the left module-algebra and right comodule-coalgebra 
compatibilities of Definition~\ref{hp} and the relativizations
\eqref{rel-mod-coalg} and \eqref{rel-comod-alg} of Lemma~\ref{strong-hd8}
we derive the third identity. The fourth equation holds by definition.

In order to complete the proof of the proposition the equation
$\beta^4_\ell=\beta^4_r$ needs to be shown. This will be done in the
following calculation.
\begin{equation*}
\hstretch 100 \vstretch 100
\divide\unitlens by 2
\def\FillCircDiam{4}
\beta^4_\ell=
\begin{tangle}
\hh\id\Step\id\step{\hstr{550}\cd}\step[3.25]{\hstr{125}\cd*}\Step\id\\
\hh{\hstr{200}\hrd}\step\x\step[2.5]{\hstr{600}\hld}\step[2.25]\sw1
    \step[1.25]\nw1\step\id\\
\hh\id\step\x\step{\hstr{250}\x}\step[3]\id\step[1.5]{\hstr{150}\cd}\step
    {\hstr{300}\hld}\step\id\\
\hh{\hstr{200}\hru}\step\x\step{\hstr{300}\cd}\step\cd\step\id\step[1.5]
    \id\step{\hstr{200}\hrd}\hstep\id\step\id\\
\hh\id\step\sw1\step\id\step\id\step[3]\x\step\id\step\id\step[1.5]\x
    \step\id\hstep\id\step\id\\
\hh\id\step\id\Step\id\step\id\Step\sw1\step{\hstr{200}\hru}\step\id\step
    \dd
\step\cu\hstep\x\\
\hh\id\step\id\Step\id\step\id\Step\id\Step{\hstr{200}\x}\step\id\Step\x
    \step\id\\
\hh\id\step\id\Step\id\step\id\Step{\hstr{200}\x}\Step\nw1{\hstr{200}\cu}
    \step\id\step\id\\
\hh\id\step\id\Step\id\step{\hstr{200}\x}\Step\nw1\Step\x\Step\id\step\id\\
\hh\id\step\id\Step\id\step\id\step{\hstr{200}\cd\hstep\cd}\step\id\step
    \id\Step\id\step\id\\
\hh\id\step\id\Step\id\step\id\step{\hstr{200}\hrd}\step\x\step
    {\hstr{200}\hld}\step\id\step\id\Step\id\step\id\\
\hh\id\step\id\Step\id\step\id\step\id\step\x\step\x\step\id\step\id
    \step\id\Step\id\step\id\\
\hh\id\step\id\Step\id\step\id\step{\hstr{200}\hru}\step\x\step
    {\hstr{200}\hlu}\step\id\step\id\Step\id\step\id\\
\hh\id\step\id\Step\id\step\id\step{\hstr{200}\cu\hstep\cu}\step\id\step
    \id\Step\id\step\id\\
\hh\id\step\id\Step\x\Step\nw1\Step{\hstr{200}\x}\step\id\Step\id\step\id\\
\hh\id\step\id\step{\hstr{200}\cd}\nw1\Step{\hstr{200}\x}\Step\id\step\id
    \Step\id\step\id\\
\hh\id\step\x\Step\id\step{\hstr{200}\x}\Step\id\Step\id\step\id\Step\id
    \step\id\\
\hh\x\hstep\cd\step\dd\step\id\step{\hstr{200}\hld}\step\sw1\Step\id\step
    \id\Step\id\step\id\\
\hh\id\step\id\hstep\id\step\x\step[1.5]\id\step\id\step\x\step[3]\id
    \step\id\step\sw1\step\id\\
\hh\id\step\id\hstep{\hstr{200}\hlu}\step\id\step[1.5]\id\step\cu\step
    {\hstr{300}\cu}\step\x\step\hstr{200}\hld\\
\hh\id\step{\hstr{300}\hru}\step{\hstr{150}\cu\step}\id\step[3]
    {\hstr{250}\x}\step\x\step\id\\
\hh\id\step\nw1\step[1.25]\sw1\step[2.25]{\hstr{600}\hru}\step[2.5]\x
    \step{\hstr{200}\hlu}\\
\hh\id\Step{\hstr{125}\cu*}\step[3.25]{\hstr{550}\cu}\step\id\Step\id
\end{tangle}
=
\def\FillCircDiam{2}
\begin{tangle}
\hh\id\step[4]\id\step[3.5]{\hstr{700}\cd}\step[5]\id\step[3.5]\id\\
\hh\id\step[4]\id\step[3.5]\id\step[1.5]{\hstr{1100}\hld}\step[5]\id
    \step[3.5]\id\\
\hh\id\step[4]\id\step[3.5]\id\step[1.5]\id\step[5]\cd\Step{\hstr{500}\cd*}
    \step\id\\
\hh\id\step[3]{\hstr{200}\cd}\Step\cd\step{\hstr{200}\hrd}\step
    {\hstr{600}\hld}\step\id\step[1.5]\cd\step[1.5]{\hstr{600}\hld}
    \step\id\\
\hh\Rd\step\id\Step\d\step[1.5]\id\step\x\step\x\step{\hstr{400}\hld}
    \step\id\step\dd\hstep\cd\step{\hstr{200}\hrd}\Step\id\step\id\\
\hh\id\step{\hstr{200}\hld}\step\rd\step\cd\step\id\step\id\step\x\step\id
    \step\id\step{\hstr{200}\hld}\step\id\step\id\step\id\step\x\step\id
    \step\ld\step\id\\
\hh\id\step\id\step\x\step\id\step\id\step\x\step\id\step\id\step\id
    \step\id\step\id\step\id\step\id\step\x\step\id\step\id\step\id\step
    \id\step\id\step\x\\
\hh\id\step\id\step\id\step\id\step\id\step\x\step\x\step\id\step\id
    \step\id\step\id\step\id\step\x\step\x\step\id\step\id\step\id\step
    \x\step\id\\
\hh\id\step\id\step\id\step\id\step\x\step\x\step\x\step\id\step\id
    \step\id\step\x\step\x\step\x\step\id\step\x\step\id\step\id\\
\hh\id\step\id\step\id\step\x\step\x\step\x\step\x\step\id\step\x\step
    \x\step\x\step\x\step\id\step\id\step\id\step\id\\
\hh\id\step\id\step\x\step\x\step\x\step\x\step\x\step\id\step\x\step
    \x\step\x\step\x\step\id\step\id\step\id\\
\hh\id\step\x\step\x\step\x\step\x\step\x\step\x\step\id\step\x\step
    \id\step\id\step\x\step\x\step\id\step\id\\
\hh\id\step\id\step\id\step\id\step\id\step\id\step\x\step\x\step\x\step
    \x\step\id\step\id\step\id\step\x\step\id\step\id\step\id\step\id\step\id\\
\hh\id\step\id\step\id\step\id\step\id\step\id\step\id\step\x\step
    \x\step\x\step\x\step\id\step\x\step\id\step\id\step\id\step\id
    \step\id\step\id\\
\hh\id\step\id\step\id\step\id\step\id\step\id\step\id\step\id\step\x
    \step\x\step\id\step\id\step\id\step\x\step\id\step\id\step\id
    \step\id\step\id\step\id\step\id\\
\hh\id\step\id\step\id\step\id\step\id\step\id\step\id\step\x\step\id
    \step\id\step\id\step\x\step\x\step\id\step\id\step\id\step\id
    \step\id\step\id\step\id\step\id\\
\hh\id\step\id\step\id\step\id\step\id\step\id\step\x\step\id\step
    \x\step\x\step\x\step\x\step\id\step\id\step\id\step\id\step\id
    \step\id\step\id\\
\hh\id\step\id\step\id\step\id\step\id\step\x\step\id\step\id\step\id\step
    \x\step\x\step\x\step\x\step\id\step\id\step\id\step\id\step\id
    \step\id\\
\hh\id\step\id\step\x\step\x\step\id\step\id\step\x\step\id\step\x\step
    \x\step\x\step\x\step\x\step\x\step\id\\
\hh\id\step\id\step\id\step\x\step\x\step\x\step\x\step\id\step\x\step\x
    \step\x\step\x\step\x\step\id\step\id\\
\hh\id\step\id\step\id\step\id\step\x\step\x\step\x\step\x\step\id\step\x
    \step\x\step\x\step\x\step\id\step\id\step\id\\
\hh\id\step\id\step\x\step\id\step\x\step\x\step\x\step\id\step\id\step
    \id\step\x\step\x\step\x\step\id\step\id\step\id\step\id\\
\hh\id\step\x\step\id\step\id\step\id\step\x\step\x\step\id\step\id
    \step\id\step\id\step\id\step\x\step\x\step\id\step\id\step\id\step
    \id\step\id\\
\hh\x\step\id\step\id\step\id\step\id\step\id\step\x\step\id\step\id\step
    \id\step\id\step\id\step\id\step\id\step\x\step\id\step\id\step\x
    \step\id\step\id\\
\hh\id\step\ru\step\id\step\x\step\id\step\id\step\id\step{\hstr{200}\hru}
    \step\id\step\id\step\x\step\id\step\id\step\cu\step\lu\step
    {\hstr{200}\hru}\step\id\\
\hh\id\step\id\Step{\hstr{200}\hlu}\step\cu\hstep\dd\step\id\step
    {\hstr{400}\hru}\step\x\step\x\step\id\step[1.5]\d\Step\id\step\Lu\\
\hh\id\step{\hstr{600}\hru}\step[1.5]\cu\step[1.5]\id\step{\hstr{600}\hru}
    \step{\hstr{200}\hlu}\step\cu\Step{\hstr{200}\cu}\step[3]\id\\
\hh\id\step{\hstr{500}\cu*}\Step\cu\step[5]\id\step[1.5]\id\step[3.5]\id
    \step[4]\id\\
\hh\id\step[3.5]\id\step[5]{\hstr{1100}\hru}\step[1.5]\id\step[3.5]\id
    \step[4]\id\\
\hh\id\step[3.5]\id\step[5]{\hstr{700}\cu}\step[3.5]\id\step[4]\id
\end{tangle}
=
\def\FillCircDiam{5}
\begin{tangle}
\hh\id\step[3.25]\id\step[5.25]{\hstr{300}\cd}\step[3.5]\id\step[3]\id\\
\hh\id\step[3.25]\id\step[4.75]\dd\Step{\hstr{200}\hld}\step[3]\cd*
    \step[2.5]\id\\
\hh\id\step[3.25]\id\step[4.25]\dd\step[1.5]\sw1\step{\hstr{300}\x}\step
    \d\Step\id\\
\hh\id\step[3.25]\id\step[3.75]\dd\step\sw1\step\sw1\step[3]\id
    \step[1.5]\d\step[1.5]\id\\
\hh\id\step[2.5]{\hstr{150}\cd}\step[2.5]\dd\step\dd\step\sw1\hstep
    {\hstr{350}\ld}\step[-1.5]{\hstr{300}\hld}\step{\hstr{200}\hld}
    \step[1.5]\id\\
\hh{\hstr{300}\hrd}\step\id\step[1.5]\id\Step\dd\step[1.5]\id
    \step\dd\step[1.5]\id\Step\id\step[1.5]\x\step\id\step[1.5]\id\\
\hh\id\hstep{\hstr{200}\hld\hstep\hrd}\hstep\id\step[1.5]
    \dd\Step\id\step\id\Step\id\Step\id\step\cd\hstep\id\step
    {\hstr{150}\x}\\
\hh\id\hstep\id\step\x\step\id\hstep\id\step\dd\step[2.5]\id
    \step\id\Step\id\Step\x\step\id\hstep\x\step[1.5]\id\\
\hh\id\hstep\cu\step\cu\hstep\x\step[2.5]\dd\step\id\step[1.5]\dd
    \step\sw1\step\cu\hstep\id\step\d\step\id\\
\hh\id\step\id\Step\x\step{\hstr{250}\x}\step\cd\hstep\cd\step
    {\hstr{200}\hrd}\step[1.5]\x\step[1.5]\id\step\id\\
\hh\id\step\id\step[1.5]\dd\step\x\step[2.5]\x\step\id\hstep\id
    \step\x\step\id\step[1.5]\id\step\id\step[1.5]\id\step\id\\
\hh\id\step\id\step[1.5]\id\step[1.5]\id\step{\hstr{250}\x}
    \step\x\hstep\cu\step{\hstr{200}\hlu}\step\dd\step\id\step[1.5]\id
    \step\id\\
\hh\id\step\id\step[1.5]\id\step[1.5]\x\step[2.5]\x\step\x
    \step[2.5]\x\step[1.5]\id\step[1.5]\id\step\id\\
\hh\id\step\id\step[1.5]\id\step\dd\step{\hstr{200}\hrd}\step\cd
    \hstep\x\step{\hstr{250}\x}\step\id\step[1.5]\id\step[1.5]\id
    \step\id\\
\hh\id\step\id\step[1.5]\id\step\id\step[1.5]\id\step
    \x\step\id\hstep\id\step\x\step[2.5]\x\step\dd\step[1.5]\id\step\id\\
\hh\id\step\id\step[1.5]\x\step[1.5]{\hstr{200}\hlu}
    \step\cu\hstep\cu\step{\hstr{250}\x}\step\x\Step\id\step\id\\
\hh\id\step\d\step\id\hstep\cd\step\sw1\step\dd\step[1.5]\id\step\dd
    \step[2.5]\x\hstep\cd\step\cd\hstep\id\\
\hh\id\step[1.5]\x\hstep\id\step\x\Step\id\Step\id\step\id
    \step[2.5]\dd\step\id\hstep\id\step\x\step\id\hstep\id\\
\hh{\hstr{150}\x}\step\id\hstep\cu\step\id\Step\id\Step\id
    \step\id\Step\dd\step[1.5]\id\hstep{\hstr{200}\hlu\hstep\hru}
    \hstep\id\\
\hh\id\step[1.5]\id\step\x\step[1.5]\id\Step\id\step[1.5]\dd\step\id
    \step[1.5]\dd\Step\id\step[1.5]\id\step{\hstr{300}\hlu}\\
\hh\id\step[1.5]{\hstr{200}\hru}\step{\hstr{300}\hru}\step[-1.5]
    {\hstr{350}\ru}\hstep\sw1\step\dd\step\dd\step[2.5]{\hstr{150}\cu}
    \step[2.5]\id\\
\hh\id\step[1.5]\d\step[1.5]\id\step[3]\sw1\step\sw1\step\dd
    \step[3.75]\id\step[3.25]\id\\
\hh\id\Step\d\step{\hstr{300}\x}\step\sw1\step[1.5]\dd\step[4.25]\id
    \step[3.25]\id\\
\hh\id\step[2.5]\cu*\step[3]{\hstr{200}\hru}\Step\dd\step[4.75]\id
    \step[3.25]\id\\
\hh\id\step[3]\id\step[3.5]{\hstr{300}\cu}\step[5.25]\id\step[3.25]\id\\
\end{tangle}
=\beta^4_r
\end{equation*}
where the first equation has been derived with the help of
\eqref{act-coact-triv}. In the second and the third equation we use 
\eqref{mod-alg-rel} and \eqref{mod-alg-rel2}.
Finally we apply \eqref{strong-hopf-combo} to obtain the fourth 
identity. 

Hence Proposition \ref{b-bialg} and therefore Theorem \ref{hp-cp} have been
proven.\end{proof}

\abs
\abs
\abs
\parbox{12cm}{\scriptsize{\textsc{Bogolyubov Institute for Theoretical
Physics, Metrologichna Str.~14-b, Kiev 252143, Ukraine.}
\textbf{E-mail:} {\tt mmtpitp@gluk.apc.org}
\\[0.5cm]
\textsc{DAMTP, University of Cambridge, Silver Street, Cambridge CB3 9EW,
UK.} \textbf{E-mail:} {\tt b.drabant@damtp.cam.ac.uk}}}


\begin{thebibliography}{MM}

\bibitem{And1:96}
N.~Andruskiewitsch and J.~Devoto, 
  \emph{Extensions of Hopf Algebras}, St.~Petersburg Math.~J.,
  \textbf{7}, No.~1, 17-52 (1996).

\bibitem{Bes1:97}
Yu.~N.~Bespalov, \emph{Crossed modules and quantum groups in braided
  categories}, Applied Categorical Structures \textbf{5}, no.~2,
  155-204 (1997).
\nl
------, \emph{Crossed modules, quantum braided groups and ribbon
 structures}, Teoreticheskaya i Matematicheskaya Fizika \textbf{103}, no.~3,
 368-387 (1995).

\bibitem{BD:95}
Yu.~N.~Bespalov and B.~Drabant, \emph{{H}opf ({B}i-){M}odules and {C}rossed
  {M}odules in {B}raided {C}ategories}, J.~Pure Appl.~Algebra
 \textbf{123}, 105-129 (1998).

\bibitem{BD:97}
Yu.~N.~Bespalov and B.~Drabant, \emph{{D}ifferential {C}alculus in
{B}raided {T}ensor {C}ategories}, preprint q-alg/9703036 (1997).
\nl
------, \emph{Bicovariant Differential Calculi and Cross
Products on  Braided Hopf Algebras}, in Proceedings ``Quantum Groups
and Quantum Spaces, 1995'', Banach Center Publications \textbf{40},
Institute of Math., Polish Acad.~Sci., eds.~R.~Budzynski et
al (Warsaw 1997).

\bibitem{BD1:98}
Yu.~N.~Bespalov and B.~Drabant, \emph{{C}ross {P}roduct {B}ialgebras
-- {P}art {I}}, J.~Algebra \textbf{219}, 466-505 (1999).

\bibitem{BD3:98}
Yu.~N.~Bespalov and B.~Drabant, \emph{{C}ross {P}roduct {B}ialgebras
-- {P}art {III}}, Preprint DAMTP-98-134 (1998).

\bibitem{BCM:86}
R.~J.~Blattner, M.~Cohen and S.~Montgomery,
\emph{Crossed products and inner actions of Hopf algebras},
Trans.~AMS \textbf{298}, 671-711 (1986).

\bibitem{Brz1:96}
T.~Brzezinski, \emph{Crossed products by a coalgebra}, Commun.\ Algebra 
 \textbf{25} (11), 3551-3575 (1997).

\bibitem{CIMZ:98}
S.~Caenepeel, B.~Ion, G.~Militaru and S.~Zhu,
\emph{Factorisation Structures of Algebras and Coalgebras}, Preprint (1998).
To appear in Algebras and Representation Theory.

\bibitem{CM:98}
A.~Connes and H.~Moscovici, \emph{Hopf Algebras, Cyclic Cohomology and
the Transverse Index Theorem}, Commun.~Math.~Phys.~\textbf{198},
199-246 (1998).

\bibitem{CZ:98}
H.-X.~Chen and S.~Zhang, \emph{The Double Biproduct in Braided Tensor
Categories}, Preprint (1998).

\bibitem{DT:86}
Y.~Doi and M.~Takeuchi, \emph{Cleft comodule algebras
for a bialgebra}, Commun.~Alg.~\textbf{14}, 801-818 (1986).

\bibitem{FY:89}
  P.~{Freyd} and D.~{Yetter},
  \emph{Braided compact closed categories with applications to low
  dimensional topology}, Adv.~Math.~\textbf{77}, 156-182 (1989).
\nl
------, \emph{Coherence theorems via knot theory}, J.~Pure Appl.~Algebra
  \textbf{78}, 49-76 (1992).

\bibitem{Hof1:90}
  I.~Hofstetter, \emph{Erweiterungen von {H}opf-{A}lgebren und ihre
  kohomologische {B}eschreibung}, Ph.D.~thesis, Universit\"at M\"unchen (1990).
\nl
  ------, \emph{Extensions of Hopf Algebras and Their Cohomological
  Description}, J.~Algebra \textbf{164}, 264-298 (1994).

\bibitem{JS:93}
A.~{Joyal} and R.~{Street},
\emph{{B}raided tensor categories},
Adv.~Math.~\textbf{102}, 20-78 (1993).
\nl
------, \emph{{B}raided tensor categories},
        Macquarie Math.~Reports 860081 (1986).

\bibitem{Kas1:95}
C.~Kassel,
\emph{{Q}uantum {G}roups}, volume \textbf{155} of {GTM},
Springer (1995).

\bibitem{Lyu1:95}
V.~V.~Lyubashenko,
\emph{Tangles and {H}opf algebras in braided categories}
J.~Pure Appl.~Algebra, \textbf{98}, 245 (1995).
\nl
------, \emph{{M}odular transformations for tensor categories}
J.~Pure Appl.~Algebra, \textbf{98}, 279 (1995).

\bibitem{Ma1:90}
S.~Majid, \emph{Physics for algebraists: {N}on-commutative and
non-cocommutative {H}opf algebras by a bicrossproduct construction},
J.~Algebra \textbf{130}, no.~1, 17-64 (1990).

\bibitem{Ma4:93}
S.~{Majid}, \emph{{B}raided {G}roups},
J.~Pure~Appl.~Algebra \textbf{86}, 187 (1993).

\bibitem{Ma6:94}
S.~Majid, \emph{More examples of bicrossproduct and double cross
product Hopf algebras}, Isr.~J.~Math.~\textbf{72}, 133-148 (1990).

\bibitem{Ma1:95}
S.~Majid, \emph{Double-bosonisation of braided groups and the
construction of $u_q(g)$}, to appear in Math.~Proc.~Camb.~Phil.~Soc.,
preprint RIMS-1047, q-alg/9511001 (1995).

\bibitem{MS:94}
 S.~Majid and Ya.~S.~Soibelman,
 \emph{{B}icrossproduct structure of the quantum {W}eyl group},
 J.~Algebra \textbf{163}, 68-87 (1994).

\bibitem{Mon1:92}
S.~Montgomery, \emph{{H}opf {A}lgebras and {T}heir {A}ctions
on {R}ings}, CBMS \textbf{82}, AMS (1992).

\bibitem{Rad1:85}
D.~E.~ Radford,
\emph{{T}he {S}tructure of {H}opf {A}lgebras with a {P}rojection},
J.~Algebra \textbf{92}, 322 (1985).

\bibitem{ReT:90}
  N.~Yu.~Reshetikhin and V.~G.~Turaev,
  \emph{{R}ibbon {G}raphs and {T}heir {I}nvariants {D}erived from {Q}uantum
  {G}roups}, Commun.~Math.~Phys.~\textbf{127}, 1 (1990).

\bibitem{Sbg1:98}
P.~Schauenburg.
\emph{{T}he {S}tructure of {H}opf {A}lgebras with a {W}eak {P}rojection},
Preprint (1998).

\bibitem{Sch1:94}
H.-J.~Schneider, \emph{{H}opf {G}alois {E}xtensions, {C}rossed
{P}roducts, and {C}lifford {T}heory}, Advances in Hopf Algebras.
Lecture Noptes in Pure Appl.~Math.~{\bf 158}, 267-297, (J.~Bergen and
S.~Montgomery, ed.), Dekker (1994).

\bibitem{Sin1:72}
W.~M.~{Singer}, \emph{Extension Theory for Connected Hopf Algebras},
J.~Algebra \textbf{21}, 1-16 (1972).

\bibitem{Swe1:68}
M.~E.~{Sweedler}, \emph{Cohomology of algebras over Hopf algebras},
Trans.\ Amer.\ Math.\ Soc.\ \textbf{133}, 205-239 (1968).

\bibitem{Tak1:81}
M.~Takeuchi, \emph{Matched pairs of groups and bismash products of {H}opf
  algebras}, Comm.\ Alg.\ \textbf{9}, 841--882 (1981).

\bibitem{Tur1:94}
V.~G.~Turaev,
\emph{{Q}uantum {I}nvariants of {K}nots and 3-{M}anifolds},
volume~\textbf{18} of {S}tudies in {M}athematics,
Walter de Gruyter (Berlin, New York 1994).
\end{thebibliography}
\end{document}